\documentclass[10pt]{article}

\usepackage{amssymb,amsmath,latexsym}
\voffset
-1.5cm

\usepackage{a4,amscd,graphicx,epsfig}
\usepackage[all]{xy}
\everymath{\displaystyle}

\newtheorem{teo}{Theorem}[section]
\newtheorem{lem}[teo]{Lemma}

\newtheorem{exa}[teo]{Example}
\newtheorem{prop}[teo]{Proposition}
\newtheorem{defi}[teo]{Definition}

\newtheorem{conj}[teo]{Conjecture}
\newtheorem{remark}[teo]{Remark}
\newtheorem{remarks}[teo]{Remarks}

\newcommand{\mr}{\mathbb{R}}
\newcommand{\mc}{\mathbb{C}}
\newcommand{\mz}{\mathbb{Z}}
\newcommand{\mh}{\mathbb{H}}

\newcommand{\mq}{\mathbb{Q}}

\newcommand{\Dd}{{\mathcal D}}

\newcommand{\Ii}{{\mathcal I}}

\newcommand{\Ll}{{\mathcal L}}

\newcommand{\Tt}{{\mathcal T}}

\newcommand{\C}{{\mathbb C}}
\newcommand{\Q}{{\mathbb Q}}
\newcommand{\Z}{{\mathbb Z}}
\newcommand{\R}{{\mathbb R}}

\title{Quantum Hyperbolic Invariants \\
Of $3$-Manifolds With $PSL(2,\C)$-Characters}

\author {St\'ephane Baseilhac$^\dagger$$^\ddagger$ and Riccardo
Benedetti$^\ddagger$}

\date {}

\begin{document}

\maketitle

\vspace{0.7cm}

\noindent $^\dagger$ 
\ Laboratoire E. Picard, UMR CNRS 5580, UFR MIG, Universit\'e Paul Sabatier, 118 route de Narbonne, 
F-31062 TOULOUSE.

\medskip

\noindent $^\ddagger$ \ Dipartimento di Matematica, Universit\`a di Pisa, Via
F. Buonarroti, 2, I-56127 PISA. Emails: baseilha@mail.dm.unipi.it, benedett@dm.unipi.it. 

\medskip

\vspace{1cm}
\begin{abstract}  

\noindent We construct {\it quantum hyperbolic invariants} (QHI) for
triples $(W,L,\rho)$, where $W$ is a compact closed oriented
$3$-manifold, $\rho$ is a flat principal bundle over $W$ with
structural group $PSL(2,\mc)$, and $L$ is a non-empty link in $W$.
These invariants are based on the Faddeev-Kashaev's {\it quantum
dilogarithms}, interpreted as matrix valued functions of suitably
decorated hyperbolic ideal tetrahedra. They are explicitely computed as
state sums over the decorated hyperbolic ideal tetrahedra of the {\it
idealization} of any fixed {\it $\Dd$-triangulation}; the $\Dd$-triangulations are simplicial $1$-cocycle descriptions of
$(W,\rho)$ in which the link is realized as a Hamiltonian subcomplex. We also
discuss how to set the Volume Conjecture for the coloured
Jones invariants $J_N(L)$ of hyperbolic knots $L$ in $S^3$ in the
framework of the general QHI theory.  

\bigskip


\medskip

\end{abstract}

\section{Introduction}\label{intro}
In the two papers \cite{K1,K5}, Kashaev proposed a new infinite family
$\{K_N\}$, $N > 1$ being any odd positive integer, of conjectural
complex valued topological invariants for pairs $(W,L)$, where $L$ is
a link in a compact closed oriented $3$-manifold $W$. These invariants
should be computed as {\it state sums} $K_N(\Tt)$ supported by some
kind of heavily decorated triangulation $\Tt$ for $(W,L)$. The main
ingredients of the state sums were the Faddeev-Kashaev's matrix
version of the {\it quantum dilogarithms} at the $N$-th root of unity
$\zeta= \exp (2\pi i/N)$, suitably associated to the decorated
tetrahedra of $\Tt$.  The nature of these decorated triangulations was
mysterious, but it was clear that they fulfilled non trivial global
constraints which made their existence not evident a priori.  Beside
this neglected {\it existence} and {\it meaning} problem, a main
question left unsettled was the {\it invariance} of the value of
$K_N(\Tt)$ when $\Tt$ varies. However, Kashaev proved the invariance
of $K_N(\Tt)$ under certain `moves' on $\Tt$, which reflect
fundamental identities verified by the quantum dilogarithms. He also showed in \cite{K2} that $K_N(S^3,L)$ is indeed a well-defined invariant, by reducing its state sum formula to one based on planar $(1,1)$-tangle presentations of $L$ (as for the Alexander polynomial) and involving a constant {\it Kashaev's $R$-matrix}.

\smallskip

\noindent On another
hand, Faddeev-Kashaev \cite{FK}, Bazhanov-Reshetikhin \cite{BR} and
Kashaev \cite{K3} had already computed the semi-classical limit of
(various versions of) the quantum dilogarithms and their five term
`pentagon' identities in terms of the classical Euler and Rogers
dilogarithm functions, which are known to be related to the
computation of the volume of spherical or hyperbolic simplices. This
definetely suggested the possibility of a deep intriguing relationship
between hyperbolic geometry and the invariants $K_N(S^3,L)$. In this
direction, the so-called Kashaev's {\it Volume Conjecture}
\cite{K4} predicts that when $L$ is a hyperbolic link in 
$S^3$, one can recover the hyperbolic volume of $S^3 \setminus L$ from
the asymptotic behaviour of $K_N(S^3,L)$, when $N\to \infty$. More
recently, Murakami-Murakami \cite{MM} proved that the Kashaev's
$R$-matrix can be enhanced into a Yang-Baxter operator which allows
one to define the coloured Jones polynomial $J_N(L)$ for links $L$ in
$S^3$ (evaluated at $\zeta=\exp(2i\pi/N)$ and normalized by
$J_N(\rm{unknot})=1$), so that $K_N(S^3,L)=J_N^N(L)$. This gave a new formulation of the Volume Conjecture, discussed in \cite{MM,Y}, in terms of those celebrated invariants of links.
 
\smallskip

\noindent 
The new formulation of $J_N(L)$ using quantum dilogarithms was an
important achievement, but it also had the negative consequence of
putting aside the initial purely $3$-dimensional and more geometric
set-up for links in an arbitrary compact closed oriented $3$-manifold
$W$, willingly forgetting the complicated and somewhat mysterious
decorated triangulations.

\noindent
In our opinion, this set-up deserved to be understood and developed as
a full Quantum Field Theory, also in the perspective of finding an
appropriate geometric framework for a well motivated more general
version of the Volume Conjecture. The present paper, which is the
first of a series, establishes some fundamental facts of our program
on this matter. The main result is the construction of so-called {\it
quantum hyperbolic invariants} (QHI) for compact closed oriented
$3$-manifolds endowed with an embedded non-empty link and a flat
principal bundle with structural group $PSL(2,\mc)$. The QHI generalize
the Kashaev's conjectural topological invariants.

\medskip

\noindent {\bf Description of the paper.}  We are mainly concerned
with pairs $(W,\rho)$ where $W$ is a compact closed oriented
$3$-manifold and $\rho$ is a flat principal bundle over $W$ with
structural group $PSL(2,\C)$. By using the hauptvermutung, depending
on the context, we will freely assume that $W$ is endowed with a
(necessarily unique) PL or smooth structure, and use differentiable or
PL homeomorphisms. The pairs $(W,\rho)$ are considered up to orientation preserving homeomorphisms of $W$ and flat bundle isomorphisms of
$\rho$. Equivalently, $\rho$ is identified with a conjugacy class of
representations of the fundamental group of $W$ in $PSL(2,\C)$,
i.e. with a $PSL(2,\C)$-character of $W$. Compact oriented hyperbolic
$3$-manifolds with their hyperbolic holonomies furnish a main example
of pairs $(W,\rho)$. There are other natural examples
$(W,\rho_\alpha)$ associated to ordinary cohomology classes $\alpha
\in H^1(W; \C)$ (see Subsection \ref{flatbungen}).
\medskip

\noindent 
In Section \ref{distingW} we introduce special combinatorial
descriptions of $(W,\rho)$ called $\Dd$-{\it triangulations}. These
are ``decorated'' triangulations $T$ of $W$, where the decoration
consists of a system $b$ of edge orientations of a special kind
(called {\it branching}), and of a `generic' $PSL(2,\mc)$-valued
$1$-cocycle $z$ on $(T,b)$. This genericity condition allows us to
define a simple explicit procedure of {\it idealization} which
converts any $\Dd$-triangulation $\Tt$ into a suitably structured
family of oriented hyperbolic ideal tetrahedra $\Tt_{\Ii}$, called an
$\Ii$-{\it triangulation} for $(W,\rho)$. Each hyperbolic tetrahedron
of $\Tt_{\Ii}$ has the vertices ordered by the branching $b$, and its
geometry is encoded by the cross-ratio moduli in $\mc \setminus
\{0,1\}$ associated to its edges. The $\Ii$-triangulations have
remarkable global properties. In particular their moduli satisfy, at
every edge, the usual compatibility condition needed when one tries to
construct hyperbolic $3$-manifolds by gluing ideal tetrahedra. This
means that given an $\Ii$-triangulation we can construct pairs
$(\tilde{\rho},s)$, where $\tilde{\rho}$ is a representative of the
character $\rho$ and $s$ is a piecewise-straight section of the flat
bundle $\widetilde{W} \times_{\tilde{\rho}} \bar{\mh}^3 \rightarrow
W$, with structural group $PSL(2,\mc)$ and total space the quotient of
$\widetilde{W} \times \bar{\mh}^3$ by the diagonal action of
$\pi_1(W)$ and $\tilde{\rho}$.

\noindent 
We also define the notions of $\Dd$ and $\Ii$-{\it transits} between
$\Dd$ and $\Ii$-triangulations of $(W,\rho)$. These are supported by
the usual elementary moves on triangulations of $3$-manifolds, but
they also include the transits of the respective extra-structures. We
prove the remarkable fact that, via the idealization, the
$\Dd$-transits dominate the $\Ii$-transits.
\medskip

\noindent 
In Section \ref{QDILOGSS}, we consider for any odd positive integer
$N>1$ certain {\it basic state sums} $\mathfrak{L}_N(\Tt_{\Ii}) \in
\mc$ supported by the idealization $\Tt_{\Ii}$ of any
$\Dd$-triangulation $\Tt$ for $(W,\rho)$. The main ingredients of
these state sums are the Faddeev-Kashaev (non symmetric) matrix
quantum dilogarithms, viewed as matrix valued functions depending on
the moduli of branched hyperbolic ideal tetrahedra. At this point some
comments are in order.

\noindent 
These matrix quantum dilogarithms (quantum dilogarithms for short)
were originally derived in \cite{K1,K5} as matrices of $6j$-{\it
symbols} for the cyclic representation theory of a Borel quantum
subalgebra $\mathcal{B}_{\zeta}$ of $U_{\zeta}(sl(2,\mc))$, where
$\zeta=\exp(2i\pi/N)$. Such matrices describe the associativity of the
tensor product in this category. Here are two key facts. First, the
isomorphism classes of irreducible cyclic representations of
$\mathcal{B}_{\zeta}$ are parametrized by the elements with non-zero
upper diagonal term in the Borel subgroup $B$ of $PSL(2,\C)$ of upper
triangular matrices. Moreover, the specific `Clebsch-Gordan'
decomposition rule into irreducibles of cyclic tensor products of such
representations relies on a (generic) $B$-valued $1$-cocycle-like
property.  These facts may be seen at hand, or alternatively they can
be deduced from the theory of quantum coadjoint action of De
Concini-Kac-Procesi
\cite{DCP}, applied to the group $B$. We recall them in the Appendix
of this paper (Section \ref{APP}), as well as the properties of the
quantum dilogarithms that we need; for full details we refer to
\cite{B}.

\noindent 
Thus, when associating irreducible cyclic representations of
$\mathcal{B}_{\zeta}$ to the edges of a branched tetrahedron
$(\Delta,b)$, generic $B$-valued $1$-cocycles on $\Delta$ seem to play
a fundamental role to associate quantum dilogarithms to it. For this
reason, we early considered the QHI only for $B$-valued characters of
$W$ (see \cite{BB1}). However, we eventually realized that the quantum
dilogarithms do in fact only depend on particular ratios of parameters
expressed in terms of the cocycle values, which may naturally be
interpreted as moduli for idealized tetrahedra. Also, the basic
identities they satisfy are only related to certain $\Ii$-transits. As
the idealization works for arbitrary $PSL(2,\C)$-characters on $W$,
this and the symmetrization procedure explained below finally leads to
the present general formulation of the theory. The quantum
dilogarithms do not appear in this way as directly related to the
whole cyclic representation theory of $U_{\zeta}(sl(2,\mc))$. Of
course it would be most useful to compute/compare explicitely the
matrices of $6j$-symbols for this theory. We expect that the theory of
quantum coadjoint action leads to generalizations of the QHI for other
semisimple Lie groups than $PSL(2,\mc)$.\footnote{The referee informed
the authors that Kashaev and Reshetikhin recently constructed new
invariants for complements of tangles in $S^3$ based on this theory,
after a preliminary version of the present paper was put on the web in
january 2001. See Kashaev and Reshetikhin, \emph{Invariants of tangles
with flat connections in their complements, I: Invariants and holonomy
R-matrices}, \emph{II: Holonomy R-matrices related to quantized
envelopping algebras at roots of $1$}, arXiv:math.AT/0202211.}

\noindent 
The value of the basic state sums $\mathfrak{L}_N(\Tt_{\Ii})$ is not
invariant with respect to the change of branching, and it is invariant
only for some specific instance of $\Ii$-transit. So, in order to
construct invariants for $(W,\rho)$ based on the quantum dilogarithms,
these should be modified in such a way that the corresponding modified
state sums are (at least) branching invariant and invariant with
respect to {\it all} instances of $\Ii$-transits. We do this via a
specific procedure of (partial) symmetrization of the quantum
dilogarithms.

\medskip

\noindent 
In Section \ref{LFIXQHI} we show that this local symmetrization leads
to fix an arbitrary non-empty link $L$ in $W$, considered up to
ambient isotopy, in order to fix one coherent globalization. So we
incorporate this {\it link-fixing} in all the discussion: we consider
triples $(W,L,\rho)$ up to orientation preserving homeomorphisms of $(W,L)$
and flat bundle isomorphisms of $\rho$, and we provide the appropriate notion of
$\Dd$-triangulation for a triple $(W,L,\rho)$. This is a
$\Dd$-triangulation $(T,b,z)$ for $(W,\rho)$ in which the link $L$ is
realized as a Hamiltonian subcomplex $H$ (i.e. $H$ contains all the
vertices of $T$). We also refine the $\Dd$-transits to
preserve this Hamiltonian property of $H$.

\noindent The globalization of the symmetrization of the quantum
dilogarithms is governed, for all odd positive integer $N>1$, by any fixed {\it
integral charge} $c$ on $(T,H)$. An integral charge is a $\Z$-valued
function of the edges of the (abstract) tetrahedra of $T$ that
satisfies suitable non-trivial global conditions, and which eventually
encodes $H$, hence the link $L$. In fact, for any fixed $N$, we rather
use the reduction mod($N$) of `half' the charge, i.e. $c'(e)=(p+1)\
c(e)$ mod($N$). This is a main point where it is important that $N$ is
odd.

\noindent 
The integral charges are a subtle ingredient of our
construction. Their structure is very close to the one of the
``flattenings'' used by Neumann in his work on Cheeger-Chern-Simons
classes of hyperbolic manifolds \cite{N1}-\cite{N3}. The main results
concerning the existence and the structure of the integral charges are
obtained by adapting some fundamental results of Neumann.

\noindent 
All this gives the notion of {\it charged} $\Dd$-triangulation
$(\Tt,c)=(T,H,b,z,c)$ for a triple $(W,L,\rho)$; we stress that their
existence is not an evident fact. The $\Dd$- and $\Ii$-transits are
extended to transits of charged triangulations. This is the final
set-up for defining the {\it quantum hyperbolic invariants} (QHI): the
idealization $(\Tt_{\Ii},c)$ of any charged $\Dd$-triangulation
supports modified state sums $H_N(\Tt_{\Ii},c) \in \mc$ based on the
symmetrized quantum dilogarithms. Up to a sign and an $N$-th root of unity
multiplicative factor, $H_N(\Tt_{\Ii},c)$ is invariant with respect to
the choice of branching and for all instances of charged
$\Ii$-transits.

\noindent 
In Subsection \ref{QHISS} we state the two main results of the present
paper, proved in Subsections \ref{DTWLEX} and \ref{QHIINV}
respectively: the existence of charged $\Dd$-triangulations for any
triple $(W,L,\rho)$, and the fact that the value of the state sums
$H_N(\Tt_{\Ii},c)$ does not depend on the choice of $(\Tt_{\Ii},c)$ up
to sign and $N$-th root of unity factors. This proof of invariance consists
in reducing the full invariance to the transit invariance mentioned
above. We eventually get the QHI $H_N(W,L,\rho)$, and
$K_N(W,L,\rho)=H_N(W,L,\rho)^{2N}$ is a well-defined complex valued
invariant for every odd integer $N>1$.

\noindent 
In Subsection \ref{qhicomp} we discuss some complements about the
QHI. In particular, we prove a {\it duality} property related to the
change of the orientation of $W$.

\medskip

\noindent 
We had presented in \cite{BB1} the construction of QHI for flat
$B$-bundles on $W$, where $B$ is the Borel subgroup of $PSL(2,\mc)$ of
upper triangular matrices. In that case we adopted a slightly
different symmetrization procedure. The resulting state sums differ
from $H_N(\Tt_{\Ii},c)$, which work for arbitrary $PSL(2,\C)$-bundles,
by a scalar factor depending on the charged $\Dd$-triangulation
$(\Tt,c)$, not only on its idealization $(\Tt_{\Ii},c)$ (see Remark
\ref{Bsym}). The topological invariants $K_N(W,L)$ conjectured by
Kashaev correspond to the particular case of these $B$-QHI, when
$\rho$ is the {\it trivial} flat bundle.

\medskip

\noindent  
In Section \ref{VOLCONJ} we discuss how to set the Volume Conjecture
for the Jones invariants $J_N(L)$ of hyperbolic links $L$ in $S^3$ in
the framework of the general QHI theory.

\medskip

\noindent 
An appropriate conceptual framework for both the QHI and the
dilogarithmic invariant defined in \cite{BB2} stems from the theory of {\it scissors
congruence classes} (see \cite{D2}, \cite{N2} and the references therein for
details on this theory). It is elaborated in
\cite{BB1} for flat $B$-bundles, and in general in \cite{BB2}.

\medskip

\noindent 
Let us conclude by saying that another idea on the background of our
work, that is at least a meaningful heuristic support, is to look at
it as part of an ``exact solution'' of the Euclidean analytic
continuation of $(2+1)$ quantum gravity with negative cosmological
constant, that was outlined in \cite{W}.  This should be a gauge
theory with gauge group $SL(2,\C)$ and an action of Chern-Simons
type. Hyperbolic $3$-manifolds are the empty ``classical solutions''.
The Volume Conjecture discussed in Section \ref{VOLCONJ} essentially
agrees with the expected ``semi-classical limits'' of the partition
functions of this theory (see page 77 of \cite{W}).

\section{$\Dd$-triangulations for a pair $(W,\rho)$} 
\label{distingW}

\noindent We first recall few generalities before defining the $\Dd$-triangulations.

\subsection{Generalities on triangulations and spines}\label{trigen}

\noindent For the fundations of this theory, including the existence of spines, the reconstruction of manifolds from them and the complete calculus of
triangulation/spine-moves, we refer to
\cite{Cas}, \cite{Mat}, \cite{Pi}. Other references are \cite{BP2}, \cite{BP5}. One finds also a clear treatment of this
material in \cite{TV} (note that sometimes the terminologies do not
agree). We shall refer to the topological space underlying a cell complex 
as its \emph{polyhedron}.

\medskip

\noindent Consider a tetrahedron $\Delta$ with its usual triangulation with $4$
vertices, and let $C$ be the interior of the $2$-skeleton of the dual
cell-decomposition. A {\it simple} polyhedron $P$ is a $2$-dimensional
compact polyhedron such that each point of $P$ has a neighbourhood
which can be embedded into an open subset of $C$. A simple polyhedron
$P$ has a natural stratification given by its singularities; $P$ is {\it
standard} (in \cite{TV} one uses the term {\it cellular}) if all the
strata of this stratification are open cells of the appropriate
dimension $\leq 2$. Depending on the dimension, we call the strata of
a standard polyhedron $P$ {\it vertices, edges} and {\it regions}.
\smallskip

\noindent Every compact $3$-manifold $Y$ (which for simplicity we assume
connected) with non-empty boundary has {\it standard
spines} \cite{Cas}, that is standard polyhedra $P$ together with an embedding
in Int($Y$) such that $Y$ is a regular neighbourhood
of $P$. Moreover, $Y$ can be reconstructed from any of its standard spines.
The standard polyhedra underlying standard spines of oriented $3$-manifolds are characterized by the property of carrying a suitable ``screw-orientation'' along the edges \cite{BP5}; a compact oriented $3$-manifold $Y$ can be reconstructed from any of its 
oriented standard spines. From now on we assume that $Y$ is oriented, and we
shall only consider oriented standard spines of it. Since we shall always work with combinatorial data encoded by triangulations/spines, which define the corresponding manifold only up to PL-homeomorphisms, we shall systematically forget the underlying embeddings.

\smallskip

\noindent A {\it singular} triangulation of a polyhedron $Q$ is a triangulation
in a weak sense, namely self-adjacencies and multiple adjacencies
of $3$-simplices along $2$-faces are allowed. For any $Y$ as above, 
let us denote by $Q(Y)$ the space obtained by collapsing each connected component of $\partial Y$ to a point. A (topological) {\it ideal triangulation} of $Y$ is a singular triangulation $T$ of $Q(Y)$ such that the vertices of $T$ are
precisely the points of $Q(Y)$ corresponding to the components of
$\partial Y$.

\noindent For any ideal triangulation $T$ of $Y$, the $2$-skeleton of the
\emph{dual} cell-decomposition of $Q(Y)$ is a standard spine $P(T)$ of
$Y$. This procedure can be reversed, so that we can associate to each
standard spine $P$ of $Y$ an ideal triangulation $T(P)$ of
$Y$ such that $P(T(P))=P$. Thus standard spines and ideal
triangulations are dual equivalent viewpoints which we will freely
intermingle. By removing small open neigbourhoods of the
vertices of $Q(Y)$, any ideal triangulation leads to a
cell decomposition of $Y$ by {\it truncated tetrahedra}, which restricts
to a singular triangulation of $\partial Y$.

\noindent Any ideal triangulation $T$ of $Y$ can be considered as a
finite family $\{\Delta_i\}$ of {\it oriented abstract} tetrahedra,
each being endowed with the standard triangulation with $4$ vertices
and the orientation induced by the one of $Y$, together with
identifications of pairs of distinct (abstract) $2$-faces. We will
often distinguish between edges and $2$-faces {\it in $T$}, that is
considered after the identifications, and {\it abstract} edges and
$2$-faces, that is considered as simplices of the abstract
$\Delta_i$'s. We view each $\Delta_i$ as positively embedded as a
straight tetrahedron in $\R^3$ endowed with the orientation specified
by the standard basis (the `right-hand screw rule').

\smallskip

\noindent Consider now a compact closed oriented $3$-manifold $W$. 
For any $r_0\geq 1$, let $Y=W_{r_0} = W\setminus r_0D^3$ be the
manifold obtained by removing $r_0$ disjoint open balls from $W$. By
definition $Q(Y)=W$ and any ideal triangulation of $Y$ is a singular
triangulation of $W$ with $r_0$ vertices; moreover, it is easily seen
that all singular triangulations of $W$ come in this way from ideal
triangulations. We shall adopt the following terminology. A singular
triangulation of $W$ is simply called a {\it triangulation}. Ordinary
triangulations (where neither self-adjacencies nor multi-adjacencies
are allowed) are said to be {\it regular}.

\begin{figure}[ht]
\begin{center}
\includegraphics[width=7cm]{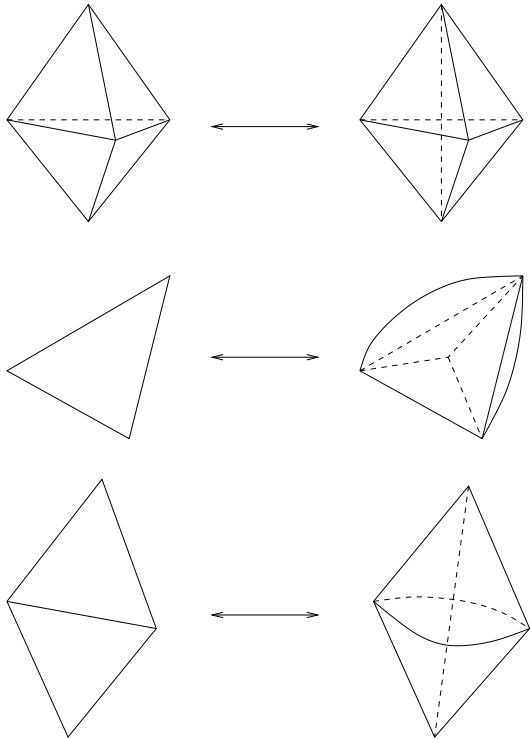}
\caption{\label{figmove1} the moves between singular triangulations.} 
\end{center}
\end{figure}

\begin{figure}[ht]
\begin{center}
\includegraphics[width=8cm]{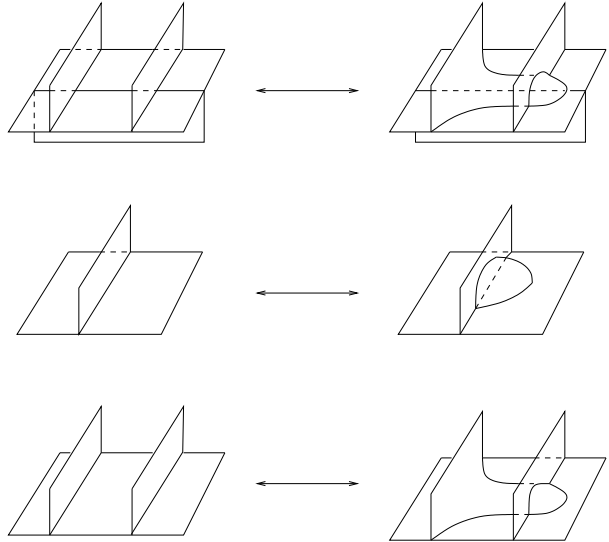}
\caption{\label{figmove2} the moves on standard spines.} 
\end{center}
\end{figure}

\noindent 
The main advantage in using singular triangulations (resp. standard
spines) instead of regular triangulations consists in the fact
that there exists a {\it finite} set of moves which are sufficient in
order to connect, by means of finite sequences of these moves,
any two singular triangulations (resp. standard spines) of the same
manifold. Let us recall some elementary moves on the triangulations
(resp. simple spines) of a polyhedron $Q(Y)$ that we shall use
throughout the paper; see Fig. \ref{figmove1} - Fig. \ref{figmove2}.

\medskip

\noindent {\bf The $2\to 3$ move.} 
Replace the triangulation $T$ of a portion of $Q(Y)$ made by the union
of $2$ tetrahedra with a common $2$-face $f$ by the triangulation made
by $3$ tetrahedra with a new common edge which connect the two
vertices opposite to $f$.  Dually this move is obtained by sliding
a portion of some region of $P(T)$ along an edge $e$, until it bumps into another region.

\smallskip

\noindent {\bf The bubble move.} 
Replace a face of a triangulation $T$ of $Q(Y)$ by the union of two
tetrahedra glued along three faces. Dually this move is done
gluing a closed $2$-disk $D$ via its boundary $\partial D$ on
the standard spine $P(T)$, with exactly two transverse intersection
points of $\partial D$ along some edge of $P(T)$. The new
triangulation thus obtained is dual to a spine of $Y \setminus D^3$,
where $D^3$ is an open ball in the interior of $Y$.

\smallskip

\noindent {\bf The $0 \to 2$ move.} 
Replace two adjacent faces of a triangulation $T$ of $Q(Y)$ by the
union of two tetrahedra glued along two faces, so that the other faces
match the two former ones. The dual of this move is the same as for
the $2 \to 3$ move, except that now we slide portions of regions {\it away}
from the edges of $P(T)$.

\noindent Standard spines of the same compact oriented $3$-manifold $Y$ with
boundary and with at least two vertices (which, of course, is a
painless requirement) may always be connected by means of a finite sequence of the (dual) $2 \to 3$ move and its inverse. In order to handle triangulations of closed
oriented $3$-manifolds we also need a move which allows us to vary the
number of vertices. The simplest way is to use the bubble move.
Note that a bubble move followed by a $2\to 3$ move with an adjacent tetrahedron gives a $1\to 4$ move: this simply consists in subdividing a 
tetrahedron $\Delta$ by the cone over its $2$-skeleton, with centre at an interior 
point of $\Delta$.  

\smallskip

\noindent Although the $2 \to 3$ move and the bubble move generate a complete
calculus for triangulations and standard spines, it is useful to introduce 
the $0\to 2$ move, or {\it lune} move. The inverse of the lune move is not always admissible because one could lose the standardness property of spines when using it. We say that a move
which increases (resp. decreases) the number of tetrahedra is \emph{positive} (resp. \emph{negative}). In some situations it may be useful to use only
positive moves. For that we need the following technical result due to Makovetskii \cite{Mak}:

\begin{prop} \label{chemin} 
Let $P$ and $P'$ be standard spines of $Y$. There exists a spine $P''$
of $Y$ such that $P''$ can be obtained from both $P$ and $P'$ via
finite sequences of positive $0 \to 2$ and $2 \to 3$ moves.
\end{prop}

\noindent In this paper we shall use a restricted class of triangulations.

\begin{defi} \label{quasireg}
{\rm A {\it quasi-regular} triangulation $T$ of a compact closed
$3$-manifold $W$ is a triangulation where all edges have
distinct vertices. A move $T \to T'$ is {\it quasi-regular} if
both $T$ and $T'$ are quasi-regular.}
\end{defi}

\begin{figure}[ht]
\begin{center}
\includegraphics[width=7cm]{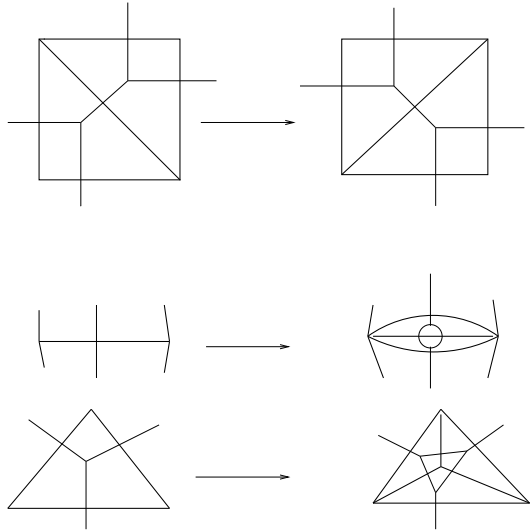}
\caption{\label{2dimmove} $2$-dimensional moves} 
\end{center}
\end{figure}

\noindent Of course any regular triangulation of $W$ is quasi-regular. We will also need the $2$-dimensional version of the above
facts. Given a compact closed surface $S$, there is a natural notion
of ideal triangulation $T$ of $S_{r_0} = S\setminus r_0D^2$ (for
arbitrary $r_0$) which corresponds to the notion of (singular)
triangulation of $S$ with $r_0$ vertices. The $1$-skeleton $P$ of the
dual cell decomposition of $T$ has only trivalent vertices and is a
standard spine of $S_{r_0}$.  In Fig. \ref{2dimmove} we show
$2$-dimensional moves on triangulations and their dual standard
spines: the $2\to 2$ ``flip'' move, which is the $2$-dimensional
analogue of the $2\to 3$ move, the $2$-dimensional bubble move, and
the $1\to 3$ move, which is the $2$-dimensional analogue of the above
$1\to 4$ move. Similarly to the $3$-dimensional case, the
$1\to 3$ move is a composition of a bubble move and a $2\to 2$ move. It
is known that any two arbitrary triangulations of $S$ with the same number
of vertices can be connected by a finite sequence of $2\to 2$ moves;
hence, to connect arbitrary triangulations of $S$ we only need a
further move which increases by one the number of vertices. Finally, we still have the notion of quasi-regular triangulations of a surface $S$.

\noindent Let $W$ be a compact closed oriented $3$-manifold, $T$ be a quasi-regular triangulation of $W$, and $v_0$ 
be a vertex of $T$. The link $S =$ Link($v_0,T$) with its natural 
triangulation $T_{v_0}$ can be identified with one of the spherical
connected component of the boundary of $Y= W_{r_0}$, triangulated, as we said before, 
by the restriction of  the natural cell decomposition of $Y$ via the truncated tetrahedra of $T$. The cone over Link$(v_0,T)$ with centre $v_0$ is the star Star$(v_0,T)$, so its natural triangulation is the cone over the triangulation $T_{v_0}$ of $S$. 
The fact that this triangulation of  Star$(v_0,T)$ is quasi-regular follows from the fact that 
$T_{v_0}$ is quasi-regular.

\noindent Note that the effect of a $2\to 3$ move on the portions of the 
involved truncated tetrahedra lying on $\partial Y$ consists of three
$2\to 2$ moves and a couple of $1\to 3$ moves. So any $2\to 2$ or
$1\to 3$ move on $T_{v_0}$ can be induced by a suitable $2\to
3$ move around $v_0$.

\subsection{Generalities on flat principal $PSL(2,\C)$-bundles of closed $3$-manifolds}\label{flatbungen}
Let $W$ be a compact closed oriented $3$-manifold, and $\rho$ be a
flat principal bundle over $W$ with structural group $PSL(2,\C)$. We consider
the pair $(W,\rho)$ up to oriented homeomorphisms of $W$
and flat bundle isomorphisms of $\rho$. Equivalently, $\rho$ is
identified with a conjugacy class of representations of the
fundamental group of $W$ in $PSL(2,\C)$, i.e. with a
$PSL(2,\C)$-character of $W$.

\noindent Let $T$ be a triangulation of $W$ with oriented
edges. Denote by $Z^1(T;PSL(2,\C))$ the set of $PSL(2,\C)$-valued simplicial 
$1$-cocycles on $T$. In particular, for such a cocycle $z$ we have $z(-e)=z(e)^{-1}$. A $0$-cochain is a $PSL(2,\C)$-valued function defined on the vertices of $T$. We denote by $[z]$ the equivalence class of $z \in Z^1(T;PSL(2,\C))$ up to cellular coboundaries: two $1$-cocycles $z$ and $z'$ are equivalent if there exists a $0$-cochain $\lambda$ such that for any oriented edge $e$ of $T$ with ordered endpoints $v_0,v_1$, we have $z'(e)=\lambda (v_0)^{-1}z(e)\lambda(v_1)$. We denote this quotient set by $H^1(T;PSL(2,\C))$. The common refinements (subdivisions) of any two triangulations $T$ and $T'$ induce isomorphisms $H^1(T;PSL(2,\C)) \cong H^1(T';PSL(2,\C))$. So $H^1(T;PSL(2,\C))$ can be identified with the set of isomorphism classes of flat principal $PSL(2,\C)$-bundles on $W$, which itself may be described as the reduction of the sheaf cohomology set $H^1(W;\mathcal{C}^{\infty}(PSL(2,\C)))$ to $H^1(W;PSL(2,\C))$ (i.e. where $PSL(2,\C)$ is endowed with the discrete topology). 

\smallskip

\noindent Compact oriented hyperbolic $3$-manifolds with their holonomy furnish a main example of pairs $(W,\rho)$. There are other natural examples $(W,\rho_\alpha)$ 
coming from the ordinary simplicial cohomology of $W$, as follows. Let
us denote by $B$ the Borel subgroup of $SL(2,\C)$ of upper triangular
matrices. There are two distinguished abelian subgroups of $B$:

\smallskip

(1) the {\it Cartan} subgroup $C=C(B)$ of diagonal matrices; it is 
isomorphic to the multiplicative group $\mc^*$ via the map which sends 
$A=(a_{ij}) \in C$ to $a_{11}$; 

\smallskip

(2) the {\it parabolic} subgroup $Par(B)$ of matrices with double 
eigenvalue $1$; it is isomorphic to the additive group $\mc$ via the map  which sends 
$A=(a_{ij})\in Par(B)$ to $x=a_{12}$. 

\smallskip

\noindent Denote by $G$ any such abelian subgroup of $B$. 
There is a natural map $H^1(T;G)\to H^1(T;B)$ induced by the 
inclusion, and $H^1(T;G)$ is endowed with the usual Abelian group structure. 
Note that $H^1(T;Par(B))=H^1(T;\mc)$ is the ordinary (singular or de Rham) first 
cohomology group of $W$.
Hence the inclusion $B\subset SL(2,\C)$ allows us to associate to each $1$-cohomology class $\alpha \in H^1(W,G)$ a pair $(W,\rho_\alpha)$. In particular, we can consider the trivial flat bundle $\rho_0$ on $W$. 
\smallskip

\noindent For our purposes, we need to 
specialize the kind of triangulations, edge orientations and
$PSL(2,\C)$-valued simplicial $1$-cocycles representing flat
$PSL(2,\C)$-bundles.

\subsection{Branchings} \label{branchings}
Let us first specialize the kind of edge orientations. We do it
for ideal triangulations of an arbitrary compact oriented $3$-manifold $Y$ with boundary. Let $P$ be a standard spine of  $Y$, and consider 
the dual ideal triangulation $T=T(P)$. Recall the
notion of abstract tetrahedron of $T$. 

\begin{defi}\label{branchdef}{\rm
A {\it branching} $b$ of $T$ is a choice of orientation for each edge
of $T$ such that on each abstract tetrahedron $\Delta$
of $T$ it is associated to a total ordering $v_0,v_1,v_2,v_3$ of the
(abstract) vertices by the rule: each edge is oriented by the arrow emanating from the smallest endpoint.}
\end{defi}

\noindent Note that for each $j=0,\ldots,3$ there are exactly $j$
$b$-oriented edges incoming at the vertex $v_j$; hence there are only one
source and one sink of the branching. This is equivalent to saying that for any $2$-face
$f$ of $\Delta$ the boundary of $f$ is not coherently oriented. In dual terms, a branching is a choice of orientation for each region of $P$ such that for each edge of $P$ we
have the same induced orientation only twice. In particular, the edges of $P$
have an induced prevailing orientation.
\smallskip

\noindent Branchings, mostly in terms of spines,
have been widely studied and applied in \cite{BP2,BP3,BP4}. A
branching of $P$ gives it the extra-structure of an
embedded and oriented (hence normally oriented) {\it branched surface}
in Int($Y$). Moreover, a branched spine $P$ carries a suitable positively
transverse {\it combing} of $Y$ (ie. a non-vanishing vector field).

\noindent Given a
branching $b$ on a oriented tetrahedron $\Delta$ (realized in $\R^3$
as stipulated in Section \ref{trigen}), denote by $E(\Delta)$ the set of $b$-oriented
edges of $\Delta$, and by $e'$ the edge opposite to $e$. We put
$e_0=[v_0,v_1]$, $e_1=[v_1,v_2]$ and $e_2=[v_0,v_2]=-[v_2,v_0]$. This
fixed ordering of the edges of the $2$-face opposite to the vertex $v_3$
will be used all along the paper. The ordered triple of edges
\begin{equation}\label{order}
(e_0=[v_0,v_1],e_2=[v_0,v_2],e_1'=[v_0,v_3])
\end{equation} 
departing from $v_0$
defines a \emph{$b$-orientation} of $\Delta$. This orientation may or
may not agree with the orientation of $Y$.  In the first case we say
that $\Delta$ is of index $*_b=1$, and it is of index $*_b=-1$
otherwise.  The $2$-faces of $\Delta$ can be named by their opposite
vertices. We orient them by working as above on the boundary of each
2-face $f$: there is a $b$-ordering of the vertices of $f$, and an
orientation of $f$ which induces on $\partial f$ the prevailing
orientation among the three $b$-oriented edges. This $2$-face
orientation corresponds dually to the orientation on the edges of $P$
mentioned above.  These considerations apply to each abstract
tetrahedron of any branched triangulation $(T,b)$ of $Y$.
\medskip

\noindent {\bf Branching's existence and transit.}  
In general a given ideal triangulation $T$ of $Y$ may not admit any
branching. 

\begin{defi}\label{btransit}{\rm 
Given {\it any} choice $g$ of edge-orientations on $T$ and
any move $T \to T'$, a transit  $(T,g)\to (T',g')$ is given by
any choice $g'$ of edge-orientations on $T'$ which agrees with
$g$ on the common edges of $T$ and $T'$. This makes a
{\it branching transit} if both $g$ and $g'$ are
branchings. We will use the same terminology for moves on
branched standard spines.}
\end{defi} 
Concerning the existence of branched triangulations there is the following result:
\begin{prop} \label{brancheable} {\rm \cite[Th. 3.4.9]{BP2}}
For any system $g$ of edge-orientations on $T$ there exists a finite
sequence of positive $2\to 3$ transits such that the final $(T',g')$
is actually branched.
\end {prop}
On another hand, any quasi-regular
triangulation $T$ of a closed $3$-manifold $W$ admits branchings of a
special type,  
defined by fixing any total ordering of its vertices and
by stipulating that the edge $[v_i,v_j]$ is positively
oriented iff $j > i$. These branchings are called total ordering branchings. Any quasi-regular move which preserves the
number of vertices also preserves the total orderings on the set of
vertices, hence it obviously induces \emph{total ordering branching transits}.  If it increases the number of vertices, one can extend to
the new vertex, in several different ways, the old total orderings of
the set of vertices. Again, any of these ways induces a
total odering branching transit.  If $(T,b)$ is an arbitrary branched
triangulation of $Y$ (ie. $T$ is not necessarily quasi-regular nor
$b$ is of total ordering type) and $T\to T'$ is either a positive $2\to 3$,
$0\to 2$ or bubble move, then it can be completed, sometimes in
different ways, to a branched transit $(T,b)\to (T',b')$. Any of these
ways is a possible transit. On the contrary, it is easily seen that a
negative $3\to 2$ or $2\to 0$ move may not be ``branchable'' at all.

\begin{figure}[ht]
\begin{center}
\includegraphics[width=12cm]{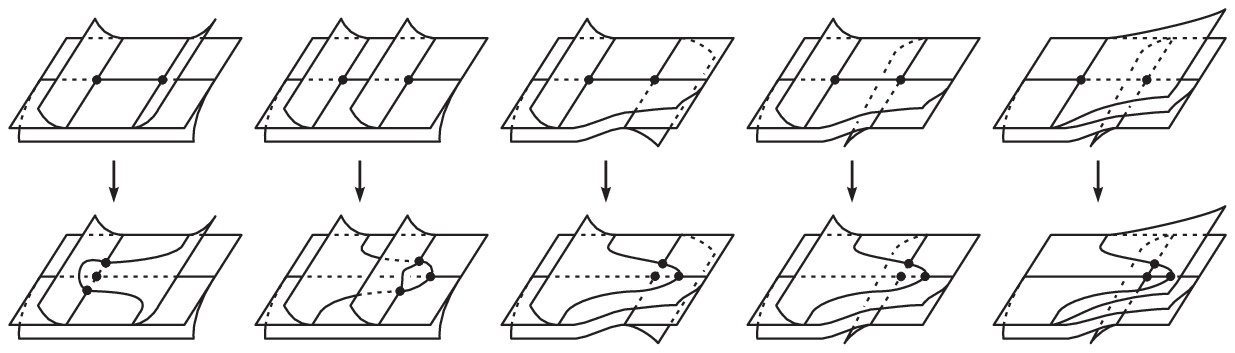}
\caption{\label{fig:slide} $2\to 3$ sliding moves.} 
\end{center}
\end{figure}

\begin{figure}[ht]
\begin{center}
\includegraphics[width=7cm]{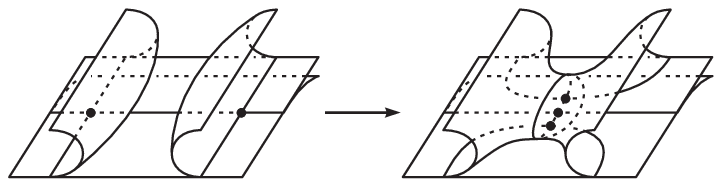}
\caption{\label{fig:bump} $2\to 3$ bumping move.} 
\end{center}
\end{figure}

\begin{figure}[ht]
\begin{center}
\includegraphics[width=8cm]{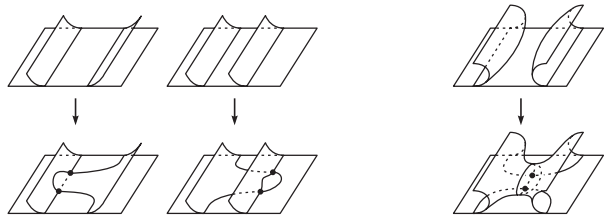}
\caption{\label{fig:lune} branched lune-moves.} 
\end{center}
\end{figure}

\noindent For the sake of clarity, we show in Fig. \ref{fig:slide} - Fig. 
\ref{fig:lune} the whole set of $2\to 3$ and $0 \to 2$ 
(dual) branched transits on standard spines, up to evident symmetries.
Note that the middle sliding move in Fig. \ref{fig:slide} corresponds
dually to the branched triangulation move shown in Fig. \ref{idealt}.
Following \cite{BP2}, one can distinguish two families of 
branched transits: the {\it sliding moves}, which actually preserve the
positively transverse combing mentioned at the beginning of this
subsection, and the {\it bumping moves}, which eventually change it. We
shall not exploit this difference in the present paper; see however
Remark \ref{phase}. 

\smallskip

\noindent Finally, we note that for the proof of the main theorem \ref{QHIinv} it is enough to use only total ordering branchings, but we need to consider general branchings to extend the construction of the QHI to other situations (see the first point in Section \ref{qhicomp}, and the discussion on cusped manifolds in Section \ref{VOLCONJ}).

\subsection{$\Dd$-triangulations for $(W,\rho)$ and their idealization}
\label{DTWID}
We will now select certain generic $PSL(2,\C)$-valued $1$-cocycles on branched
quasi-regular triangulations $(T,b)$ of $W$, so as to define the $\Dd$-triangulations. 

\begin{defi}\label{Dtriang} 
{\rm A $\Dd$-{\it triangulation} for the pair $(W,\rho)$ consists of a
triple $\Tt=(T,b,z)$ where: $T$ is a quasi-regular triangulation
of $W$; $b$ is a branching of $T$; $z$ is a $PSL(2,\C)$-valued
$1$-cocycle on $(T,b)$ representing $\rho$, such that $(T,b,z)$ is {\it
idealizable} (see Def. \ref{idealizable}).}
\end{defi}
The name `$\Dd$-triangulation' refers to the fact that they are
``decorated'' by the branching and the cocycle (and, later in Def. \ref{dist}, ``distinguished'' by an Hamiltonian link). If $z$ is $PSL(2,\C)$-valued $1$-cocycle on $(\Delta,b)$, we write $z_j
= z(e_j)$ and $z'_j = z(e'_j)$. Then, one reads for instance
the cocycle condition
on the $2$-face opposite to $v_3$ as $ z_0z_1z_2^{-1} = 1$.  This
holds for each abstract tetrahedron of any branched triangulation
$(T,b)$ of $W$ and for (the restrictions of) any $PSL(2,\C)$-valued
$1$-cocycle $z$ on $(T,b)$.

\noindent Consider the half space model of the hyperbolic space
$\mh^3$. We orient it as an open set of $\R^3$. The natural boundary
$\partial \bar{\mh}^3=\mc\mathbb{P}^1 = \C \cup\{ \infty \}$ of $\mh^3$ is
oriented by its complex structure. We realize $PSL(2,\C)$ as the group of orientation preserving isometries of $\mh^3$, with the corresponding conformal action on $\mc\mathbb{P}^1$.

\begin{defi}\label{idealizable}
{\rm Let $(\Delta,b,z)$ be a branched tetrahedron endowed with a
$PSL(2,\C)$-valued $1$-cocycle $z$. It is {\it idealizable} iff
$$ u_0=0,\ u_1= z_0(0),\ u_2= z_0z_1(0),\ u_3= z_0z_1z'_0(0) $$ are 4
distinct points in $\C \subset \mc\mathbb{P}^1= \partial
\bar{\mh}^3$. These 4 points span a (possibly degenerate) hyperbolic
ideal tetrahedron with ordered vertices. A triangulation $(T,b,z)$ is
{\it idealizable} iff all its abstract tetrahedra $(\Delta_i,b_i,z_i)$ are
idealizable.}
\end{defi} 

\noindent 
If $(\Delta,b,z)$ is idealizable, for all $j=0,1,2$ one can associate
to $e_j$ and $e'_j$ the same {\it cross-ratio} modulus $w_j\in
\C\setminus \{0,1\}$ of the hyperbolic ideal tetrahedron spanned by 
$(u_0,u_1,u_2,u_3)$; we refer to \cite[Ch. 5]{BP1} for details on the meaning and the role of moduli in hyperbolic geometry. We have (indices mod($\mz/3\mz$)):
\begin{equation}\label{mod1}
w_{j+1} = 1/(1-w_j)
\end{equation}
and
$$w_0 = ( u_2- u_1)u_3/u_2(u_3- u_1)\ .$$ 
Let us write $p_0= u_1(u_3-u_2)$, $p_1=( u_2- u_1)u_3$, and $p_2= -u_2(u_3-u_1)$. Then 
\begin{equation}\label{mod2}
w_j = - p_{j+1}/p_{j+2}\ .
\end{equation} 
\noindent Set $w=(w_0,w_1,w_2)$ and call it a {\it modular triple}. The ideal
tetrahedron spanned by $(u_0,u_1,u_2,u_3)$ is non-degenerate iff the imaginary parts of the $w_j$'s are
not equal to zero; in such a case they share the same sign $*_w = \pm 1$.

\begin{defi}\label{idealizations}
{\rm We call $(\Delta,b,w)$ the {\it idealization} of the idealizable
$(\Delta,b,z)$, and identify it with the branched tetrahedron  in $\bar{\mh}^3$ spanned by $(u_0,u_1,u_2,u_3)$. For any $\Dd$-triangulation $\Tt=(T,b,z)$ of
$(W,\rho)$, its {\it idealization} $\Tt_{\Ii}=(T,b,w)$ is given by the
family $\{(\Delta_i, b_i, w_i)\}$ of idealizations of the
$(\Delta_i,b_i, z_i)$'s. We say that $\Tt_{\Ii}$ is an
$\Ii$-triangulation for $(W,\rho)$. It is non-degenerate if each
$\{(\Delta_i, b_i, w_i)\}$ is non-degenerate.}
\end{defi}

\begin{remarks}\label{idealrem}
{\rm 
(1) We could incorporate the non-degeneracy assumption into
the notion of idealizable tetrahedron. All the constructions of the present
paper would run in the same way. The non-degenerate assumption simplifies 
the exposition and also certain proofs concerning the dilogarithmic invariant of 
$(W,\rho)$ considered in \cite{BB2}.

\smallskip

\noindent (2) Since $PSL(2,\mc)$ acts on $\mc\mathbb{P}^1$ via Moebius transformations $z_j: x \mapsto (a_jx +b_j)/(c_jx+d_j)$, it is immediate to formulate for any given quasi-regular triangulation $T$ a simple system of algebraic equalities on the entries of the $z_j$'s, whose zero set describes non-idealizable cocycles.

\smallskip

\noindent (3) In \cite{BB1} we have used so-called \emph{full} $B$-valued
$1$-cocycles $z$ to construct the QHI for $B$-characters. `Full' means that for any edge
$e$ the upper diagonal entry $x(e)$ of $z(e)$ is non-zero. It is easy to verify that a $B$-valued $1$-cocycle is full iff it is idealizable. The idealization we proposed for such cocycles was in fact a specialization of the present general procedure. We can simply write the moduli for the idealization
of a $\Dd$-tetrahedron with a full $B$-valued $1$-cocycle as $w_j = - q_{j+1}/q_{j+2}$, where $q_j=x(e_j)x(e_j')$ for $j=0,1$, and
$q_2=-x(e_2)x(e_2')$ (beware that $p_i \ne q_i$).

\smallskip

\noindent (4) It follows from the cocycle condition or from the relation (\ref{mod1}) that $p_0 + p_1 + p_2 =0$ (and also that $q_0+q_1+q_2=0$ - see remark (3)).
}
\end{remarks}

\noindent The following lemma is immediate:

\begin{lem} \label{esistenza} 
For any $PSL(2,\C)$-character $\rho$, any quasi-regular branched
triangulation $(T,b)$ of $W$ can be completed to a $\Dd$-triangulation
$(T,b,z)$ for the pair $(W,\rho)$.
\end{lem}
\noindent 
In fact, given any $PSL(2,\mc)$-valued $1$-cocycle, one can perturb it by the coboundaries
of generic $0$-cochains which are injective on the vertices of $T$,
so that we get idealizable $1$-cocycles.
\medskip

\noindent{\bf Tetrahedral symmetries.} The
idealization has a good behaviour with respect to a change of branching (the
`tetrahedral symmetries'). Indeed, we have:

\begin{lem} \label{Iisym}  
Denote by $S_4$ the permutation group on four elements.  A permutation
$p \in S_4$ of the vertices of an idealizable tetrahedron
$(\Delta,b,z)$ gives another idealizable tetrahedron $(\Delta,b',z')$.
The permutation turns the idealization $(\Delta,b,w)$ into an isometric $(\Delta,b',w')$, where for each edge $e$ of $\Delta$ we have $w'(e)=w(e)^{\epsilon(p)}$ and $\epsilon(p)$ is the signature of $p$.
\end{lem}
\noindent {\it Proof.} Consider for instance the transposition $(0,1)$. 
It turns the (ordered) set of points $ 0, z_0(0), z_0z_1(0), z_0z_1z'_0(0)$ into $0, (z_0)^{-1}(0), z_1(0), z_1z'_0(0)$.  By applying on this
second set the hyperbolic isometry $z_0$, one gets the
first set after the transposition of $0$ and $z_0(0)$. Things go
similarly for any other permutation.  Then the lemma follows immediately, due
to the behaviour of cross-ratios with respect to vertex permutations. \hfill $\Box$

\medskip

\noindent {\bf Hyperbolic edge compatibility.} 
We are now concerned with an important global property of the
idealized triangulations $\Tt_{\Ii}$. Before to state it, let us
stress that when dealing with modular triples one has to be careful
with the orientations. Recall that every $\Ii$-tetrahedron
$(\Delta,b,w)$ is oriented by definition; in the case of an
$\Ii$-triangulation this is given by the orientation of $W$. There is
also the $b$-orientation encoded by the sign $*=*_b$. The idealization
`physically' realizes the vertices of $\Delta$ on $\partial
\bar{\mh}^3$, with the ordering induced by $b$. 
When the spanned ideal tetrahedron is non-degenerate, the $b$-orientation may
or may not agree with the one induced by the fixed orientation of $\mh^3$, which is encoded by the
sign $*_w$ of the modular triple.

\noindent 
Let $\Tt_{\Ii}=(T,b,w)$ be an $\Ii$-triangulation. The preceding discussion shows that the contribution of each $(\Delta_i,b_i,w_i)$ to any computation with
the moduli is given by the $w(e)^{*}$'s, where $e$ is any edge of
$\Delta_i$ and $*= *_{b_i}$.  The next Lemma \ref{compat} is a first
concretization of this fact (see also the notion of $\Ii$-transit
below). Denote by $E(T)$ the set of edges of $T$, by $E_{\Delta}(T)$ the whole
set of edges of the associated abstract tetrahedra $\{\Delta_i\}$, and
by $\epsilon_T:E_{\Delta}(T) \longrightarrow E(T)$ the natural
identification map.

\begin{figure}[ht]
\begin{center}
\includegraphics[width=3.5cm]{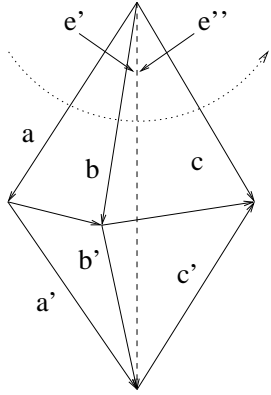}
\caption{\label{compatrel} the compatibility relation around an edge.} 
\end{center}
\end{figure}

\begin{lem} \label{compat} For any edge $e \in E(T)$ we have 
$\textstyle \prod_{a \in \epsilon_T^{-1}(e)} w(a)^{*} =1$, where $*=\pm
1$ according to the $b$-orientation of the tetrahedron $\Delta_i$ that
contains $a$.
\end{lem}

\noindent{\it Proof.} Looking at Star$(e,T)$ we see that up to a sign two consecutive moduli partially compensate along the common face of the corresponding tetrahedra. For instance, in Fig. \ref{compatrel} the left (resp. right) tetrahedron is negatively (resp. positively) $b$-oriented; we have $w(e')^{-1}w(e'') = (ab')(-bc')/(ba')(b'c) = -ac'/a'c$. This and Lemma \ref{Iisym} show that the same holds true when the two tetrahedra are simultaneously positively or negatively $b$-oriented. Continuing this way around $e$, we end up with $\textstyle \prod_{a \in \epsilon_T^{-1}(e)} w(a)^* =\pm 1$. Each $-1$ contribution comes from a tetrahedron where the $b$-orientations of the two faces containing $e$ are opposite (that is, when the corresponding $a$ is $e_2$ or $e_2'$). Since $W$ is orientable, a short closed loop about $e$ may only meet an even number of such tetrahedra. This gives the result.\hfill $\Box$

\medskip

\noindent This Lemma means that around each edge the signed moduli verify the usual compatibility condition needed when one
tries to construct hyperbolic $3$-manifolds by gluing hyperbolic ideal
tetrahedra. So the $\Ii$-triangulations have the following geometric interpretation. Given an $\Ii$-triangulation $\Tt_{\Ii}=(T,b,w)$ of $(W,\rho)$, lift $T$ to a cellulation $\widetilde{T}$ of the universal cover $\widetilde{W}$, and fix a base point $\tilde{x}_0$ in the $0$-skeleton of $\widetilde{T}$; denote by $x_0$ the projection of $\tilde{x}_0$ onto $W$. Then, for any tetrahedron in $\widetilde{T}$ that contains $\tilde{x}_0$, use the moduli of the corresponding $(\Delta_i,b_i,w_i) \in \Tt_{\Ii}$ to define an hyperbolic ideal tetrahedron. Do this by respecting the gluings in $\widetilde{T}$. Starting from the vertices adjacent to $\tilde{x}_0$ and continuing in this way, we construct an image in $\bar{\mh}^3$ of a complete lift of $T$ in $\widetilde{T}$, having one tetrahedron in each $\pi_1(W)$-orbit. The key point is that Lemma \ref{compat} implies that for any two paths of tetrahedra in $\widetilde{T}$ having a same starting point, we get the same end point. This construction extends to a piecewise-linear map $D: \widetilde{W} \rightarrow \bar{\mh}^3$, equivariant with respect to the action of $\pi_1(W)$ and $PSL(2,\mc)$. So we eventually find: a representation $\tilde{\rho}: \pi_1(W,x_0) \rightarrow PSL(2,\mc)$ with character $\rho$ and satisfying $D(\gamma(x))=\tilde{\rho}(\gamma) D(x)$ for each $\gamma \in PSL(2,\mc)$; a piecewise-straight continuous section of the flat bundle $\widetilde{W} \times_{\tilde{\rho}} \bar{\mh}^3 \rightarrow W$, with structural group $PSL(2,\mc)$ and total space the quotient of $\widetilde{W} \times  \bar{\mh}^3$ by the diagonal action of $\pi_1(W)$ and $\tilde{\rho}$. The map $D$ behaves formally as a developing map for a $(PSL(2,\mc),\mh^3)$-structure on $W$ (see e.g. \cite[Ch. B]{BP1} for this notion).   

\medskip

\noindent {\bf $\Dd$- and $\Ii$-transits.} 
We consider now moves on $\Dd$-triangulations $\Tt=(T,b,z)$ and
$\Ii$-triangulations $\Tt_{\Ii}=(T,b,w) $ for the pair $(W,\rho)$,
called $\Dd$-{\it transits} and $\Ii$-{\it transits} respectively.
They are supported by the bare triangulation moves mentioned in Section
\ref{trigen}, but they also include the transits of the respective
extra-structures. First of all we require that they are quasi-regular moves.
We stress that this is not an automatic fact; on the contrary this
leads to one main technical complication in the proofs. Then we require that
$(T_0,b_0) \leftrightarrow (T_1,b_1)$ 
is a branching transit in the sense of Def. \ref{btransit}. 
 
\begin{defi}\label{ztransit} 
{\rm Let $(T_0,b_0)$, $(T_1,b_1)$ be branched quasi-regular triangulations and $z_k \in
Z^1(T_k;PSL(2,\C))$, $k=0,1$. We have a {\it
cocycle transit} $(T_0,z_0)
\leftrightarrow (T_1,z_1)$ if $z_0$ and $z_1$ agree on the common
edges of $T_0$ and $T_1$. This makes an {\it idealizable cocycle transit} if both $z_0$ and $z_1$ are idealizable $1$-cocycles, and in this case we say that $(T_0,b_0,z_0) \leftrightarrow (T_1,b_1,z_1)$ is a $\Dd$-transit.}
\end{defi}
It is not hard to see that $z_0$ and $z_1$ as above represent the same
flat bundle $\rho$. Note that for $2\rightarrow 3$ and $0 \rightarrow
2$ moves, given $z_k$ there is only one (resp. at most one) $z_{k+1}$
with this property. We stress that in some special cases a $2\to 3$
transit of an idealizable cocycle can actually not preserve the
idealizability, but generically this does not hold. For positive
bubble moves there is always an infinite set of possible (idealizable)
cocycle transits. The following lemma shows that the $\Dd$-transits are generic.

\begin{lem} \label{generique} Let $(T,b)$ be a branched quasi-regular 
triangulation of $W$. Suppose that $(T,b)=(T_1,b_1)\to \dots \to
(T_s,b_s)=(T',b')$ is a finite sequence of quasi-regular
$2\leftrightarrow 3$ branching transits. Then for each $T_i$ there
exists a dense open set $U_i$ of $PSL(2,\C)$-valued $1$-cocycles, in the quotient topology of $PSL(2,\C)^{r_1(T_i)}$ as a space of matrices ($r_1(T_i)$ being the
number of edges of $T_i$), such that for every $z_i\in U_i$, $(T_i,b_i,z_i)$ is a $\Dd$-triangulation, and the transit $T_i\to
T_{i+1}$ maps $U_{i}$ into $U_{i+1}$.  Moreover each class $\alpha \in H^1(W;PSL(2,\C)) \cong H^1(T_i;PSL(2,\C))$ can be represented by cocycles in $U_i$.
\end{lem}

\noindent {\it Proof.} Each $2\leftrightarrow 3$, $0\leftrightarrow 2$ or negative bubble transit 
$(T_i,b_i,z_i) \to (T_{i+1},b_{i+1},z_{i+1})$ defines an algebraic surjective map from $Z^1(T_i;PSL(2,\C))$ to $Z^1(T_{i+1};PSL(2,\C))$. 
Since all edges of $T_{i+1}$ have 
distinct vertices, there are no trivial (two term) cocycle relations on $T_{i+1}$. 
Hence the set of idealizable cocycles for which a transit fails to be idealizable is contained in a proper algebraic subvariety 
of $Z^1(T_i;PSL(2,\C))$. 
Working by induction on $s$ we get the conclusion.\hfill $\Box$
\medskip

\noindent Let us now consider the transits for the idealized triangulations.
Consider the convex hull of five distinct points $u_0,u_1,u_2,u_3,u_4
\in \partial \bar{\mathbb{H}}^3$, with the two possible triangulations
$Q_0$ $Q_1$ made of the oriented hyperbolic ideal
tetrahedra $\Delta^i$ obtained by omitting $u_i$, $i=0,\ldots,4$. An edge $e$ of $Q_i
\cap Q_{i+1}$ belongs to one tetrahedron of $Q_i$ iff it belongs to two 
tetrahedra of $Q_{i+1}$. 
Then, the modulus of $e$ in $Q_i$ is the product of the two moduli of $e$ in
$Q_{i+1}$. Also, the product of the moduli on the central edge of
$Q_1$ is equal to $1$.

\noindent Let $T \to T'$ be a $2\to 3$ move. 
Consider the two (resp. three) abstract tetrahedra of $T$ (resp. $T'$)
involved in the move. They determine subsets $\widetilde{E}(T)$ of
$E_{\Delta}(T)$ and $\widetilde{E} (T')$ of $E_{\Delta}(T')$. Put
$\widehat{E} (T) = E_{\Delta}(T) \setminus \widetilde{E}(T)$ and $\widehat{E} (T') = E_{\Delta}(T') \setminus \widetilde{E}(T')$. Clearly one can identify $\widehat{E} (T)$ and $\widehat{E}
(T')$. The above configurations $Q_0$ and $Q_1$ and
the considerations made before Lemma \ref{compat} lead to the
following definition:

\begin{figure}[ht]
\begin{center}
\includegraphics[width=11cm]{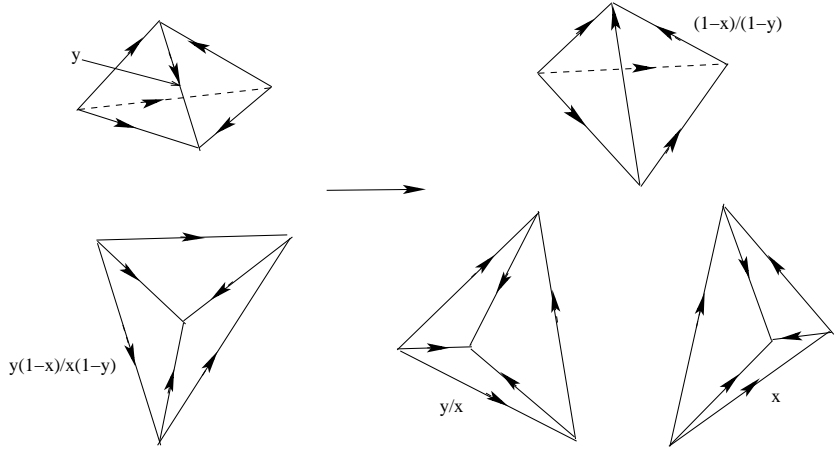}
\caption{\label{idealt} a $2 \leftrightarrow 3$ ideal transit.} 
\end{center}
\end{figure}

\begin{defi} \label{id-transit} 
{\rm A $2 \to 3 \ $ {\it $\Ii$-transit} $(T,b,w) \rightarrow
(T',b',w')$ of $\Ii$-triangulations for a pair $(W,\rho)$ is such that:

\smallskip

\noindent 1) $w$ and $w'$ agree on $\widehat{E} (T) = \widehat{E} (T')$;

\noindent 
2) for each common edge $e\in \epsilon_T(\widetilde{E} (T))\cap
   \epsilon_{T'}(\widetilde{E} (T'))$ we have}
\begin{equation}\label{ideqmod} \prod_{a\in \epsilon_T^{-1}(e)}w (a)^*=
 \prod_{a'\in \epsilon_{T'}^{-1}(e)}w' (a')^*
\end{equation}
\noindent {\rm where $*=\pm 1$ according to the $b$-orientation of the abstract tetrahedron containing $a$ (resp. $a'$). 
We have a $0 \to 2$ (resp. \emph{bubble}) \emph{$\Ii$-transit} if the
above first condition is satisfied, and we replace the second by:

\smallskip

\noindent 
2') for each edge $e\in \epsilon_{T'}(\widetilde{E} (T'))$ we have}
\begin{equation}\label{ideq2mod} \prod_{a'\in \epsilon_{T'}^{-1}(e)}w' 
(a')^*=1\ .
\end{equation}
\end{defi}  

\noindent 
$\Ii$-transits for negative $3 \to 2$ moves are defined in exactly the
same way, and for negative $2 \to 0$ and bubble moves $w'$ is defined
by simply forgetting the moduli of the two disappearing
tetrahedra. The condition (1) above implies that the product of the
$w'(a')^*$'s around the new edge is equal to $1$.  A $2
\leftrightarrow 3$ $\Ii$-transit is shown in Fig. \ref{idealt}; we
only indicate the first component `$w_0$' of each modular triple. In general,
the relations (\ref{ideqmod}) may imply that $w$ or $w'$ equals $0$ or
$1$ on some edges.  In that case, the $2 \leftrightarrow 3$
$\Ii$-transit fails. In particular, in Fig. \ref{idealt} we assume
that $x \ne y$.

\noindent 
Note that for $2 \leftrightarrow 3$ $\Ii$-transits $w'$ is uniquely
determined by $w$. On the contrary, there is one degree of freedom for (positive) $ 0 \to 2$ and bubble $\Ii$-transits. The relation (\ref{ideq2mod}) simply means that such transits give \emph{the same}
modular triples to the two new tetrahedra, for their $b$-orientations
are opposite.

\smallskip

\noindent The next proposition states the remarkable fact that $\Dd$-transits
and $\Ii$-transits together with the idealization make commutative
diagrams, that is the $\Dd$-transits {\it dominate} the
$\Ii$-transits.

\begin{prop}\label{DdomI} Consider a fixed pair $(W,\rho)$, and denote by $\Ii$
the idealization map $\Tt \to \Tt_{\Ii}$ on its $\Dd$-triangulations. For any $\Dd$-transit $\mathfrak{d}$ there exists an $\Ii$-transit $\mathfrak{i}$ (resp. for any
$\mathfrak{i}$ there exists $\mathfrak{d}$) such that $\mathfrak{i}
\circ \Ii = \Ii \circ \mathfrak{d}$.
\end{prop}
For $2 \leftrightarrow 3$ transits there is also an uniqueness
statement. 
\smallskip

\noindent{\it Proof.} By using the tetrahedral symmetries
of Lemma \ref{Iisym}, it is enough to show the proposition for one 
branching transit configuration (for instance the one of Fig.
\ref{idealt}). The idealization map defines embeddings of the $\Dd$-tetrahedra of this configuration as branched ($b$-oriented) ideal tetrahedra in $\mh^3$. Since orientation preserving isometries do not alter the moduli, the union of these tetrahedra may be viewed as the convex hull of five distinct ordered points on $\mc\mathbb{P}^1$, such that the ordering induces the branching on each tetrahedron. Then the verification follows immediately from the
definition of the idealization. \hfill $\Box$
\medskip

\noindent 
Note that the possible failures of $2\to 3$ transits of idealizable
cocycles that we mentioned after Def. \ref{ztransit} exactly correspond to the failures of $2\to
3$ $\Ii$-transits (for instance when $x=y$ in Fig.  \ref{idealt}).

\section{Quantum dilogarithms and basic state sums for pairs $(W,\rho)$}
\label{QDILOGSS}
\noindent 
Let $N=2p+1>1$ be a fixed odd positive integer, and put
$\zeta=\exp(2i\pi/N)$. The quantum algebraic origin of
the Faddeev-Kashaev's matrix quantum
dilogarithms is discussed in the Appendix
(Section \ref{APP}) of the present paper. 
Here we forget this origin, and, for the reader's
convenience, we simply introduce the special functions needed for
defining {\it basic state sums}
supported by the $\Ii$-triangulations $\Tt_{\Ii}$ of any pair $(W,\rho)$.
The main property of these basic state sums is to be invariant for some {\it specific}
instances of $\Ii$-transits.

\smallskip

\noindent We denote by $g$ the analytic function defined for any complex number $x$ with $\vert x \vert <1$ by
$$g(x) := \prod_{j=1}^{N-1}(1 - x\zeta^j)^{j/N}$$
\noindent and set $h(x) := x^{-p}g(x)/g(1)$ when $x$ 
is non-zero (one computes that $\vert g(1) \vert = N^{1/2}$).  We
shall still write $g$ for its analytic continuation to the complex
plane with cuts from the points $x = \exp (i\epsilon)\zeta^k,\ k
=0,\ldots,\ N-1$, $\epsilon \in \mathbb{R}$, to infinity. Hereafter we
will implicitly assume that $\epsilon$ is such that the cuts are away
from the points where $g$ is evaluated (things will not depend on this choice).

\noindent Consider the curve $\Gamma=\{x^N + y^N = z^N\} \subset \mathbb{C}P^2$
(homogeneous coordinates), and the rational functions on $\Gamma$ given for any $
n \in \mathbb{N}$ by
\begin{eqnarray} \label{omeg0}
\omega(x,y,z \vert n) = \prod_{j=1}^n \frac{(y/z)}{1-(x/z)\zeta^j} \quad .
\end{eqnarray}
\noindent 
These functions are periodic in their integer argument, with period
$N$. Denote by $\delta$ the $N$-periodic Kronecker symbol,
i.e. $\delta(n) = 1$ if $n \equiv 0$ mod($N$), and $\delta(n) = 0$
otherwise. Set $[x] = N^{-1}\ (1-x^N)/(1-x)$.

\smallskip

\noindent The elementary building blocks of the basic state sums are the $N^2 \times N^2$-matrix valued quantum dilogarithms and their
inverses, whose matrix entries are the rational
functions defined on the curve $\Gamma$ by
$$R(x,y,z)_{\alpha,\beta}^{\gamma,\delta} = h(z/x)\
\zeta^{\alpha\delta+\frac{\alpha^2}{2}}\ \omega(x,y,z \vert \gamma-\alpha) \ \delta(\gamma + 
\delta - \beta)$$
$$\bar{R}(x,y,z)_{\gamma,\delta}^{\alpha,\beta} =  
\frac{[x/z]}{h(z/x)}\ \zeta^{-\alpha\delta-\frac{\alpha^2}{2}}\ \frac{\delta(\gamma + 
\delta - \beta)}{\omega(\frac{x}{\zeta},
y,z\vert \gamma-\alpha)}\ .$$

\noindent We can interpret these matrices as functions of 
$\Ii$-tetrahedra as follows. Let $(\Delta,b,w)$ be an
$\Ii$-tetrahedron. The $1$-skeleton of the cell decomposition of
$\Delta$ dual to the canonical triangulation with 4 vertices is made
of 4 edges incident at an interior point of $\Delta$. As we said in
Subsection \ref{branchings}, the orientations of these edges are
complementary to the $b$-orientations of the dual
$2$-faces of $\Delta$. Two of them are pointing inwards $\Delta$, and
the others are pointing outwards. So they form two distinguished
pairs. Let us order the two edges of each pair as the corresponding $2$-faces of $\Delta$ (ordered by the opposite vertices). We can associate to both
ordered pairs a copy of $\C^N \otimes \C^N$ (with the standard basis),
which we denote respectively by $I_1 \otimes I_2$ (for `inwards') and
$O_1 \otimes O_2$ (for `outwards').

\noindent Write $w_i = -p_{i+1}/p_{i+2}$ (indices
mod($\mz/3\mz$)) as in (\ref {mod2}). Recall from Remark \ref{idealrem} (4)
that $p_0 + p_1 + p_2 =0$. Fix common determinations of the $N$-th roots of the $p_i$'s, and denote them by $p_i'$. We define a matrix $\mathfrak{L}_N(\Delta,b,w): I_1 \otimes I_2 \longrightarrow O_1 \otimes O_2$ by
\begin{eqnarray}
\mathfrak{L}_N(\Delta,b,w) = \left\lbrace \begin{array}{l}
R(p_1',p_0',-p_2') \quad {\rm if} \quad *=1 \\ \\ 
\bar{R}(p_1',p_0',-p_2') \quad {\rm if} \quad *=-1  
\end{array} \right. \nonumber
\end{eqnarray}
\noindent 
\noindent where $*=\pm1$ according to the $b$-orientation of $\Delta$.
Since $\mathfrak{L}_N(\Delta,b,w)$ is homogeneous in the $p_i'$'s, it only
depends on $(b,w)$. 

\medskip

\noindent Let $\Tt_{\Ii} = (T,b,w)$ be any $\Ii$-triangulation for
$(W,\rho)$. Let us consider the $1$-skeleton $C$ of the cell
decomposition dual to $T$, with the edges oriented as above. By
associating to each $(\Delta_i,b_i,w_i)$ the corresponding operator
$\mathfrak{L}_N(\Delta_i,b_i,w_i)$, one gets an {\it operator network}
whose complete contraction gives a scalar
$\mathfrak{L}_N(\Tt_{\Ii})\in \C$ (note that there is no edge with
free ends in $C$). This has an explicit expression as a {\it state
sum} as follows. A {\it state} is a function defined on the edges of $C$, with
values in $\{0,\dots, N-1\}$. Any state
$\alpha$ determines an entry (a $6j$-symbol)  
$\mathfrak{L}_N(\Delta_i,b_i,w_i)_\alpha$ of $\mathfrak{L}_N(\Delta_i,b_i,w_i)$,  for each $(\Delta_i,b_i,w_i)$.
Set 
$$\mathfrak{L}_N(\Tt_{\Ii})_\alpha = 
\prod_i \mathfrak{L}_N(\Delta_i,b_i,w_i)_\alpha $$ and
\begin{equation}\label{statesum}
\mathfrak{L}_N(\Tt_{\Ii})= 
\sum_\alpha \mathfrak{L}_N(\Tt_{\Ii})_\alpha \quad .
\end{equation} 

\noindent Given any maximal tree $\tau$ in $C$, we can consider the polyhedron $P_{\tau}$ obtained by cutting $T$ along the faces dual to the edges of $C \setminus \tau$. Fix an ordering of these faces. Then we can write $\mathfrak{L}_N(\Tt_{\Ii})$ as the trace of the operator obtained by composing the $\mathfrak{L}_N(\Delta_i,b_i,w_i)$'s along the faces dual to the edges of $\tau$. The domain (resp. target) space of this operator is the tensor product of one copy of $\mc^N$ for each face of $\partial P_{\tau}$ whose dual edge points inwards (resp. outwards) $P_{\tau}$, with the same ordering. We have the following key facts:

\medskip

(a) Straightforward computations show that $\mathfrak{L}_N(\Delta,b,w)$ does not respect the tetrahedral
symmetries, i.e. it is not invariant if we change the branching.
\medskip

(b) $\mathfrak{L}_N(\Tt_{\Ii})$ is invariant {\it only} for some
peculiar instances of $2 \leftrightarrow 3$ $\Ii$-transits. One among them is shown in Fig. \ref{idealt}. This instance corresponds to the {\it basic pentagon identity} satisfied by $\mathfrak{L}_N(\Delta,b,w)$ (see (\ref{basicpentagon}) in the Appendix).
 
\medskip

\noindent Before we overcome these problems, let us disgress a bit to motivate and explain the approach we will follow.

\medskip

\noindent {\bf Quantum vs. classical dilogarithms.} There is a `classical' analogue of $\mathfrak{L}_N(\Tt_{\Ii})$, which we now describe. We refer to \cite{BB2} for details. Denote by $\log$ the standard branch of the logarithm, with arguments
in $]-\pi,\pi]$.  Put $\mathfrak{D} = \mc \setminus \{(-\infty;0) \cup
(1;+\infty) \}$.  The \emph{Rogers dilogarithm} is the complex
analytic function defined over $\mathfrak{D}$ by
\begin{equation}\label{Rdilog}
 {\rm L}(x) = -\frac{\pi^2}{6} -\frac{1}{2} \int_0^x \biggl( \frac{\log(t)}{1-t} + 
\frac{\log(1-t)}{t} \biggr) \ dt \ ,
\end{equation}
\noindent  where we integrate first along the path $[0;1/2]$ on the real axis and then along any path in $\mathfrak{D}$ from $1/2$ to $x$. Here we add $-\pi^2/6$ so that ${\rm L}(1)=0$. It is well-known that L verifies the fundamental Schaeffer's identity:
\begin{equation}\label{Rfivet}
{\rm L}(x) - {\rm L}(y) +{\rm L}(y/x) - {\rm L}(\frac{1-x^{-1}}{1-y^{-1}}) + 
{\rm L}(\frac{1-x}{1-y})=0
\end{equation}
which for real $x$, $y$ holds when $0 < y < x < 1$. In fact, this
identity characterizes the Rogers dilogarithm: if $f(0;1) \to \mr$ is
a 3 times differentiable function satisfying (\ref{Rfivet}) for all $0
<y <x <1$, then $f(x) = k{\rm L}(x)$ for a suitable constant $k$.
By analytic continuation, the
relation (\ref{Rfivet}) also holds true for complex parameters $x$,
$y$ with ${\rm Im}(y) \ne 0$, providing that $x$ lies inside the
triangle formed by $0$, $1$ and $y$.  Note that for such $x$, $y$ all
the arguments of L in (\ref{Rfivet}) have imaginary parts with the
same sign. For every non-degenerate $\Ii$-tetrahedron $(\Delta,b,w)$
set ${\rm L}(\Delta,b,w) = {\rm L}(w_0)$, and for every
$\Ii$-triangulation $\Tt_{\Ii}$ set $\textstyle {\rm L}(\Tt_{\Ii}) =
\sum_i {\rm L}(\Delta_i,b_i,w_i)$.  

\noindent We note that ${\rm L}(\Delta,b,w)$ does not respect
the tetrahedral symmetries. Also, with the above restriction on the
moduli, the Schaeffer's identity implies the invariance of ${\rm
L}(\Tt_{\Ii})$ for the same specific instance of $\Ii$-transit shown
in Fig. \ref{idealt}, and considered in (b) above. On another hand, $\mathfrak{L}_N(\Delta,b,w)$ is a peculiar matrix representation of a specific operator $\Phi$ acting on a suitable completion of the $\mc$-algebra generated by two elements $a$, $b$ satisfying $ab=\zeta ba$ (see \cite{B}). The operator $\Phi$ may be defined by an $N$-dependent power series whose dominant term for $N \rightarrow \infty$ essentially involves dilogarithms. It satisfies a non-commutative version of the relation (\ref{Rfivet}), which induces the basic pentagon identity (\ref{basicpentagon}) in the particular matrix representation defining $\mathfrak{L}_N(\Delta,b,w)$. The dominant term for $N \rightarrow \infty$ of that `quantum Schaeffer's identity' satisfied by $\Phi$ is the exponential of (\ref{Rfivet}) up to a multiplicative constant times $N$, where L is expressed as its power series expansion for $\vert x - 1/2 \vert < 1$ \cite{BR,K3}. It turns out that this result also holds for the matrix entries of $\mathfrak{L}_N(\Delta,b,w)$.

\noindent These facts justify the following name: $\mathfrak{L}_N(\Delta,b,w)$ is the $N^2$-dimensional {\it non symmetric quantum dilogarithm}, computed on the given 
$\Ii$-tetrahedron. 

\noindent In order to
construct invariants for $(W,\rho)$ based on
$\mathfrak{L}_N(\Delta,b,w)$, these should be modified so that the
corresponding modified state sums are invariant with respect to the
{\it whole} set of instances of $\Ii$-transits, as well as the choice
of branching. This is done as follows. Formally similar problems have
been solved in \cite{BB2} to define a dilogarithmic invariant ${\rm
R}(W,\rho)$ based on ${\rm L}(\Delta,b,w)$.

\medskip

\noindent {\bf Symmetrized quantum dilogarithms.} Let $(\Delta,b,w)$ be an $\Ii$-tetrahedron. The notion of integral charges on hyperbolic ideal tetrahedra 
that we are going to define is strictly related to that of {\it flattenings}, introduced by Neumann in his work on 
Cheeger-Chern-Simons classes of hyperbolic manifolds \cite{N1}-\cite{N3}. Flattenings also emerge
straightforwardly in \cite{BB2}, to repare the non-invariance of 
${\rm L}(\Delta,b,w)$ with respect to a change of branching, that we discussed above. In a similar way, the integral charges are 
going to be used in order to (partially) repare the same non-invariance of the 
quantum dilogarithms $\mathfrak{L}_N(\Delta,b,w)$.
The main difference between integral charges and flattenings is that the charges do not
depend on the moduli; a charge defines a flattening on a non-degenerate 
$\Ii$-tetrahedron only if $*_w = -1$.  
\begin{defi}\label{Deltacharge}{\rm An {\it integral
charge} on $(\Delta,b,w)$ is a $\mz$-valued map defined on the edges
of $\Delta$ such that $c(e)=c(e')$ for opposite edges $e$ and $e'$,
and $c_0 + c_1 + c_2 = 1$ (where $c_i=c(e_i)$). We call $c(e)$ the
charge of $e$.}
\end{defi}
\noindent Write $N=2p+1$, and for
each edge $e$ of $\Delta$ set $c'(e)=(p+1)\ c(e)$ mod($N$), viewed as a point in
$\{0,\ldots,N-1\}$. Recall the notation $p_i'$ for the determinations of the $N$-th roots of the $p_i$'s.
\begin{defi}\label{symqd}
{\rm The $N^2$-dimensional 
{\it symmetrized quantum dilogarithm} is the matrix valued
function $\mathfrak{R}_N(\Delta,b,w,c): I_1 \otimes I_2 \rightarrow O_1 \otimes O_2$ defined on the set of charged $\Ii$-tetrahedra
$(\Delta,b,w,c)$ and given by}
\begin{eqnarray}\label{formsym}
\mathfrak{R}_N(\Delta,b,w,c) = \left\lbrace \begin{array}{l}
\bigl( (-p_1'/p_2')^{-c_1}\ (-p_2'/p_0')^{c_0} \bigr)^p\ R'(w \vert c) \quad {\rm if} \quad *=1 \\
\bigl( (-p_1'/p_2')^{-c_1}\ (-p_2'/p_0')^{c_0} \bigr)^p\ \bar{R}'(w \vert c) \quad {\rm if} \quad *=-1 
\end{array} \right.
\end{eqnarray}
\noindent 
{\rm where $*=\pm1$ according to the $b$-orientation of $\Delta$, and
the matrix entries of $R'(w \vert c)$ and $\bar{R}'(w \vert c)$ are respectively}
\begin{eqnarray}\label{symmatrix}
R'(w \vert c)_{\alpha,\beta}^{\gamma,\delta} = \zeta^{c_1'(\gamma - \alpha)} \ 
R(p_1',p_0',-p_2')_{\alpha,\beta - c_0'}^{\gamma - c_0',\delta}\nonumber \hspace{0.2cm}\\ \\
\bar{R}'(w \vert c)_{\gamma,\delta}^{\alpha,\beta} = \zeta^{c_1'(\gamma - \alpha)}\ 
\bar{R}(p_1',p_0',-p_2')_{\gamma+c_0',\delta}^{\alpha,\beta+c_0'}\ .\nonumber
\end{eqnarray}
\end{defi}
\noindent 
As for $\mathfrak{L}_N(\Delta,b,w)$, we see from (\ref{omeg0}) that $\mathfrak{R}_N(\Delta,b,w,c)$ only depends on $(b,w,c)$, and not on the choice of the $N$-th roots $p_i'$ of the $p_i$'s. 

\noindent Write $\nu = g(1)/\vert g(1) \vert$. Let $S$ and $T$ be the $N \times N$ invertible square matrices with matrix entries 
$$T_{m,n} = \nu\ \zeta^{m^2/2}\delta(m+ n)\quad ,\quad S_{m,n} = N^{-1/2}\zeta^{mn}\ .$$
\noindent We have
$$S^4 = id\quad ,\quad S^2 = \zeta'(ST)^3$$ 
\noindent for some root of unity $\zeta'$. Hence the matrices $S$ and $T$ define a projective $N$-dimensional representation $\Theta$ of $SL(2,\mathbb{Z})$. The following lemma describes the tetrahedral symmetries of $\mathfrak{R}_N$ in terms of $\Theta$. Recall that the symmetry group on four elements numbered from $0$ to $3$ is generated by the transpositions $(01)$, $(12)$ and $(23)$. 

\begin{lem} \label{6jsym} 
Let $(\Delta,b,w,c)$ be a charged $\Ii$-tetrahedron with $*_b=+1$. If we change $b$ via the transpositions $(01)$, $(12)$ and $(23)$ of the vertices we have respectively
\begin{eqnarray}
\mathfrak{R}_N\bigl((01)(\Delta,b,w,c)\bigr) & \equiv_{N} & 
\pm \ T_1^{-1}\ \mathfrak{R}_N(\Delta,b,w,c) \ T_1 
\nonumber \\
\mathfrak{R}_N\bigl((12)(\Delta,b,w,c)\bigr) & \equiv_{N} &  
\pm \ S_1^{-1}\ \mathfrak{R}_N(\Delta,b,w,c)\ T_2 
\nonumber \\
\mathfrak{R}_N\bigl((23)(\Delta,b,w,c)\bigr) & \equiv_{N} &  
\pm \ S_2^{-1}\ \mathfrak{R}_N(\Delta,b,w,c)\ S_2  
\nonumber 
\end{eqnarray}
\noindent  
where $\equiv_{N}$ means equality up to multiplication by $N$-th roots of unity. Here we write $T_1=T \otimes 1$, etc.
\end{lem}
 
\noindent {\it Proof.}  Use the relations $w_0w_1w_2=-1$ between the moduli and $c_0 + c_1 +c_2 = 1$ between the charges to rewrite the scalar factors in both sides of each equality in terms of the same variables. For instance for the first equality we have on the left hand side:
$$\left((w_0')^{-1})^{-c_2}((w_2')^{-1})^{c_0}\right)^p =\left((w_0')^{-c_1+1}((w_0'w_2')^{-c_0})\right)^p$$
where $w_0'=-p_1'/p_2'$, $w_1'=-p_2'/p_0'$ and $w_2'=-p_0'/p_1'$. As $w_0'w_1'w_2'=-1$, up to a sign this is equal to $(w_0')^p$ times $((w_0')^{-c_1}((w_1')^{c_0}))^p$. This last scalar is exactly the one appearing on the right hand side. Then the result follows from Prop. \ref{symmetry} in the Appendix. We do the same for the other transpositions. \hfill $\Box$

\medskip

\noindent {\bf Complete pentagon relations.} \noindent Let us say that an $\Ii$-triangulation $\Tt=(T,b,w)$ of $(W,\rho)$ is
{\it roughly charged} if every abstract tetrahedron
$(\Delta_i,b_i,w_i)$ is equipped with an integral charge $c_i$. We say
`roughly' because in Section \ref{LFIXQHI} it shall be necessary to
specialize to integral charges satisfying global constraints. By
replacing in (\ref{statesum}) the non-symmetric quantum dilogarithms
with the symmetrized ones, we obtain new state sums
\begin{equation}\label{newSS}
\mathfrak{R}_N(\Tt_{\Ii},c) = \sum_\alpha \prod_i \mathfrak{R}_N(\Delta_i,b_i,w_i,c_i)_\alpha\quad .
\end{equation} 
The next step is to introduce a
suitable notion of {\it charged $\Ii$-transit}, such that
$\mathfrak{R}_N(\Tt_{\Ii},c)$ is invariant for all instances of
$2\leftrightarrow 3$ charged $\Ii$-transit. As the integral charges do
not depend on the moduli, also a charged $\Ii$-transit is obtained
by completing a usual $\Ii$-transit with a moduli-independent {\it
charge transit}. We use the notations of Def. \ref{id-transit}.

\begin{defi}\label{chargetransit} 
{\rm We say that there is a {\it charge transit} $(T,c)
\leftrightarrow (T',c')$ if $c'$
equals $c$ on the edges of the abstract tetrahedra of $T$ not involved
in the move, and for any other edge $e$ we have the {\it transit of
sum} condition:}
\begin{equation}\label{sumtransitcond}
\sum_{a \in \epsilon_T^{-1}(e)} c(a) = \sum_{a' \in
\epsilon_{T'}^{-1}(e)} c'(a')\quad .
\end{equation}
\end{defi}

\noindent Note that for positive $2 \to 3$
transits this relation implies that the sum of the charges around the new edge
after the move is equal to $2$.

\begin{prop}\label{6jtransit} 
For any charged $2 \leftrightarrow 3$ $\Ii$-transit $(T,b,w,c) \leftrightarrow
(T',b',w',c')$ we have
$$\prod_{\Delta_i \subset T} \ \mathfrak{R}_N(\Delta_i,b_i,w_i,c_i)
\equiv_{N}\pm \ \prod_{\Delta_i' \subset T'}
\ \mathfrak{R}_N(\Delta_i',b_i',w_i',c_i')\quad .$$
\end{prop}

\noindent {\it Proof.} Denote by $f(\Delta,b,w,c)$ the scalar factor in front of the matrices $R'$ and $\bar{R}'$ in (\ref{formsym}). By Prop. \ref{EP} in the Appendix we see that the statement is true if the $\Ii$-transit is the one shown in Fig. \ref{idealt}, and if we remove $f(\Delta_j,b_j,w_j,c_j)$ from $\mathfrak{R}_N(\Delta_j,b_j,w_j,c_j)$, for each tetrahedron $\Delta_j$ involved in the move. We claim that we also have
\begin{equation}\label{sca}
\prod_{\Delta_j \subset T} \ f(\Delta_j,b_j,w_j,c_j) \equiv_{N}\pm \ \prod_{\Delta_j' \subset T'} f(\Delta_j,b_j,w_j,c_j)\quad .
\end{equation}
\noindent Indeed, denote by $c^i$ the integral charge of the tetrahedron opposite to the $i$-th vertex (for the ordering of the vertices induced by the branching), and rewrite the moduli as in Fig. \ref{idealt}. Let log be the standard branch of the logarithm. Up to $N$-th roots of unity the left hand side of (\ref{sca}) is
\begin{eqnarray}
\exp\left( \frac{p}{N}(-c_1^1\log(y)+c_0^1\log(1-y))\right) \nonumber \hspace{7cm}\\
\hspace{2cm}\exp\left( \frac{p}{N}(-c_1^3\log(y(1-x)/x(1-y)) +c_0^3\log((y-x)/x(1-y))\right)\nonumber 
\end{eqnarray}
\noindent and the right hand side is
\begin{eqnarray}
\exp\left( \frac{p}{N}(-c_1^0\log(x)+c_0^0\log(1-x))\right) \ \exp\left( \frac{p}{N}(-c_1^2\log(y/x)+c_0^2\log(1-y/x)) \right) \nonumber \hspace{3cm}\\
 \exp\left( \frac{p}{N}(-c_1^4\log((1-x)/(1-y))+c_0^4\log((x-y)/(1-y))\right)\quad . \nonumber \hspace{2.5cm}
\end{eqnarray} 
\noindent Consider the exponents in these formulas. An elementary computation using the relations (\ref{sumtransitcond}) shows that they are equal up to $2i\pi/N$. For instance the coefficient of ${\rm log}(y)$ in the left hand side is $-c_1^1-c_1^3=-c_1^2$, whereas in the right hand side it is $-c_1^2$. Things go similarly for the coefficients of ${\rm log}(1-y)$, etc. Hence the statement is true for the $\Ii$-transit shown in Fig. \ref{idealt}. We get the result for all the instances of $\Ii$-transit by using Lemma \ref{6jsym}, together with the fact that the action of the matrices $S$ and $T$ cancel on a common face of two tetrahedra (details on this claim are given in the proof of Lemma \ref{branchinv}).\hfill $\Box$ 

\section{Link-fixing and QHI for triples $(W,L,\rho)$ }\label{LFIXQHI}
We first refine the notion of charged $\Ii$-triangulation so as
to make it {\it stable} for charged $\Ii$-transits. A naive
idea would be to require that the sum of the charges around each edge of
$T$ is equal to $2$. But simple combinatorial considerations show that
such tentative global integral charges do not exist. A way to
overcome this difficulty is to fix an arbitrary non empty link $L$ in $W$,
considered up to ambient isotopy, and to incorporate this link-fixing in 
all the constructions. This eventually leads to the definition of the QHI
for triples $(W,L,\rho)$.
\begin{defi}\label{dist}
{\rm A {\it distinguished triangulation} of $(W,L)$ is a pair $(T,H)$
such that $T$ is a triangulation of $W$ and $H$ is a {\it Hamiltonian}
subcomplex of the $1$-skeleton of $T$ which realizes the link $L$
(Hamiltonian means that $H$ contains all the vertices of $T$) .}
\end{defi} 

\begin{defi}\label{DTWL}{\rm  A $\Dd$-{\it triangulation}
$\Tt=(T,H,b,z)$ {\it for a triple}  $(W,L,\rho)$ consists of 
a $\Dd$-triangulation $(T,b,z)$ for $(W,\rho)$ such that 
$(T,H)$ is a distinguished triangulation of $(W,L)$.}
\end{defi}
So a $\Dd$-triangulation for $(W,L,\rho)$ is
a distinguished and quasi-regular triangulation of $(W,L)$, 
decorated by a branching $b$ and a $PSL(2,\mc)$-valued simplicial 
$1$-cocycle $z$.
We postpone to Subsection \ref{DTWLEX} the proof of the existence 
of such $\Dd$-triangulations for any triple $(W,L,\rho)$.

\begin{figure}[ht]
\begin{center}
\includegraphics[width=4cm]{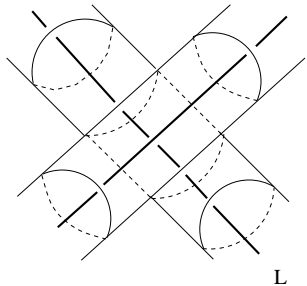}
\caption{\label{tunnel} A tunnel junction over a diagram crossing.} 
\end{center}
\end{figure}

\begin{figure}[ht]
\begin{center}
\includegraphics[width=9cm]{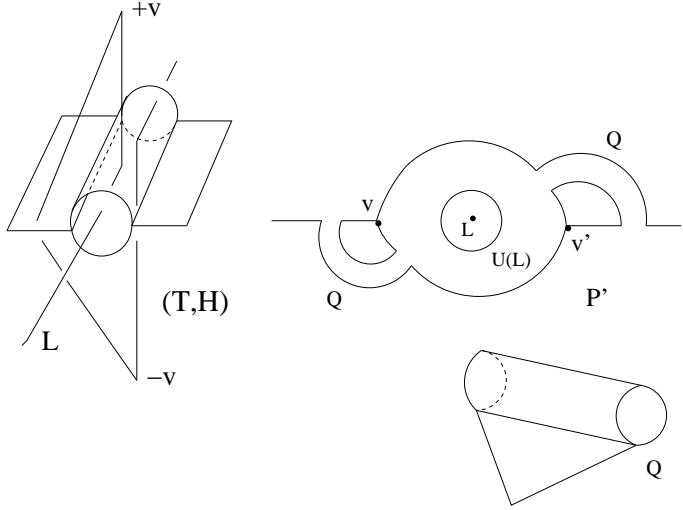}
\caption{\label{variespine} Final steps of the constructions of
$(T,H)$ and $T'$.} 
\end{center}
\end{figure}

\begin{exa} \label{tunnelconst} 
{\rm {\it The tunnel construction.} Here is a simple construction of
distinguished and quasi-regular triangulations of $(S^3,L)$ derived from link
diagrams. 

\noindent Remove two ordered open $3$-balls $B^3_{\pm}$
from $S^3$ away from the link $L$. We get a manifold homeomorphic to
$S^2 \times [-1,1]$ with the embedded $2$-sphere $\Sigma = S^2 \times
\{0\}$ as a simple spine, and two ordered spherical boundary
components $\Sigma_{\pm}$.  Consider a generic projection $\pi(L)$ of
$L \subset S^2 \times [-1,1]$ onto $\Sigma$ such that every connected
component of $\Sigma \setminus \pi(L)$, called a {\it diagram region},
is an open $2$-disk (for instance, this is automatic if $L$ is a
knot). Then, as usual, encode $L$ by a link diagram on $\Sigma$ with
support $\pi(L)$, by specifying the under/over crossings with respect
to the direction normal to $\Sigma$ and going from $\Sigma_-$ towards
$\Sigma_+$. Dig tunnels on $\Sigma$ around $\pi(L)$, by respecting the
under/over crossings, as in Fig. \ref{tunnel}.  Glue $2$-disk walls
inside the tunnels, one between each of the tunnel junctions, such
that their boundaries span meridians. So there is one wall for each
arc of the link diagram. In this way we get a standard spine $P$
corresponding to a quasi-regular triangulation $T$ of $S^3$. To obtain
a distinguished triangulation $(T,H)$ of $(S^3,L)$ do as
follows. There are two distinguished vertices $\pm v$ in $T$, at the
interior of the balls $B^3_{\pm}$ we have initially removed. The edges
of $T$ which are dual to the walls realize $L$ and contain all the
vertices of $T$ except $\pm v$. Select one wall $D$, and remove from
$L$ the interior of the edge dual to $D$. We get an arc with two
vertices of $T$ as endpoints. Connect one of these vertices to $+v$
and the other with $-v$, by means of the edges of $T$ dual to the two
opposite regions contained in the boundary of the tunnel around the
removed edge (see the left side of Fig. \ref{variespine}). Finally
connect $+v$ with $-v$ by another edge dual to an adjacent region of
$P$ contained in $\Sigma$. This construction gives a distinguished and
quasi-regular triangulation of $(S^3,L)$.  Note that we can define a
very particular branching $b$ on $T$ as follows. Fix an orientation of
$L$. Then, the walls are positively oriented in accordance with the
orientations of $L$ and $S^3$, and the other regions of $P$ are
positively oriented with respect to the flow transverse to $P$ and traversing $S^2 \times
[-1,1]$ from $\Sigma_-$ towards $\Sigma_+$.

\noindent Next we show how to modify sligthly the above construction in
order to obtain an ideal triangulation for $Y= S^3 \setminus U(L)$,
where $U(L)$ is an open tubular neighbourhood of $L$. Remove from $P$
all the $2$-disk walls. We get a standard spine $P''$ of $Y' = S^3 \setminus \{
U(L)\cup B_+ \cup B_-\}$. Then modify $P''$ near the (removed) wall $D$ as shown on the right side of Fig. \ref{variespine}. In fact we attach to $P''$ two copies of the $2$-dimensional polyhedron $Q$, and then we remove 4 open $2$-disks at its extremities, on $P''$. The effect is to remove the interior of two
$1$-handles connecting $\partial U(L)$ with $\Sigma_{\pm}$, so that
the so obtained $P'$ is a standard spine of $Y$.

\noindent Note that for both $P$ and $P'$ there is the same pattern
of $4$ vertices at each diagram crossing (Fig. \ref{tunnel}).  It
corresponds to an octahedron of $T$ or $T'$ made of $4$ tetrahedra.
In $P$, there are 2 more vertices for each wall (hence for each arc in
the diagram). In $P'$ there are just 2 further vertices (indicated as
$v$ and $v'$ in Fig. \ref{variespine}). The non-tunnel regions of $P$
which are contained in $\Sigma$ exactly correspond to the original
diagram regions. For both constructions, the diagram arc corresponding to
the selected wall $D$ plays a peculiar
role. Also, the adjacent regions are modified by the respective final
steps. One can obviously orient the regions of $Q$ so that the branching $b$
of $P$ extends to a branching of $P'$.}
\end{exa}
\noindent We have to refine the notion of $\Dd$-transit in order to incorporate
the fixed link $L$. Roughly speaking, a {\it $\Dd$-transit
$(T,H,b,z)\to (T',H',b',z')$ of $\Dd$-triangulations for $(W,L,\rho)$}
consists of a $\Dd$-transit $(T,b,z)\to (T',b',z')$ of
$\Dd$-triangulations for $(W,\rho)$ such that the two Hamiltonian
subcomplexes $H$ and $H'$ which realize $L$ coincide on the tetrahedra
not involved by the underlying move. Precisely:
\smallskip

1)  Any positive $0 \to 2$ or $2 \to 3$ move $T\to T'$ naturally specializes 
to a move $(T,H)\to (T',H')$; in fact $H'=H$ is still Hamiltonian.
The inverse moves are defined in the same way. In particular, for  
negative $3\to 2$ moves we require that the disappearing edge of $T$
belongs to $T \setminus H$; 
\smallskip

2) For positive bubble moves, we assume that an edge $e$ of $H$ lies
in the boundary of the involved face $f$; then $e$ lies in the
boundary of a unique $2$-face $f'$ of $T'$ containing the new vertex
of $T'$. We define the Hamiltonian subcomplex $H'$ of $T'$ just
by replacing $e$ with the other two edges of $f'$. The inverse move is
defined in the same way.
\begin{defi}\label{distmove}{\rm The above moves 
make sense for (non-necessarily quasi-regular) 
distinguished triangulations of $(W,L)$.
We will refer to them as {\it distinguished moves}.}
\end{defi}

\subsection{ Integral charges on $(T,H)$}\label{CHARGTH}
\noindent 
Let $(T,H)$ be a distinguished triangulation of $(W,L)$. Let us recall
the notations already used in Lemma \ref{compat}. We 
denote by $E(T)$ the set of edges of $T$, by $E_{\Delta}(T)$ the whole
set of edges of the associated abstract tetrahedra $\{\Delta_i\}$, and
by $\epsilon_T:E_{\Delta}(T) \longrightarrow E(T)$ the natural
identification map.

\noindent  Let $s$ be a simple closed curve in $W$ in
general position with respect to $T$. We say that $s$ {\it has no
back-tracking} if it never departs a tetrahedron
of $T$ across the same $2$-face by which it entered. Thus each time $s$
passes through a tetrahedron, it selects the edge between the entering
and departing faces.

\begin{defi} \label{defcharges}  
{\rm An {\it integral charge} on a distinguished triangulation $(T,H)$ 
of $(W,L)$ is a map $c: E_{\Delta}(T)
\to \mz$ such that the restriction of $c$ to each abstract tetrahedron 
$\Delta$ of $T$ is an integral charge (see Def. \ref{Deltacharge}), 
and such that the following global properties are satisfied:
\smallskip

\noindent (1) for each $e\in E(T)\setminus E(H)$ we have 
$ \sum_{e'\in \epsilon ^{-1}(e)}c(e') = 2$,

\smallskip

\noindent 
\ \ \ \ \ for each $e\in E(H)$ we have 
$ \sum_{e'\in \epsilon ^{-1}(e)}c(e') = 0$.

\smallskip

\noindent (2) Let $s$ be any curve which has no back-tracking with 
respect to $T$. 
Each time $s$ enters a tetrahedron of $T$ the map $c$ associates an
integer to the selected edge.  Let $c(s)$ be the sum of these
integers. Then, for each $s$ we have $c(s) \equiv 0 \ {\rm mod} \ 2$.

\medskip  

\noindent We call $c(e)$ the \emph{charge} of the edge $e$.}
\end{defi}

\noindent A map $c: E_{\Delta}(T) \to \mz$ inducing a charge on each tetrahedron of $T$ and satisfying Def. \ref{defcharges} (1) defines an 
element $[c] \in H^1(W;\mz/2\mz)$. The meaning of Def. \ref{defcharges}
$(2)$ is that we prescribe $[c]=0$. Note that any integral
charge $c$ on $(T,H)$ eventually encodes $H$, hence the link $L$.

\begin{defi}
{\rm A {\it charged $\Dd$-triangulation} for a triple 
$(W,L,\rho)$ consists of a couple $(\Tt,c)$ where $\Tt=(T,H,b,z)$ 
is $\Dd$-triangulation for $(W,L,\rho)$, and $c$ is an integral charge on
$(T,H)$.}
\end{defi}

\begin{teo}\label{chargeex} For every distinguished triangulation
$(T,H)$ of $(W,L)$ there exist integral charges. In particular, 
every $\Dd$-triangulation $\Tt$ of a triple $(W,L,\rho)$ can be charged.
\end{teo}
This theorem is obtained by adapting, almost
verbatim, Neumann's proof of the existence of combinatorial
flattenings of ideal triangulations of compact $3$-manifolds whose
boundary is a union of tori (Th. 2.4.(i) and Lemma 6.1 of \cite{N1}).
In Neumann's situation there is no link but the manifold has a
non-empty boundary; only the first condition of Def. \ref{defcharges}
(1) is present, and there is a further condition in
Def. \ref{defcharges} (2) about the behaviour of the charges on the
boundary. In our situation, as $W$ is a closed manifold, this further
condition is essentially empty. The second condition in
Def. \ref{defcharges} (1) together with the fact that $H$ is Hamiltonian 
replace the role of the non-empty boundary
in the combinatorial algebraic considerations that lead to the
existence of combinatorial flattenings. All the details of this adaptation are
contained in \cite[Prop. 2.2.5]{B}.

\medskip

\noindent Next we describe the structure of the set of integral charges on $(T,H)$, which is an affine space over an integer lattice. Again, this is an adaptation to the present situation of a result of \cite{N1}.
Let $(T,H)$ be a distinguished triangulation of $(W,L)$, and choose an abtract tetrahedron $\Delta$ of $T$. By definition, there are only two degrees of freedom in choosing the
charges of the edges of $\Delta$.  Assume for
simplicity that 
$T$ is branched, and use the branching to order the edges of $\Delta$ as
in (\ref{order}), from $e_0$ to $e_2'$. Hence, given a
branching on $T$ there is a preferred ordered pair of charges
$(c_1^{\Delta},c_2^{\Delta})=(c(e_0),c(e_1))$ for each abstract tetrahedron $\Delta$.

\smallskip

\noindent Set $w_1^{\Delta} := c_1^{\Delta}$ and $w_2^{\Delta}=-c_2^{\Delta}$.
Let $r_0$ and $r_1$ be respectively the number of vertices and edges
of $T$.  An easy computation with the Euler characteristic shows that
there are exactly $r_1 - r_0$ tetrahedra in $T$. If we order the
tetrahedra of $T$ in a sequence $\{ \Delta^i \}_{i=1,\ldots,r_1-r_0}$,
one can write down an integral charge on $(T,H,b)$ as a vector $c=c(w)
\in \mz^{2(r_1-r_0)}$, with
$$c=(w_1^{\Delta^1},\ldots,w_1^{\Delta^{r_1-r_0}},w_2^{\Delta^1},\ldots,
w_2^{\Delta^{r_1-r_0}})^t\
.$$
\begin{prop} \label{lattice} \emph{\cite[Cor. 2.2.7]{B}}
The difference between any two integral charges $c$ and $c'$ on $(T,H)$ is a linear combination 
with integer coefficients of determined vectors $w(e) \in \mz^{2(r_1-r_0)}$ associated to the edges 
$e \in E(T)$: $\textstyle c'=c+\sum_e \lambda_ew(e)$.
\end{prop}

\noindent The vectors $w(e)$ have the following form. For any abstract
tetrahedron $\Delta^i$ glued along a specific edge $e$, define $r_1^{\Delta^i}(e)$
(resp. $r_2^{\Delta^i}(e)$) as the number of occurences of
$w_1^{\Delta^i}$ (resp. $w_2^{\Delta^i}$) in $\epsilon^{-1}(e) \cap
\Delta^i$. Then
$$w(e)=(r_2^{\Delta^1},\ldots,r_2^{\Delta^{,r_1-r_0}},-r_1^{\Delta^1},
\ldots,-r_1^{\Delta^{r_1-r_0}})^t
\in \mz^{2(r_1-r_0)}\ .$$ 

\begin{figure}[h]
\begin{center}
\includegraphics[width=10cm]{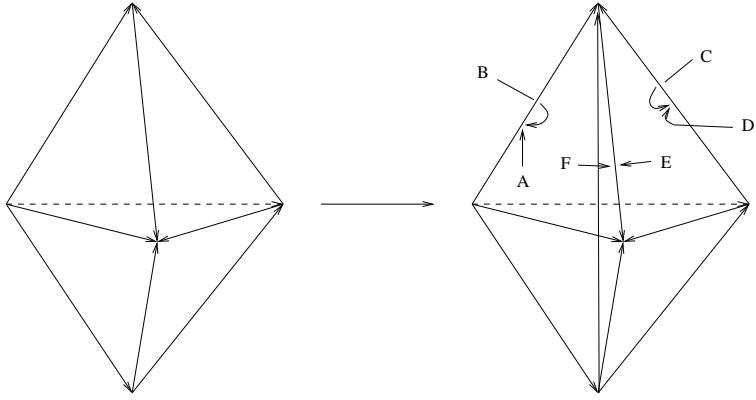}
\caption{\label{w(e)} $2 \to 3$ charge transits 
are generated by Neumann's vectors $w(e)$.}
\end{center}
\end{figure}

\begin{exa} \label{exa1} 
{\rm Consider the situation depicted in the right of
Fig. \ref{w(e)}. Denote by $\Delta^j$ the tetrahedron opposite to the $j$-th vertex. We have
$$\begin{array}{lll}
r_1^{\Delta^0}(e)=-1 & \ \ r_1^{\Delta^2}(e)=0  & \ \ r_1^{\Delta^4}(e)=-1 \\
r_2^{\Delta^0}(e)=1  & \ \ r_2^{\Delta^2}(e)=-1 & \ \ r_2^{\Delta^4}(e)=1
\end{array}$$ 
where $e$ is the central edge. Then $w(e) = (1,-1,1,1,0,1)^t$. 
}
\end{exa}

\noindent {\bf Charge transit.} 
Charge transits for roughly charged triangulations of $(W,\rho)$
have been described in Def. \ref{chargetransit}. 
We have to prove that they specialize well to integral charges on 
$(T,H)$.

\begin{lem} \label{charge-transit} 
Let $(T_1,H_1) \to (T_2,H_2)$ be any
distinguished move between distinguished triangulations of
$(W,L)$. Assume that $c_1$ is an integral charge
on $(T_1,H_1)$, and that $c_1$ transits as a rough charge $c_2$
on $(T_2,H_2)$. Then $c_2$ is actually an integral charge on
$(T_2,H_2)$.
\end{lem}

\begin{defi}\label{chargedDT}{\rm We have a {\it charged $\Dd$-transit} $(\Tt_1,c_1)\to (\Tt_2,c_2)$ between charged $\Dd$-triangulations of a triple $(W,L,\rho)$ if $\Tt_1\to \Tt_2$ is a $\Dd$-transit and $(T_1,H_1,c_1)\to (T_2,H_2,c_2)$ is a transit of integral charges as in Lem. \ref{charge-transit}.}
\end{defi}

\begin{lem}\label{charge-transit2} 
Suppose that  $(T_1,H_1) \to (T_2,H_2)$ is a $2 \to 3$ move between distinguished triangulations of
$(W,L)$. Fix integral charges $c_1$, $c_2$ on $(T_1,H_1)$ and $(T_2,H_2)$ respectively, and put
$$C(e,c_2,T)=\{ c_2' = c_2 + \lambda w(e), \lambda \in \mathbb{Z} \}$$
where $e$ is the edge that appears and $w(e)$ is as in Example
\ref{exa1}. The integral charges $c_2'$ obtained by varying the
charge transit $(T_1,H_1,c_1) \to (T_2,H_2,c_2')$ exactly span $C(e,c_2,T)$.
\end{lem}

\noindent \emph{Proof of \ref{charge-transit}.} First consider the $2 \to 3$ moves. It follows from Def. \ref{chargetransit} that we can restrict our attention
to Star($e,T_2$).  Consider the situation of Fig. \ref{w(e)}, and denote by $\Delta^i$ the tetrahedron opposite to the $i$-th vertex. Let $c^i$ be the
integral charge on $\Delta^i$ and $c_{jk}^i$ the value of $c^i$ on
the edge with vertices $v_j$ and $v_k$. The relation (\ref{sumtransitcond}) implies that the sum of the charges around each of the edges of $T_1 \cap T_2$ stays equal. Moreover it gives:
$$c_{02}^1 + c_{24}^1 + c_{40}^1 = (c_{02}^4 - c_{02}^3) + (c_{24}^0 - 
c_{24}^3) + (c_{04}^2 - c_{04}^3) = c_{13}^4 + c_{13}^0 + c_{13}^2 - 
(c_{02}^3 + c_{24}^3 + c_{40}^3)$$
\noindent 
where in the second equality we use the fact that opposite edges of a
tetrahedron share the same charge.  Since $c_{02}^1 + c_{24}^1 + c_{40}^1 = c_{02}^3 + c_{24}^3 + c_{40}^3 = 1$ we have $c_{13}^4 + c_{13}^0 + c_{13}^2 = 2$.
\noindent Similar computations show that (\ref{sumtransitcond}) forces $c_2$ to induce an integral charge on each abstract tetrahedron of $T_2$. As $H_1$ is not altered by a $2 \to 3$ move, we conclude that
$c_2$ verifies Def. \ref{defcharges} (1).

\noindent Next consider the $0 \to 2$ moves. Any non-branched $0 \to 2$ move
$(T_1,H_1) \rightarrow (T_2,H_2)$ is a composition of $2 \to 3$ and $3
\to 2$ moves \cite{Pi}. In particular, the negative moves in this
composition do not involve the edges of $E(T_1) \cap E(T_2)$. Also,
the integral charges do not depend on branchings. Then our previous conclusion
for $2 \to 3$ charge transits (which obviously holds for $3 \to 2$ ones)
is still true for $0 \to 2$ charge transits. For such a transit
$(T_1,H_1,c_1) \rightarrow (T_2,H_2,c_2)$, denote by $\Delta'$ and
$\Delta''$ the new tetrahedra. It is easy to verify that it is defined
by $s_1 := c_2(\epsilon^{-1}(e) \cap \Delta') + c_2(\epsilon^{-1}(e)
\cap \Delta'') = 0$ for each $e \in E(T_1) \cap E(T_2)$, by $s_2 :=
c_2(\epsilon^{-1}(e_c) \cap \Delta') + c_2(\epsilon^{-1}(e_c)\cap
\Delta'') = 2$ on the new interior edge $e_c$, and by $s_3 :=
c_2(\epsilon^{-1}(e')\cap \Delta') + c_2(\epsilon^{-1}(e'')\cap
\Delta'') = 2$ on the edges $e'$ and $e''$ opposite to $e_c$ in
$\Delta'$ and $\Delta''$ respectively.

\noindent 
Finally consider the bubble moves. Remark that a distinguished bubble
move $(T_1,H_1)\rightarrow (T_2,H_2)$ is \emph{abstractly} obtained
from the final configuration of a $0 \to 2$ move by gluing two more
faces. The resulting face contains the two
new edges of $H_2$. Define a charge transit for a distinguished
bubble move via
the very same formulas as for a $0 \to 2$ move. This makes sense, because
the sum of the charges is equal to $s_1=0$ along each of the two new
edges of $H_2$, to $s_2=2$ along the other interior edge of $\Delta'
\cap \Delta''$, and to $s_3=2$ along the former edge of $H_1$. Hence for bubble charge transits $c_2$ also satisfies Def. \ref{defcharges} (1).

\begin{figure}[h]
\begin{center}
\includegraphics{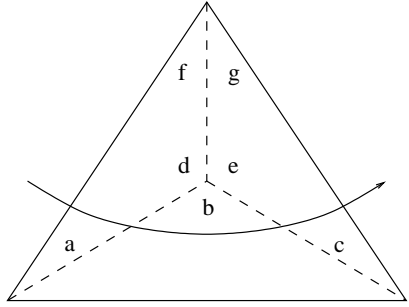}
\caption{\label{triangle} proof of \ref{defcharges} (2) for $c_2$.}
\end{center}
\end{figure}

\noindent 
Let us show that $c_2$ also verifies Def. \ref{defcharges} (2). As
above it is enough to consider a $2 \to 3$ move $(T_1,H_1)
\rightarrow (T_2,H_2)$.  Denote by $e$ the edge that appears. We have
to prove that for any simple closed curve $s$ without back-tracking
with respect to $T_1$ and $T_2$ we have $c_2(s) \equiv 0$
mod($2$). Fig. \ref{triangle} shows an instance of $s$ in a section of
the three tetrahedra of $T_2$ glued along $e$. In this picture the
charges $a,\ldots,g$ are attached to the dihedral angles of the
tetrahedra. Using the first two conditions in Def. \ref{defcharges} for
$c_2$, we see that
$$-a+b-c = (d + f -1) + (2 -
d - e) + (e + g -1) = f + g\quad .\nonumber$$
\noindent Then $c_2(s) = c_1(s) \equiv 0$. \hfill $\Box$

\medskip

\noindent \emph{Proof of \ref{charge-transit2}.} 
Again consider Fig. \ref{w(e)}. The symbols $E,D,F,A,C,B$ denote the charges
on the top edges of $\Delta^0,\Delta^2$ and $\Delta^4$ respectively.
The space of solutions of the linear system (\ref{sumtransitcond}) of relations which define $c_2$ from $c_1$ is one-dimensional. Hence there is a single
degree of freedom in choosing these charges. Fix a particular choice
for them, hence for $c_2$. If $c_2'$ is defined by decreasing $B$ by
$1$, we have
$$c_2'(w) - c_2(w) = (1,-1,1,1,0,1)^t = w(e) \in \mz^{2(r_1-r_0)}.$$
\noindent 
This shows that the integral charges on $T_2$ obtained by varying the
charge transit may only differ by a $\mz$-multiple of $w(e)$.\hfill $\Box$

\subsection{The QHI state sums}\label{QHISS}
We are ready to state the main results of the present paper.

\begin{teo}\label{existenceDWL} For every triple $(W,L,\rho)$ there
exist charged $\Dd$-triangulations $(\Tt,c)$.
\end{teo}

\noindent Fix a triple $(W,L,\rho)$, and let $(\Tt,c)=((T,H,b,z),c)$ be any charged $\Dd$-triangulation of it, with associated charged $\Ii$-triangulation $(\Tt_{\Ii},c)$. Denote by $n_0$ the number of vertices of $T$. Recall the state sums defined in (\ref{newSS}).

\begin{teo}\label{QHIinv} For every odd integer $N>1$, the value of the (normalized) state sum $H_N(\Tt_{\Ii},c) = N^{-n_0}\ \mathfrak{R}_N(\Tt_{\Ii},c)$ does not depend on the choice of $(\Tt,c)$, up to sign and multiplication by $N$-th roots of unity. Hence, up to this ambiguity, it
defines a \emph{quantum hyperbolic invariant} $H_N(W,L,\rho) \in \mc$.
\end{teo}
This shows that $K_N(W,L,\rho)=H_N(W,L,\rho)^{2N}$ is a well-defined complex valued invariant of $(W,L,\rho)$. We can prove immediately the invariance of the QHI state sums $H_N(\Tt_{\Ii},c)$ with respect to the choice of branching and the charged $\Ii$-transits. Recall that Lemma \ref{Iisym} describes how vary the moduli when we change the branching of an $\Ii$-triangulation.

\begin{lem}\label{branchinv} Suppose that $(\Tt',c)$ is obtained from $(\Tt,c)$ be changing the branching. Then 
$$H_N(\Tt_{\Ii},c)\equiv_{N} \pm \ H_N(\Tt'_{\Ii},c)\quad .$$
\end{lem}

\noindent {\it Proof.} Any change of branching on a fixed triangulation translates on each of its abstract tetrahedra $\Delta_i$ as a composition of transpositions of the vertices. By Lemma \ref{6jsym} such transpositions induce an equivariant projective action of $SL(2,\mz)$ on the carrying spaces $I_1 \otimes I_2$ and $O_1 \otimes O_2$ of $\mathfrak{R}_N(\Delta_i,b_i,w_i,c_i)$, which are associated to pairs of faces of $\Delta_i$. This action is defined via matrices $S^{\pm 1}$ and $T^{\pm 1}$. For each (branched) face, it depends on the $b$-orientation of $\Delta_i$: the action is turned into its inverse if we change the agreement between the $b$-orientation of the face and the orientation induced as a boundary of $\Delta_i$. We can see this by simply changing in Lemma \ref{6jsym} the side where the above matrices act. Since a face is always given opposite boundary orientations by the two adjacent tetrahedra, a change of branching may only alter $\mathfrak{R}_N(\Tt_{\Ii},c)$ by the projective factor, which is a sign or an $N$-th root of unity.\hfill $\Box$

\begin{remark}\label{branchnec}{\rm Note that the branching is a necessary ingredient for defining the state sums. Moreover, the branching invariance results from global considerations, as the individual quantum dilogarithms have been only partially symmetrized. This makes a difference, for instance, with respect to the state sums used for the Turaev-Viro invariants.}
\end{remark} 
 
\begin{lem}\label{TINVQHI} Let $(\Tt,c) \to (\Tt',c')$ be any transit of
charged $\Dd$-triangulations for $(W,L,\rho)$. Then
$$H_N(\Tt_{\Ii},c)\equiv_{N} \pm \ H_N(\Tt'_{\Ii},c')\quad .$$
\end{lem}
\noindent {\it Proof.} We use the fact that the $\Dd$-transits dominate $\Ii$-transits (see Prop. \ref{DdomI}). For $2\leftrightarrow 3$ transits, the transit invariance 
of the QHI state sums has been already proved in Prop. \ref{6jtransit}. For the other transits it is obtained as follows. 

\noindent Consider the abstract $2\leftrightarrow 3$ $\Ii$-transit shown in Fig. \ref{idealt}. Denote by $\Delta^i$ the tetrahedron opposite to the $i$-th vertex. Do a further $2 \rightarrow 3$ $\Ii$-transit on $\Delta^0$ and $\Delta^2$. A mirror image of $\Delta^4$ appears, which together with $\Delta^4$ forms the final configuration of a $0 \rightarrow 2$ $\Ii$-transit. Moreover, the other two new $\Ii$-tetrahedra have exactly the same decorations and gluings than $\Delta^1$ and $\Delta^3$. Hence Prop. \ref{6jtransit} implies that, after a trivial simplification, such sequences of $\Ii$-transits (varying the branching and using Lemma \ref{6jsym}) translate as the following {\it orthogonality relations} for the $0 \leftrightarrow 2$ $\Ii$-transits (above for $\Delta^4$): 
$$\mathfrak{R}_N(\Delta,b,w,c)\ \mathfrak{R}_N(\Delta,\bar{b},w,\bar{c}) \equiv_{N}\pm \ id \otimes id\quad .$$
Here $\bar{b}$ and $\bar{c}$ denote the branching and the integral charge mirror to $b$ and $c$, as given by a $0 \to 2$ branched charged move (the explicit formulas for $c$ are given in the proof of Lemma \ref{charge-transit}). The mirror moduli are the same. By taking the trace over one of the tensor factors in the orthogonality relations, we get the {\it normalization relations} corresponding to the bubble $\Ii$-transits. In these relations there is an $N$ in factor; we compensate it by normalizing with $N^{-n_0}$ in $H_N(\Tt_{\Ii},c)$.\hfill $\Box$

\medskip

\noindent The rest of this section shall be mainly devoted to the proof
of Theorems \ref{existenceDWL} and \ref{QHIinv}.

\subsection{Existence of $\Dd$-triangulations for $(W,L,\rho)$}\label{DTWLEX}

\noindent We prove Theorem \ref{existenceDWL}. As the existence of integral charges 
has been already settled, it remains to show the existence 
of $\Dd$-triangulations for any triple $(W,L,\rho)$. 

\noindent Recall Def. \ref{dist} and \ref{DTWL}.
We prove at first the existence of distinguished triangulations for
pairs $(W,L)$. Let us describe these triangulations $(T,H)$ in terms
of dual spines.  Let $M=W \setminus U(L)$, where $U(L)$ is an open
tubular neighbourhood of $L$ in $W$, and $S$ be the union of $t_i\geq
1$ parallel copies on $\partial M$ of the meridian $m_i$ of the
component $L_i$ of $L$, $i=1,\dots,n$. Set $\textstyle r=\sum_i t_i$.

\begin{defi}\label{Ladapted}{\rm  We say that a spine $Q$ of $M$ is 
{\it quasi-standard and adapted} to $L$ 
of {\it type} $t=(t_1,\dots,t_n)$ if:
\smallskip

(i) $Q$ is a simple polyhedron with boundary $\partial Q$ consisting
of $r$ circles. These circles bound (unilaterally) $r$ annular regions
of $Q$. The other regions are cells.
\smallskip

(ii) $(Q,\partial Q)$ is properly embedded in $(M,\partial M)$ and 
transversely intersects $\partial M$ at $S$ (we also say that $Q$ is relative to $S$).
\smallskip

(iii) $Q$ is is a spine of $M$.}
\end{defi}
\noindent Let $Q$ be a spine of $M$ adapted to $L$. Filling each 
boundary component of $Q$ by a $2$-disk we get a standard spine 
$P=P(Q)$ of $W_r = W\setminus r D^3$. 
The dual triangulation $T(P)$ of $W$ contains $L$ as a Hamiltonian 
subcomplex. Conversely, starting from any 
distinguished triangulation $(T,H)$ and removing 
an open disk in each of the regions dual to an edge of $H$, we pass from $P=P(T)$ to a 
quasi-standard spine 
$Q=Q(P)$ of $M$ adapted to $L$. So adapted spines and  distinguished triangulations 
$(T,H)$ are equivalent objects. 

\begin{lem} \label{existadapte} Quasi-standard spines of $M$ adapted to 
$L$ and of arbitrary type,
hence distinguished triangulations of $(W,L)$ 
with an arbitrary number of vertices, do exist.
\end{lem}

\noindent {\it Proof.} Let $\widetilde{P}$ be any standard spine of $M$. 
Consider a {\it normal} retraction $h:M\to \widetilde{P}$. Recall that $M$ is
the mapping cylinder of $h$. For each region $R$ of $\widetilde{P}$,
$h^{-1}(R) = R\times I$; for each edge $e$, $h^{-1}(e) = e\times \{$a
``tripode''$\}$; for each vertev $v$, $h^{-1}(v) = \{$a
``quadripode''$\}$. We can assume that $S$ is in general position
with respect to $h$, so that the mapping cylinder of the restriction
of $h$ to $S$ is a simple spine of $M$ relative to $S$. 
Possibly after doing some $0\to 2$
moves, far from the boundary curves, we obtain a quasi-standard spine $Q$ 
adapted to $L$. \hfill $\Box$
\medskip

\noindent We get the stronger existence result we need with the help of more distinguished moves (see Def. \ref{distmove}).
 
\begin{prop} \label{existqr} For any pair $(W,L)$ there exist distinguished and quasi-regular triangulations. 
\end{prop}
\noindent {\it Proof.} Let $(T,H)$ be any distinguished triangulation of $(W,L)$. It is not quasi-regular 
if some edge $e$ of $T$ is a loop, i.e. if the ends of $e$ are
identified. In the cellulation $D(T)$ of $W$ dual to $T$, this means
that the spine $P=P(T)$ contains some region $R=R(e)$ which has the
same $3$-cell $C$ on both sides: the boundary of $C$ is a sphere $S$
immersed at $R$. Let us say that $R$ is \emph{bad}. We construct a
distinguished and quasi-regular triangulation $(T',H')$ by doing
some distinguished bubble moves on $(T,H)$ (thus adding new $3$-cells to $D(T)$). Then we slide portions of their ``capping'' disks until they cover the bad regions, thus desingularizing all the boundary $2$-spheres.

\noindent Let us formalize this argument. Any (dual) bubble move $P \to P'$ is obtained by gluing a 
closed $2$-disk $D^2$ along its boundary $\partial D^2$, with two transverse intersection points of 
$\partial D^2$ with an edge $e$ of $P$ (see the second move in Fig. \ref{figmove2}). 
Denote by $A$ and $B$, $A \cup B = \partial D^2$, the two arcs thus defined. 
The bubble move is distinguished if at least one of $A$ or $B$ lies on a region $R_H$ of $P$ dual to 
an edge of $H$. The two new regions of $P'$ dual to edges of $H'$ are $D^2$ and the region bounded 
by $\partial D^2$ and adjacent to $R_H$. We call $D^2$ the \emph{capping disk} of the bubble move. 
Note that a bubble move does not increase the number of bad regions, and that any $2 \leftrightarrow 3$ move done by sliding a portion of the capping disk also has this property as long as $\partial D^2$ is embedded.

\begin{figure}[ht]
\begin{center}
\includegraphics[width=9cm]{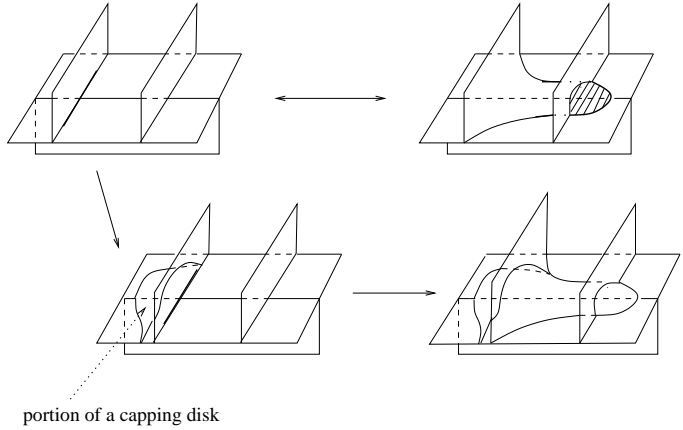}
\caption{\label{badmove} a non-quasi-regular move, 
and how to repair it by capping off the sector of immersion of the corresponding $3$-cell.} 
\end{center}
\end{figure}

\noindent Let now $R \in S$ be a bad region (dashed in the top right of Fig. \ref{badmove}), where $S$ is a singular sphere as above. 
Using distinguished bubble moves we may 
always assume that each connected component of $H$ has at least two vertices. 
Since $(T,H)$ is  distinguished, there are exactly two regions $R_H$ and $R_H'$ in the 
cellular decomposition of $S$ which are dual to edges of $H$. As above, do a bubble move 
that involves $R_H$ (for instance), and slide a portion of its capping disk $D^2$ via $2 \leftrightarrow 3$ 
moves along the $1$-skeleton of $S$, until it reaches a vertex of $R$. 
This is obviously always possible. The only thing is to keep track of the region initially bounded by $\partial D^2$ and adjacent to $R_H$; we cannot remove it, for it is dual to an edge of $H$. Also, if $\partial D^2$ was no longer embedded after this sequence of moves, we could find a shorter sequence leading to the same vertex of $R$. So at each step we still have (dual) distinguished triangulations with no more bad regions. Next expand $D^2$ over $R$ by doing further $2 \leftrightarrow 3$ moves along the edges of $\partial R$, possibly arranged so that they give $0 \leftrightarrow 2$ moves.
If $R$ is embedded in $S$, we can choose such a sequence of moves so that $D^2$ is embedded at 
each step and finally covers $R$ completely (see the bottom right of Fig. \ref{badmove}). 
Both $R$ and $D^2$ are in the boundary of the
$3$-cell introduced by the bubble move. Thus we eventually finish with a spine dual to a distinguished triangulation and having one 
less bad region than $P$. 

\noindent If $R$ is immersed on its boundary (for example if it looks like an annulus
 with one edge that joins the boundary circles), note that it is contained inside a disk embedded in $S$, and as above we may still find a sequence of  $2 \leftrightarrow 3$ moves ending with a spine dual to a distinguished triangulation and having one less bad region than $P$. Iterating this procedure, 
we get the conclusion.\hfill$\Box$

\medskip

\noindent By using Lemma \ref{esistenza} and, for instance, a total ordering branching, we can complete any distinguished and quasi-regular triangulation $(T,H)$ of $(W,L)$ to a $\Dd$-triangulation for $(W,L,\rho)$. So we have achieved the proof of Theorem \ref{existenceDWL}.

\subsection{Invariance of the QHI state sums} \label{QHIINV}
\noindent As bundle preserving oriented homeomorphisms of triples $(W,L,\rho)$ transfer charged $\Dd$-triangulations, 
we can fix a model of $W$ and a flat $PSL(2,\C)$-bundle $\rho$ on $W$, with the link 
$L\subset W$ considered up to ambient isotopy.
 
\begin{figure}[ht]
\begin{center}
\includegraphics[width=9cm]{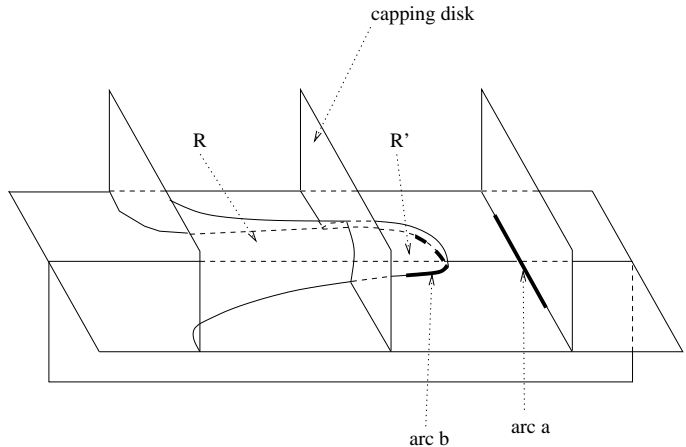}
\caption{\label{noobst} Capping disks are no obstructions for moves.} 
\end{center}
\end{figure}

\noindent We need to show that the set of distinguished and 
quasi-regular triangulations of $(W,L)$ is "connected". In a sense
this is the main technical point. 
As for the existence of $\Dd$-triangulations, let us prove at first a weaker result for  
distinguished triangulations.
Let $(T,H)$ and $(T',H')$ be distinguished triangulations of $(W,L)$ such that the associated
quasi-standard spines $Q$, $Q'$ of $M$ adapted to $L$ are relative to $S$ and $S'$ and are 
of the same type $t$. 
Up to isotopy, we can assume that $S=S'$ and that the ``germs'' of $Q$ and $Q'$ at $S$ coincide. 
By using Th. 3.4.B of \cite{TV} we have the following relative version of Lemma \ref{chemin} for  adapted spines (this follows also from the argument depicted in Fig. \ref{noobst} and used in Prop. \ref{transitqr}):

\begin{lem} \label{chemintype} Let $P$ and $P'$ be quasi-standard spines 
of $M$ adapted to $L$ and relative to $S$. There exists 
a spine $P''$ of $M$ adapted to $L$ and relative to $S$, such that $P''$ can be obtained 
from both $P$ and $P'$ via finite sequences of positive $0 \to 2$ and $2 \to 3$ moves, where at each step the spines are adapted to $L$ and relative to $S$.
\end{lem}

\noindent By possibly using distinguished bubble moves, we deduce from 
Lemma \ref{chemintype} and the correspondence between adapted spines and 
distinguished triangulations that: 

\begin{lem} \label{exist1} Given any two  distinguished 
triangulations $(T,H)$ and $(T',H')$ of $(W,L)$ there exists a distinguished 
triangulation $(T'',H'')$ which may be obtained from both $(T,H)$ and $(T',H')$ via finite sequences of positive bubble, $0 \to 2$ and $2 \to 3$ distinguished moves, where at each step the triangulations of $(W,L)$ are distinguished.
\end{lem}
\noindent Finally we have:
\begin{prop}\label{transitqr} Any two distinguished and quasi-regular 
triangulations $(T,H)$ and $(T',H')$ of $(W,L)$ can be connected 
by means of a finite sequence of distinguished and quasi-regular $2 \to 3$ moves, bubble moves and 
their inverses, where at each step the triangulations of $(W,L)$ are distinguished and quasi-regular. 
\end{prop}

\noindent {\it Proof.}  We use the same terminology as in Prop. \ref{existqr}. 
Let $s: (T,H) \rightarrow \ldots \rightarrow (T',H')$ be a sequence of moves as in Lemma \ref{exist1}.
We may assume, up to further sudivisions of $s$, that there are no $0 \to 2$ moves. 
We divide the proof in two steps. 
We first prove that there exists a sequence $s': P=P(T) \rightarrow \ldots \rightarrow P''$ 
with only quasi-regular moves and such that the spine $P''$ is obtained from $P'=P(T')$ 
by gluing some $2$-disks $\{ D_i^2\}$ along their boundaries. 
Then we show that we may construct $P''$ from $P'$ just by using distinguished bubble moves 
and quasi-regular moves. By combining both sequences we will get the conclusion.

\noindent Bubble moves are always quasi-regular. 
Consider the first non quasi-regular move $m$ in $s$. It produces a bad region $R$; 
see the top of Fig. \ref{badmove}, where we indicate $R$ by dashed lines and we underline the sliding arc $a$. 
Alternatively, a step before $m$ we may do a distinguished bubble move and slide a portion of its capping disk $D^2$ as in 
Prop. \ref{existqr}, until it covers $a$. Next, make the arc $a$ sliding as in $s$; see the bottom of 
Fig. \ref{badmove}. These two moves are quasi-regular and their dual triangulations are distinguished. Starting 
with the moves of $s$ and turning $m$ into this sequence, we define the first part of $s'$. 
We wish to complete it with the following moves of $s$, applying the same procedure each time a 
non quasi-regular move would be done. But suppose that one of these moves would have affected $a$, 
and let $b$ be the sliding arc responsible for it. Then in $s'$ we just have first to slide $b$ ``under'' $D^2$, 
pushing it up. We can do so because all the moves are purely local. This puts $b$ in the same position w.r.t. $a$ than it has in $s$; see Fig. \ref{noobst}. With this rule there are no obstruction to complete the desired sequence $s'$. 
The images in $T''=T(P'')$ of all the capping disks form the set $\{ D_i^2\}$. Remark that there are as many $D_i^2$'s as there were distinguished bubble moves used to construct the sequence $s'$; in other words, the capping disks stay connected all along $s'$. This is due to the fact that in situations such as depicted in Fig. \ref{noobst}, once the region $R$ has bumped into the capping disk the rest of the move is done as in $s$, by sliding the region $R'$. 

\noindent Let us now turn to the second claim. In the dual cellulation $D(T')$
 of $W$ consider the boundary spheres $S_j$ obtained by removing the disks $D_i^2$ one after the other. Fix one of them, $S$, and reversing this procedure let $D^2 \in \{ D_i^2\}$
 (considered with its gluings) be the first disk glued on it. By the above remark, we can do a distinguished bubble move on $S$ and let a portion of its capping disk slide isotopically via $2 \leftrightarrow 3$ moves along the $1$-skeleton of $S$, 
so that it finally reaches the position of $D^2$ in $P''$. We may repeat this argument inductively on the $D_i^2$'s. 
Since all these moves are quasi-regular, this proves our claim.\hfill$\Box$
\medskip

\begin{figure}[h]
\begin{center}
\includegraphics[width=12cm]{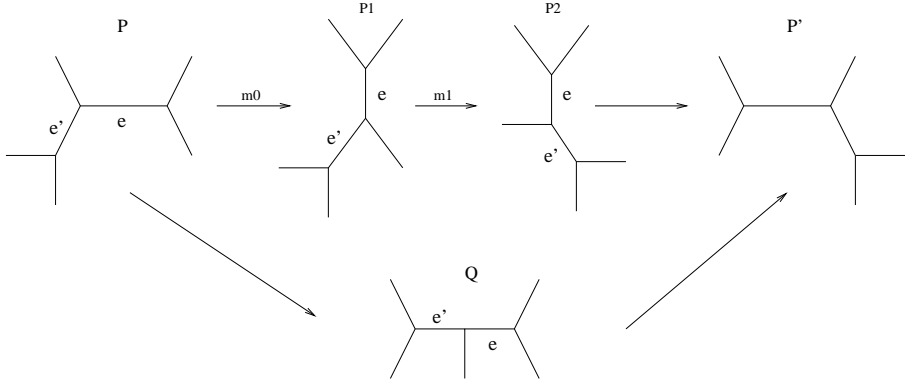}
\caption{\label{2dimtransitqr} the $2$-dimensional analogue of Prop. \ref{transitqr}.}
\end{center}
\end{figure}

\begin{figure}[h]
\begin{center}
\includegraphics[width=9cm]{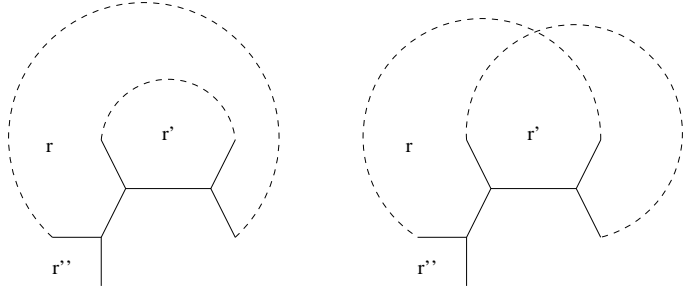}
\caption{\label{2dimcommute} the proof that $r_{P}(m_1)$ is quasi-regular.}
\end{center}
\end{figure}

\noindent We will also use a $2$-dimensional analogue of the previous proposition. 
The main general facts about triangulations and spines of surfaces have been recalled at 
the end of Subsection \ref{trigen}.
\begin{lem}\label{2Dtransitqr} Any two quasi-regular triangulations $T$ and $T'$ of a compact closed
surface $S$ can be connected by a finite sequence of quasi-regular $2\to2$ or $1\to 3$ moves.
\end{lem}
\noindent {\it Proof.} The proof is similar to the one of Prop. \ref{transitqr}. In fact it is simpler as it
uses an argument of commutation of moves which is peculiar to the $2$-dimensional situation.
Consider any sequence $s$ of moves $m_i$ connecting $T$ and $T'$.
View it a sequence 
$$\xymatrix{\relax s:  \ \ldots \ar[r]  & P \ar[r]^{m_0} & P_1  \ar[r]^{m_1} & P_2 \ar[r]^{m_2} & \ldots}$$
\noindent between the ($1$-dimensional) dual spines. On a $1$-dimensional standard spine dual to a quasi-regular 
triangulation of $S$, a move which is not quasi-regular is the flip of an edge that makes it the frontier of a same region. 
Let $m_0$ be the first non quasi-regular flip in $s$, and denote by $e$ the corresponding edge.
A step before $m_0$ let us first apply the ``relative'' $r_P(m_1)$ of $m_1$ on $P$, 
where by ``relative'' we mean the flip of the same edge $e'$; we get $Q$. 
Then apply $r_{Q}(m_0)$; see the bottom sequence of Fig. \ref{2dimtransitqr}. (Beware that in this figure, the notations for $e$ and $e'$ are interchanged when following the upper or the lower sequence of flips; this is why we introduce the notion of `relative'). Note that $r_Q(m_0)$ is necessarily quasi-regular, for otherwise $m_0$ would not be the first 
non quasi-regular flip in $s$, since the horizontal edge below $e'$ in the top left picture of Fig. \ref{2dimtransitqr} would have the same region on both sides. We claim that $r_P(m_1)$ is also quasi-regular. 
Indeed, in $P$ we necessarily have one of the two situations of Fig. \ref{2dimcommute}, where 
the dotted arcs represent boundary edges. In the first situation, $r'=r''$ is impossible. In the second one, 
if $r'=r''$ then $r' = r$ and $m_0$ is not the first non quasi-regular flip in $s$, thus giving a contradiction. 
Hence the sequence $r_{Q}(m_0) \circ r_P(m_1)$ is quasi-regular. 
Moreover we have:
$$P' = r_{P_2}(m_0) \circ m_1 \circ m_0\ (P) =r_{Q}(m_0) \circ r_P(m_1)\left( P \right)\quad .$$

\noindent This implies that we can modify $s$ locally so as to obtain
$$\xymatrix@!{\relax s': \ \ldots \ar[r] & P \ar[r]^{r_P(m_1)} & Q  
\ar[r]^{r_{Q}(m_0)} & P' \ar[r]^{r_{P'}(m_0)} & P_2 \ar[r]^{m_2} & \ldots}$$
\noindent where the first possible non quasi-regular move 
is $r_{P'}(m_0)$. The length of $s'$ after $r_{P'}(m_0)$ is less than the length of $s$ after $m_0$. 
Then, working by induction on the length, replacing each non quasi-regular flip as above and noting 
that $1 \to 3$ moves are always quasi-regular, we get a quasi-regular sequence $s'$. \hfill $\Box$
\medskip

\noindent {\bf Full invariance of the QHI state sums}
Let $(T,H)$ and $(T',H')$ be two arbitrary distinguished and quasi-regular triangulations of $(W,L)$. 
Let $(T,H)\to \dots \to (T',H')$ be a finite sequence of distinguished and quasi-regular moves
which connects $(T,H)$ to $(T',H')$, as in Prop. \ref{transitqr}.
Any total ordering branching $b$ on $T$ 
(see Subsection \ref{branchings}) transits through total ordering branchings to a branching $b'$ on $T'$.
By Lemma \ref{charge-transit}, any integral charge $c$ on $(T,H)$ transits to an integral charge $c'$ on $(T',H')$. Applying Lemma \ref{generique}, we know that for generic $1$-cocycles $z$ on $(T,b)$ these transits can be completed to a sequence a charged $\Dd$-transits which connects    	
the charged $\Dd$-triangulation $(\Tt,c)=(T,H,b,z,c)$ to another $(\Tt',c')=(T',H',b',z',c')$. So, by using the transit invariance of Prop. \ref{TINVQHI}, 
we have proved:
\begin{lem}\label{partialinv} For any triple $(W,L,\rho)$ and every odd integer $N>1$, given two arbitrary distinguished and quasi-regular triangulations $(T,H)$ and $(T',H')$
of $(W,L)$, there exist charged $\Dd$-triangulations $(\Tt,c)$ and $(\Tt',c')$ for $(W,L,\rho)$, supported by $(T,H)$ and $(T',H')$ respectively, such that
$$H_N(\Tt_{\Ii},c)\equiv_{N} \pm \ H_N(\Tt'_{\Ii},c')\quad .$$
\end{lem}

\noindent This statement can be complemented as follows.
\begin{lem}\label{zOK} Assume that $(\Tt,c)$ and $(\Tt',c')$ are charged $\Dd$-triangulations for 
$(W,L,\rho)$ which are connected by a finite sequence of $\Dd$-transits, with the possible exception of some bad cocycle transits for which 
the idealizability condition is lost. Nevertheless we have
$$H_N(\Tt_{\Ii},c)\equiv_{N}\pm \ H_N(\Tt'_{\Ii},c')\quad .$$
\end{lem}
\noindent{\it Proof.} Thanks again to Lemma \ref{generique}, we can replace $z$ and $z'$ with 
{\it arbitrarily close} $1$-cocycles $z_1$ and $z_2$ respectively, such that the corresponding new 
charged $\Dd$-triangulations
$(\Tt'',c)$ and $(\Tt''',c')$ for $(W,L,\rho)$ are actually connected by charged $\Dd$-transits. Then $H_N(\Tt_{\Ii}'',c)\equiv_{N}\pm \ H_N(\Tt_{\Ii}''',c')\ $. Since $z_1$ and $z_2$ are 
arbitrarily close to $z$ and $z'$, and $H_N$ is continuous as a function of idealizable $1$-cocycles, we get the required conclusion. \hfill$\Box$

\medskip

\noindent In the rest of the proof we will tacitely use this genericity/continuity argument, so that we can always assume that the idealizability condition is never lost. So, in order to complete the proof of Theorem \ref{QHIinv}, it is enough to prove the following proposition.

\begin{prop}\label{partialinv2} For any triple $(W,L,\rho)$ and every odd integer $N>1$, given two charged $\Dd$-triangulations  
$(\Tt,c)$ and $(\Tt',c')$ of $(W,L,\rho)$ which only differ by the
respective decorations of a same distinguished and quasi-regular
triangulation $(T,H)$ of $(W,L)$, we have
$$H_N(\Tt_{\Ii},c)\equiv_{N} \pm \ H_N(\Tt'_{\Ii},c')\quad .$$
\end{prop}

\noindent {\bf Proof.} The invariance with respect to the choice
of branching has been already obtained in Lemma \ref{branchinv}. So, from now on, we will use only
total ordering branchings as they do not pose any problems of transit.
\smallskip

\noindent{\bf Charge invariance.}  
Let us localize the problem. Fix a triple $(W,L,\rho)$, a $\Dd$-triangulation $(\Tt,c)=(T,H,b,z,c)$ of $(W,L,\rho)$, and an arbitrary edge $e$ of $T$. 
Consider the set of integral charges which differ from $c$ only on Star$(e,T)$. It is of the form (we use the notations of Prop. \ref{lattice})

$$C(e,c,T) = \{ c'= c+ \lambda w(e),\ \ \lambda \in \mz \}\quad .$$

\noindent Thanks to Prop. \ref{lattice}, it is enough to prove that 
$H_N(\Tt_{\Ii},c)\equiv_{\zeta^{\Z/2}} H_N(\Tt'_{\Ii},c')$ when $c'$ varies in $C(e,c,T)$. Assume that $e\in T\setminus H$; the charge invariance is a consequence of the following facts:

 \smallskip

(1) Let $(\mathcal{ \Tt},c) \to (\mathcal{ \Tt}'',c'')$ be any $2\to 3$ 
charged $\Dd$-transit such that $e$ is a common edge of $T$ and $T''$. Then the result holds 
for $C(e,c,T)$ if and only if it holds for $C(e,c'',T'')$.
\medskip

(2) There exists a sequence of distinguished quasi-regular $2\to3$ moves  
which connects $(T,H)$ to $(T'',H'')$, such that $e$ persists at each step and Star$(e,T'')$ is like the final configuration of 
a $2\to 3$ move, with 
$e$ playing the role of the central common edge of the 3 tetrahedra.

\medskip

(3) If Star$(e,T)$ is like Star$(e,T'')$ in (2), then the result holds for $C(e,c,T)$.

\medskip

\noindent  By Lemmas \ref{charge-transit} and \ref{charge-transit2} we know that $C(e,c,T)$ transits to $C(e,c'',T'')$. As the value of the QHI state sums is 
not altered by charged $\Dd$-transits (Lemma \ref{TINVQHI}), the fact (1) follows.

\smallskip

\begin{figure}[h]
\begin{center}
\includegraphics[width=12cm]{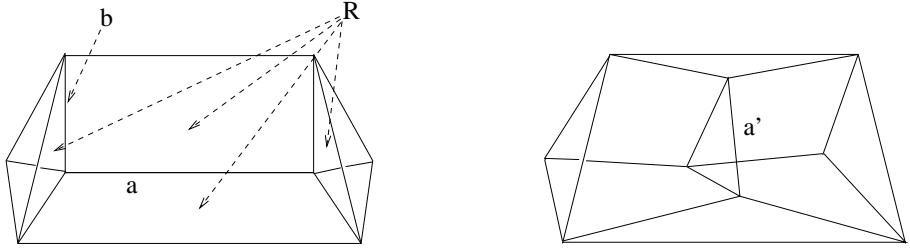}
\end{center}

\caption{\label{embed} evolution of the geometric multiplicity of $R$ when blowing-up $a$.}
\end{figure}

\noindent To prove (2) it is perhaps easier to think, for a while, 
in dual terms. Consider the dual region $R=R(e)$ in $P=P(T)$. The final
configuration of $e$ in $T''$ corresponds dually to the case when $R$ is an embedded
{\it triangle}. More generally, there is a natural notion of {\it geometric multiplicity} $m(R,a)$ of $R$ at 
each edge $a$ of $P$, and $m(R,a)\in \{0,1,2,3\}$. We say that $R$ is {\it embedded} in $P$ if for each 
$a$, $m(R,a)\in \{0,1 \}$. Call {\it proper} an edge with two distinct vertices. If $R$ has a loop in its boundary, a suitable $2\to 3$ move at a proper edge of 
$P(T)$ having a common vertex with the loop puts proper edges in place of the loop. Each time $R$ has 
a proper edge $a$ with $m(R,a) \in \{ 2,3 \}$, the (non-branched) $2\to 3$ 
move along $a$ puts new edges $a'$ with $m(R,a') \leq 2$ in place of $a$. In the situation where this is an equality,
 remark that if we first blow up an edge $b$ adjacent to $a$ and such that $m(R,b) = 2$, and then
 we apply the $2 \to 3$ move along $a$, we get $m(R,a') = 1$ (look at Fig. \ref{embed}). 
By induction, up to $2\to 3$ moves, we can assume that $R$ is an embedded polygon. 
To obtain the final configuration of $e$ in $T''$ let us come back to the dual situation. 
We possibly have more than 3 tetrahedra around $e$. It is not hard to reduce the number to 3, 
via some further $2\to 3$ moves. In the above construction we could accidentely do some non 
quasi-regular moves, which we would like to avoid. For this, do appropriate distinguished bubble moves and slide portions of their 
capping disks as in Prop. \ref{transitqr}. This is always possible because 
these moves may not increase the geometric multiplicity of the edges of the region $R$ under consideration. In this way we eventually find sequences of distinguished and quasi-regular moves which transform $R$ into an embedded triangle.

\noindent Concerning fact (3), do first a $3 \to 2$ $\Dd$-transit on 
$e$ and then a $2 \to 3$ $\Dd$-transit, varying the charge transit $(T,c)\to (T',c'')$. 
By Lemma \ref{charge-transit2} we know that the charges $c''$ exactly describe $C(e,c',T')$. 
Since the value of the state sums is not altered by $\Dd$-transits, this concludes (note that, as we are using total ordering branchings, there is no problem of transit with the negative moves).
\smallskip

\noindent Suppose now that $e \in H$. The analogue of (1) for distinguished bubble 
moves is true for the same reasons. Then, applying a distinguished bubble move on a face of $T$
containing $e$ we are brought back to the previous situation. The charge invariance is thus proved. 
\medskip 

\noindent{\bf Cocycle-invariance.} Let $(\Tt,c)$ and $(\Tt',c)$ be
two charged $\Dd$-triangulations of $(W,L,\rho)$ which only
differ by the $1$-cocycles $z$ and $z'$ representing $\rho$. 
We have to prove that $H_N(\Tt_{\Ii},c)\equiv_{N} \pm\ H_N(\Tt'_{\Ii},c)$.
The two cocycles $z$ and $z'$ differ by a coboundary $\delta \lambda$, and it is enough to consider the elementary case 
when the $0$-cochain $\lambda$ is supported by one vertex $v_0$ of $T$. Again we have localized the problem. The invariance of the value of the state sums for bubble 
$\Dd$-transits gives us the 
result in the special situation when $v_0$ is the new vertex after the move. Let us reduce 
the general case to this special one by means of $\Dd$-transits. We use the notations
and the facts stated at the end of Section \ref{trigen}.
It is enough to show that we can modify Star$(v_0,T)$ to reach the star-configuration 
of the special situation. Recall that Star$(v_0,T)$ is the cone over $S={\rm Link}(v_0,T)$, 
which is homeomorphic to $S^2$. So Star$(v_0,T)$ is determined by the triangulation $T_{v_0}$ of 
$S$. By Lemma \ref{2Dtransitqr} we know that $T_{v_0}$ is connected to the
triangulation of $S$ corresponding to the special situation by a sequence of quasi-regular
$2\to 2$ or $1\to 3$ moves. These can be obtained as trace of quasi-regular $2\to 3$ moves,
by applying inductively the last remark in Section \ref{trigen}. Hence also the cocycle invariance is proved, and finally the proof of Theorem \ref{QHIinv}
is complete. \hfill$\Box$

\subsection{Complements on the QHI}\label{qhicomp}

\noindent {\bf State sums over non quasi-regular triangulations.} Let $\Tt=(T,H,b,z)$ be any branched distinguished (not necessarily quasi-regular) triangulation of a pair $(W,L)$, with an idealizable cocycle $z$ representing a bundle $\rho$. As $T$ is not necessarily quasi-regular, the existence of such a $z$ depends on $\rho$. For instance, if $\rho=\rho_0$ is the trivial flat bundle, it implies that $T$ is quasi-regular. We know that $\Tt$ can be charged, by $c$ say, so the state sum $H_N(\Tt_{\Ii},c)$ is still defined. We claim that in fact
$$H_N(W,L,\rho) \equiv_{N}\pm \ H_N(\Tt_{\Ii},c)\quad .$$
\noindent Indeed, the proof of Prop. \ref{existqr} shows that any distinguished triangulation of $(W,L)$ can be made quasi-regular just by using suitable distinguished bubble moves together with $2 \leftrightarrow 3$ moves done by sliding portions of their capping disks. We can complete such a sequence of moves with arbitrary charge and branching transits; the branchings may not be induced by total orderings on the vertices, so we use the general definition \ref{btransit}. Since $H_N(\Tt_{\Ii},c)$ is invariant for $\Ii$-transits (Lemma \ref{TINVQHI}), we are left to prove that we can complete the above transits with idealizable cocycle transits starting from $z$. For that, remark that the cocycle transits are generically idealizable for bubble moves. Moreover, there is only a finite set of cocycle values on the capping disks that lead to non idealizable cocycle transits for the moves to follow.

\noindent Combining this argument with those used in the proof of Prop. \ref{transitqr} (essentially Fig. \ref{noobst}) and \ref{partialinv2}, with some work we get the following proposition. Although we do not need it for proving Th. \ref{TINVQHI}, it shows that we can bypass the genericity argument of Lemma \ref{zOK}. It is necessary for proving the existence of scissors congruence classes for $(W,L,\rho)$ in \cite{BB2}. Note that it holds in greater generality, replacing $PSL(2,\mc)$ with any algebraic group $G$, and the idealizability condition by demanding that the cocycles take their values outside of some proper algebraic subvarieties of $G$. 

\begin{prop} \label{fullt} Any two $\Dd$-triangulations of a same triple $(W,L,\rho)$ may be connected by a sequence of $\Dd$-transits.
\end{prop}

\noindent The fact that $H_N(W,L,\rho) \equiv_{N} \pm\ H_N(\Tt_{\Ii},c)$ for decorated triangulations which are not quasi-regular but support idealizable $1$-cocycles representing $\rho$ is of practical interest. Indeed, explicit computations are easier with non quasi-regular triangulations, since they contain in general much lesser tetrahedra than quasi-regular ones. 

\medskip

\noindent {\bf Duality.} There are two natural involutions on the arguments of a triple $(W,L,\rho)$: the first consists in changing the orientation of $W$, and the second is defined by passing from $\rho$ to the complex conjugate bundle. The QHI {\it duality}
property relates these involutions. Let $(\Tt,c)$ be a charged
$\Dd$-triangulation for $(W,L,\rho)$. Denote by $z^*$ the complex
conjugate of the $1$-cocycle $z$ of $\Tt$, and by $(\Tt^*,c)$ the corresponding charged $\Dd$-triangulation for
$(W,L,\rho^*)$, where $\rho^* =[z^*]$. We write $-W $ for the manifold $W$ with the opposite orientation. Recall the notation $\equiv_{2N}$ from Lemma \ref{6jsym}. 
\begin{prop}\label{duality} 
We have $(H_N(W,L,\rho))^* \equiv_{N} \pm \ H_N(- W,L,\rho^*)$.
\end{prop}
\noindent {\it Proof.} If we change the orientation of $W$, the $b$-orientation of each tetrahedron $\Delta_i$ turns into the opposite, so that the pairs of faces associated to the carrying spaces $I_1 \otimes I_2$ and $O_1 \otimes O_2$ of $\mathfrak{R}_N(\Delta_i,b_i,w_i,c_i)$ are exchanged. Hence $\mathfrak{R}_N(\Delta_i,b_i,w_i,c_i)$ becomes $^T\mathfrak{R}_N(\Delta_i,\bar{b}_i,w_i,c_i)$, where $^T$ is the transposition of matrices and $\bar{b}_i$ is the branching $b_i$  for the opposite ambient orientation. But Prop. \ref{unitarite} in the Appendix implies that
$$^T\mathfrak{R}_N(\Delta_i,\bar{b}_i,w_i,c_i)_\alpha = \left(\mathfrak{R}_N(\Delta_i,b_i,w_i^*,c_i)_{-\alpha}\right)^*\quad .$$
\noindent Here $\alpha$ is a state, as defined in Section \ref{QDILOGSS}. Since $H_N(\Tt_{\Ii},c)$ does not depend on the states, this yields the conclusion. \hfill $\Box$
\medskip

\noindent {\bf Some natural triples $(W,L,\rho)$.} 1) Let $M$ be a cusped complete hyperbolic $3$-manifold of finite volume. We know from Thurston's hyperbolic Dehn fillings theorem (see e.g. \cite[Ch. E]{BP1}) that there exist sequences $(W_n,L_n,\rho_n)$ of compact hyperbolic Dehn fillings
of $M$ converging geometrically to $M$. Here, $L_n$ denotes the link made of the simple short geodesics in $W_n$ forming the {\it cores} of the fillings, and $\rho_n$ is the holonomy of the hyperbolic manifold $W_n$. These are our favorite examples of triples $(W,L,\rho)$.
\smallskip

\noindent  2) Consider a compact oriented $3$-manifold $M$ with non empty boundary made of tori. Fix $\alpha \in H^1(M;\C)$ and consider the associated flat bundle $\rho_\alpha$ on $M$,
as in Subsection \ref{flatbungen}. The class $\alpha$, hence the holonomy of $\rho_\alpha$, may be non 
trivial on the boundary of $M$. Here is an elementary procedure reminescent of 
the hyperbolic Dehn surgery, which allows us to define triples $(W,L,\rho)$ from these pairs $(M,\alpha)$.
\noindent To simplify the notations, let us assume that $Z = \partial M$ consists of one torus.
It is well-known that the kernel of the map 
$$i_*: H_1(Z; \mq)\to  H_1(M; \mq)$$ is a Lagrangian subspace $\Ll$ of $H_1(Z; \mq)$ w.r.t. 
the intersection form.  Then there exists a basis $(m,l)$ of $\pi_1(Z)\cong H_1(Z;\mz)$ such that
$\Ll$ is generated by the homology class of $pm+ql$, where $p,q \in \mz$ and ${\rm gcd}(p,q)=1$.
Let us denote by $W$ the closed manifold obtained from $M$ by the Dehn filling of $Z$ with 
coefficient $(p,q)$ w.r.t. the basis $(m,l)$. The bundle $\rho_\alpha$ extends to the whole of $W$. If $L$ denotes the core of the filling, 
then $(W, L,\rho_\alpha)$ is a triple canonically associated to $(M,\rho_\alpha)$. 

\noindent For example, if $L$ is a knot in $S^3$ there are two families of QHI that give natural
topological invariants of the knot. The first one is $K_N(S^3,L,\rho_0)=H_N^{2N}(S^3,L,\rho_0)$, where $\rho_0$ is
the trivial flat bundle on $S^3$. The second one is obtained by applying the above
procedure to $M = S^3 \setminus U(L)$ and a generator $\alpha$ of $H^1(M;\Z)\cong \Z$, where $U(L)$ is an open tubular neighbourhood of $L$. Similar considerations apply to links in $S^3$, or
more generally in $\Z$-homology spheres.
\smallskip

\noindent 3) Finally, note that we can specialize the choice of the link. For example, 
we may take $L$ as the trivial knot embedded in an open ball of $W$. In this way we formally
obtain QHI for pairs $(W,\rho)$.
\medskip

\noindent Here are some further remarks.
\begin{remark}\label{phase}{\rm 
{\it About the QHI phase factor.} We have prudently defined $H_N(W,L,\rho)$ only up to sign and multiplication by $N$-th roots of unity, which depend on the branching and the charge of the $\Ii$-triangulations used to compute it. This is due to Lemma \ref{6jsym} and Prop. \ref{symmetry} in the Appendix. It is natural to ask wether this phase ambiguity is in fact not present, due to some systematic global compensations between the roots of unity coming from each tetrahedron, for a given change of branching on an $\Ii$-triangulation.

\noindent Alternatively, it is known that branchings and suitably 
restricted sets of branching transits can be used to encode several 
extra-structures on $3$-manifolds, such as combings, framings, spin and Euler structures \cite{BP2, BP3}. So we wonder about the existence of a suitable extra-structure on the pair $(W,L)$ which, in our setup, would reflect itself in the branchings, and could serve to dominate the phase ambiguity.  The models we have in mind are the Euler structures on $W$ for which $L$ is a pseudo-Legendrian link. As Turaev discovered, the Euler structures dominate the ambiguity, due to the action of the fundamental group on the universal covering, in the definition of Reidemeister torsions (see \cite{T} and also \cite{BP4}).}
\end{remark} 

\begin{remark} \label{Bsym} {\rm {\it On the $B$-QHI.} We already considered in \cite{BB1} the QHI restricted to $B$-characters. In that paper we used state sum formulas differing from those in
Th. \ref{QHIinv} by a scalar factor depending on the cocycle $z$ of the 
$\Dd$-triangulation $\Tt$ (not only on the associated $\Ii$-triangulation
$\Tt_{\Ii}$).  This was a consequence of a slightly different
symmetrization procedure of the quantum dilogarithms, which consisted in replacing in (\ref{formsym}) the scalar factor in front of the matrices $R'$ and $\bar{R}'$ by $(-q_2')^p$. (The $q_j$'s have been defined in Remark
\ref{idealrem} (3), and $'$ denotes the determinations
of the $N$-th roots of the $q_j$'s induced by a common determination of the $N$-th roots of the cocycle values). Let us write $\mathfrak{R}_N^B(\Tt,c)$ for the associated state sums. 

\noindent 
Then, the statement of Lemma \ref{6jsym} is unchanged, except that
the ambiguity is only up to $N$-th roots of unity. However, in Prop. \ref{6jtransit} we have to multiply both sides by the respective $\textstyle Q_2 := \prod_i
(-q_2')_i^p$.  It is a remarkable but somewhat fortuitous fact that,
for $B$-characters and for any positive $2 \to 3$ $\Dd$-transit $\Tt
\rightarrow \Tt'$, we have $Q_2(\Tt')/Q_2(\Tt)=x(e)^{2p}$, where
$x(e)$ is the upper-diagonal value of the cocycle $z$ on the new edge
in $T' \setminus H'$.  Normalizing $\mathfrak{R}_N^B(\Tt,c)$ by dividing
it with $\textstyle \prod_{e \in T \setminus H} x(e)^{2p}$, we
eventually get a well defined invariant $H_N^B(W,L,\rho)$ up to $N$-th roots of
unity. The same procedure for general $PSL(2,\mc)$-characters (using
the $p_2'$'s instead of the $q_2'$'s) does not seem to work, because
the explicit formula for $P_2(\Tt')/P_2(\Tt)$ heavily
depends on the branching. Moreover, we believe that it is conceptually
relevant that the QHI for arbitrary $PSL(2,\C)$-characters can be
computed only in terms of the idealized $\Ii$-triangulations
$\Tt_{\Ii}$.}
\end{remark}

\section{On the Volume Conjecture} \label{VOLCONJ}

\noindent Denote by $J_N(L)$ the coloured Jones polynomial of the link $L$ in
$S^3$, with colour $N$ on each component of $L$, normalized by dividing
it with the value on the unknot, and evaluated at
$\zeta=\exp(2i\pi/N)$. By combining the results of \cite{K2} and
\cite{MM} we know that
\begin{teo}\label{KMM}
For every link $L$ in $S^3$ we have $J_N(L)\equiv_{N}H_N^B(S^3,L,\rho_0)$,
where $\rho_0$ is the necessarily trivial character of $S^3$, and $H_N^B$ is the QHI for $B$-characters discussed in Remark \ref{Bsym}.
\end{teo}
\noindent By using Th. \ref{KMM} we can state the {\it Volume Conjecture} of Kashaev \cite{K4} as:
\begin{conj}\label{Jconj}
For every hyperbolic link $L$ in $S^3$ we have
$$ \lim_{N\to \infty} (2\pi/N) \log ( |J_N(L)|) = {\rm Vol}(M)$$
where $M$ is the cusped complete hyperbolic manifold (unique up to isometry) 
homeomorphic to the complement of $L$ in $S^3$.
\end{conj}
\noindent 
Recall that this conjecture has been rigorously confirmed at
least for the celebrated figure-8 knot (see the references in \cite{Y}). In this Section we try to set Conjecture \ref{Jconj} against the background of the general QHI theory we have developed, also in order to find a geometric motivation for it. Our leading idea is

\smallskip

\noindent {\it The hyperbolic geometry is a constitutive element of the QHI,
because they are defined as state sums over the hyperbolic ideal
tetrahedra of any $\Ii$-triangulation. So their asymptotic behaviour
should be expressable in terms of suitable 'classical' invariants of
hyperbolic nature, computable over the same $\Ii$-triangulations
and sharing with the QHI some basic structural features.}

\smallskip

\noindent This idea cannot be implemented straightforwardly. Indeed, in the case of $(S^3,L)$, the hyperbolic geometry associated to the trivial character $\rho_0$ of $S^3$ by the idealization is trivial. On the other hand, Th. \ref{KMM} shows that $H_N(S^3,L,\rho_0)$ actually reflects the non trivial geometry of $S^3\setminus L$. In the general case (for instance when $W$ is hyperbolic and $\rho$ is its holonomy) we expect that $H_N(W,L,\rho)$ combines, in a not yet understood way, the non trivial contributions coming from both $W \setminus L$ and $(W,\rho)$. For $S^3\setminus L$, we can still implement our leading idea, as follows.

\medskip

\noindent {\bf QHI for cusped $3$-manifolds.} 
The technology we have developed in this paper can be applied to the
hyperbolic manifold $M=S^3\setminus L$, and more generally to any
non-compact complete hyperbolic $3$-manifold $M$ of finite volume. Let
us call it a {\it cusped} manifold.

\noindent Consider a geometric triangulation of $M$ by geodesically embedded
ideal tetrahedra of non-negative volume. It is well-know that such
triangulations do exist \cite{EP}. 
The manifold $M$ is homeomorphic to the interior of a
compact manifold $Y$ with non-empty boundary made of tori, and the
above triangulation, forgetting the hyperbolic structure, is a
topological ideal triangulation of $Y$ in the sense of Section
\ref{trigen}. Assume that this triangulation admits a branching
$b$. This is a rather mild assumption.  The hyperbolic ideal
tetrahedra can be encoded as usual by the cross-ratio moduli. This gives
an $\Ii$-triangulation $\Tt_{\Ii}$ of $M$ with possibly some (but
not all) degenerate tetrahedra, such that for each non-degenerate
$(\Delta_j,b_j,w_j)$ of $\Tt_{\Ii}$ we have $*_j = *_{w_j}$.  We can
endow $\Tt_{\Ii}$ with an integral charge $c$ as in \cite{N1} (see the
discussion after Prop. \ref{chargeex}). So the formula (\ref{newSS})
defines a state sum $\mathfrak{R}_N(\Tt_{\Ii},c)$.
Prop. \ref{branchinv} and the statement in Prop. \ref{TINVQHI}
concerning the $2\to3$ and $0\to 2$ transits do apply to these state
sums.  

\noindent In spite of these facts, there are some technical problems to
prove that $\mathfrak{R}_N(\Tt_{\Ii},c)$ defines
an invariant $H_N(M)$. For instance, it was important in the proof of
Th. \ref{QHIinv} that the $\Ii$-transits were dominated by
$\Dd$-transits. On the other hand, it may happen (as for an hyperbolic
knot in $S^3$) that the ideal triangulation of $Y$ only admits the
trivial constant $1$-cocycle, which is not idealizable. Anyway, let us
postulate here that $H_N(M)$ is well defined; the details about its
construction and invariance shall be worked out in a paper in
preparation. Alternatively, the reader can replace $H_N(M)$ with
$\mathfrak{R}_N(\Tt_{\Ii},c)$ without effecting seriously the rest of
the discussion.

\noindent For every cusped manifold $M$, set
\begin{equation}\label{VCS}
{\rm R}(M) := {\rm CS}(M) + i\ {\rm Vol}(M)\quad {\rm
mod}(\pi^2\Z)
\end{equation}
where ${\rm CS}(M)$ and ${\rm Vol}(M)$ are respectively the metric
Chern-Simons invariant and the hyperbolic volume of the cusped
manifold $M$. We propose the following generalization of Conjecture
\ref{Jconj}, that gives it a strong geometric motivation.
\begin{conj}\label{volconj} 
(1) For every cusped manifold $M$, there exist $C \in \mc^*$ and $D
\in \mc$ such that
$$ H_N(M)^{2N} = \left[ CN^{D}\exp\left( \frac{N\ {\rm
R}(M)}{i\pi}\right)\left(1 + \mathcal{O}(1/N)\right) \right]^{2N}\ .$$
\noindent (2) If $L$ is a hyperbolic link in $S^3$ and $M=S^3\setminus L$, then
$$H_N(S^3,L,\rho_0) \equiv_{N} \pm \ H_N(M) \ .$$
\end{conj}

\noindent Clearly, both assertions are interesting on their own. We can relax the second, still in a meaningful way, by stating the equality up to a different normalization of $H_N(S^3,L,\rho_0)$, or even that it holds only asymptotically. Note that point (1) implies
\begin{equation}\label{volH}
\lim_{N\to \infty} (2\pi/N) \log \left(| H_N(M)|\right) = {\rm
Vol}(M) \quad .
\end{equation}
\noindent Together with point (2) this generalizes Conjecture
\ref{Jconj}, because the QHI for $B$-characters have the same asymptotic behaviour than those for $PSL(2,\mc)$-characters. The conjecture \ref{volconj} says at first that $H_N(M)^{2N}$ has an asymptotic power series expansion with, in general, an exponential
growth rate. Assuming it, the invariance of $H_N(M)^{2N}$ and the
uniqueness of coefficients of asymptotic expansions imply that
$\exp({\rm R}(M)/i\pi)$, $C$ and $D$ are well-determined invariants of
$M$.  Then, it predicts that ${\rm R}(M)$ is of the form
(\ref{VCS}). We have expressed the conjecture in terms of the $2N$-th
power of $H_N(M)$ so as to kill an eventual multiplicative ambiguity up
to $2N$-th roots of unity (which is present in $H_N(W,L,\rho)$). Point (2) would make manifest the hyperbolic geometry of $M$ hidden in $H_N(S^3,L)$.

\medskip

\noindent Let $M$ be a cusped manifold and $(W_n,L_n,\rho_n)$ be a sequence of
compact hyperbolic Dehn fillings of $M$ converging geometrically to
$M$, thanks to Thurston's hyperbolic Dehn filling theorem.  Here,
$L_n$ denotes the link made of the simple short geodesics in $W_n$
forming the cores of the fillings, and $\rho_n$ is the holonomy of the
hyperbolic manifold $W_n$.  Recall that ${\rm Vol}(W_n) \to {\rm Vol}(M)$ when
$n\to \infty$. We also propose:
\begin{conj}\label{2volconj}
For every fixed $N$, when $n\to \infty$ we have 
$$ H_N(W_n,L_n,\rho_n)^{2N} \longrightarrow H_N(M)^{2N} \quad .$$
\end{conj}
\noindent By taking a double limit, this and (\ref{volH}) imply
that, when $n,\ N \to \infty$, we have
$$ (2\pi/N) \log \left(|H_N(W_n,L_n,\rho_n)|\right) \longrightarrow
{\rm Vol}(M) \quad .$$

\medskip

\noindent {\bf Motivations and comments.} (1) Set ${\rm R}(W,\rho) := {\rm
CS}(\rho) + i\ {\rm Vol}(\rho)$ mod$(\pi^2\Z)$, where ${\rm
CS}(\rho)$ and ${\rm Vol}(\rho)$ are respectively the Chern-Simons invariant and the volume of the character $\rho$ (see \cite{D2} and the references therein for these notions). For every pair $(W,\rho)$, we have proved in \cite{BB2} that $\exp((1/i\pi){\rm R}(W,\rho))$ has strong
structural relations with the QHI. For instance, as ${\rm
R}(-W,\rho)=-{\rm R}(W,\rho)$, ${\rm CS}(\rho^*)= {\rm CS}(\rho)$ and
${\rm Vol}(\rho^*)= -{\rm Vol}(\rho)$, we see that $\exp((1/i\pi){\rm
R}(W,\rho))$ formally verifies the duality property stated in
Prop. \ref{duality}. More substantially, ${\rm
R}(W,\rho)$ can be computed over any $\Ii$-triangulation $\Tt_{\Ii}$
for $(W,\rho)$ endowed with a so-called ``flattening'' $f$ as
\begin{equation}\label{Rogstate}
\textstyle {\rm R}(W,\rho) = {\rm R}(\Tt_{\Ii},f) = \sum_j *_j{\rm
R}(\Delta_j,b_j,w_j,f_j)
\end{equation}
where the sum runs over the branched hyperbolic ideal tetrahedra of
$\Tt_{\Ii}$ with induced flattenings $f_i$, $*_j$ is the index of the
branching $b_j$, and ${\rm R}(\Delta,b,w,f)$ is a suitably
`uniformized' and symmetrized version of the Rogers dilogarithmic
function ${\rm L}(\Delta,b,w)$, defined in Section \ref{QDILOGSS}. So
$\exp((1/i\pi){\rm R}(W,\rho))=
\exp((1/i\pi){\rm R}(\Tt_{\Ii},f))$ looks very like a QHI state sum
$H_N(\Tt_{\Ii},c)$ (here it should be with $N=1$). This formula
refines a description {\rm mod}($\pi^2\Q$) of the universal second
Cheeger-Chern-Simons class on $BPSL(2,\mc)$ due to Dupont-Sah
\cite{DS,D1}, and is in agreement with the results of \cite{N1} and
\cite{NY}, stated for cusped and
closed hyperbolic $3$-manifolds (in the particular case when $\rho$ is
their holonomy).

\noindent The symmetrized quantum dilogarithms $\mathfrak{R}_N(\Delta,b,w,c)$
and the symmetrized Rogers dilogarithm ${\rm R}(\Delta,b,w,f)$ verify
the same fundamental identities, that is they are invariant for all
instances of charged (resp. flattened) $\Ii$-transits. Moreover, as
mentioned in Section \ref{QDILOGSS}, the Rogers dilogarithm (also the
symmetrized one) is the unique solution of these functional
identities, up to a multiplicative scalar factor. Finally, the
classical dilogarithms play the main role in the asymptotic expansion
of the quantum dilogarithms, whence of the QHI.

\noindent On another hand, the construction of the QHI includes a
link-fixing while the
one of ${\rm R}(W,\rho)$ is link-free.  This corresponds to the fact
that the integral charges do not depend on the cross-ratio moduli, in
contrast with the flattenings.  This is a crucial difference because
we know that the QHI are sensitive to the link, even asymptotically.
However, this
discrepancy vanishes when we work with cusped $3$-manifolds, so that
Conjecture \ref{volconj} (1) looks as an appropriate implementation of
the leading idea stated at the beginning.

\smallskip

\noindent (2) The presence of the link $L$ in $H_N(W,L,\rho)$ as well as its
ambiguity up to sign and multiplication by $N$-th roots of unity are entirely
a consequence of the specific symmetrization procedure of the basic
state sums $\mathfrak{L}_N$ for $(W,\rho)$, that we have adopted in
Section \ref{QDILOGSS}. Suitable variations of this
procedure based on {\it moduli-dependent charges}, similar to the
flattenings, should allow us to define the QHI directly for $(W,\rho)$.
The asymptotic behaviour of such ``absolute'' 
QHI should be dominated by $R(W,\rho)$, similarly to Conjecture \ref{volconj} (1). 

\smallskip

\noindent (3) Here we outline a possible way to approach Conjecture \ref{volconj} (2). We can use the triangulations $(T,H)$ of
$(S^3,L)$ and $T'$ of $Y = S^3 \setminus U(L)$ constructed in
Example \ref{tunnelconst} to compute $H_N(S^3,L,\rho_0)$ and $H_N(M)$ respectively.
In both cases we have a complete decoration including an appropriate
integral charge, and cross-ratio moduli of the involved $\Ii$-tetrahedra. In the first case we use as usual the idealization of an idealizable cocycle
representing the trivial character $\rho_0$. In the second case we
assume that the moduli are obtained via a sequence of $\Ii$-moves
connecting $T'$ with an hyperbolic geodesic triangulation of $M$.
Recall that both constructions of $(T,H)$ and $T'$ include the
selection of a same link-diagram arc, hence the selection of a
$(1,1)$-tangle presentation of $L$. Then, developing the contributions
of the diagram crossings to the state sums, we obtain for
$H_N(S^3,L,\rho_0)$ and $H_N(M)$ very close expressions in terms of suitable
$R$-matrices {\it depending on parameters}, and supported by that
$(1,1)$-tangle presentation of $L$. But the values
of the parameters of each $R$-matrix are specified by the respective
global decorations (the charges give ``discrete''
parameters, and the cross-ratio moduli ``continuous'' ones).

\noindent On another hand, to compute $J_N(L)$ we can use {\it bare} tangle
presentations of $L$, and, as shown in \cite{MM}, a single {\it constant} Kashaev's $R$-matrix which corresponds to one \emph{fixed} particular choice in the
parameters. The proof of Theorem \ref{KMM} includes a reduction of the
above expression for $H_N^B(S^3,L,\rho_0)$ to an expression which involves only that constant $R$-matrix. This is due
to Kashaev and is not a trivial fact. The main ingredients are indicated in \cite{K2}.\footnote{The second author thanks Kashaev for having
explained him the details.}

\noindent So Conjecture \ref{volconj} (2) would be achieved if the same
reduction to the constant $R$-matrix holds also for the formally
similar non-constant $R$-matrix expressions of $H_N(S^3,L,\rho_0)$ and $H_N(M)$.
This cannot be a simple adaptation of the $H_N^B(S^3,L,\rho_0)$ case,
because the {\it global homological} properties of the integral charges as well as the fact that the moduli satisfy both edge compatibility and boundary completeness necessarily enter the proof. We believe that even eventually
disproving this reduction should be very instructive.

\section{Appendix: quantum dilogarithms}\label{APP}

\noindent In this Section we present the definition of the $N^2 \times N^2$-matrix valued quantum dilogarithms as matrices of $6j$-symbols, which is originally due to Kashaev \cite{K1}. We also state their fundamental functional/symmetry relations needed for the present paper. We refer to \cite[Ch. 3]{B} for details and for the proofs. 

\medskip

\noindent Recall that $\zeta = \exp(2i\pi/N)$ and that $N > 1$ is an odd positive integer. Set $N= 2p+1,\ p \in \mathbb{N}$. We shall henceforth denote $1/2 := p + 1$ mod($N$). Fix the determination $\zeta^{1/2} = \zeta^{p+1}=-\exp(i\pi/N)$ of the square root of $\zeta$.

\bigskip

\noindent {\bf Cyclic representations of $\mathcal{B}_{\zeta}$.} Consider the $\mc$-algebra $\mathcal{B}_{\zeta}$ with unity generated by elements $E$, $E^{-1}$ and $D$ such that $ED = \zeta DE$. It is well-known that $\mathcal{B}_{\zeta}$ can be endowed with a structure of Hopf algebra isomorphic to the simply-connected (non-restricted) integral form of a Borel subalgebra of $U_q(sl(2,\mathbb{C}))$ specialized in $q=\zeta$ \cite[\S 9]{CP}. Thus it has the following co-multiplication, co-unit and antipode maps :
$$\Delta(E)  = E \otimes E\ ,\quad \Delta(D) =  
E \otimes D + D \otimes 1$$
$$\epsilon(E)=1\ ,\quad \epsilon(D)=0\ ,\quad S(E)  = E^{-1}\ ,\quad S(D) = -E^{-1}D\quad .$$
\noindent Given a representation $\rho$ of $\mathcal{B}_{\zeta}$, denote by $V_\rho$ the associated $\mathcal{B}_{\zeta}$-module. It is easily seen that if $\rho$ is irreducible, then ${\rm dim}_{\mc}(V_\rho)\leq N$. We say that $\rho$ is \emph{cyclic} if $\rho(D) \in {\rm GL}(V_{\rho})$, i.e. if ${\rm dim}_{\mc}(V_\rho)=N$. Recall that the tensor product of two representations $\rho$ and $\mu$ is defined by 
\begin{equation}\label{tensrep}
\textstyle (\rho \otimes \mu) (a) = \sum_i \rho(a_i') \otimes \mu (a_i'')
\end{equation}
where $a \in \mathcal{B}_{\zeta}$, $\textstyle \Delta(a) = \sum_i a_i' \otimes a_i''$, and the tensor product on $V_{\rho} \otimes V_{\mu}$ is over $\mc$. We say that a sequence $\rho_1,\ldots,\rho_n$ of irreducible cyclic representations of $\mathcal{B}_{\zeta}$ is \emph{regular} if $\rho_i \otimes \ldots \otimes \rho_{i+j}$ is cyclic, for any $1 \leq i \leq n,\  1 \leq j \leq n -i$. Two representations $\rho$ and $\mu$ are {\it equivalent} if there exists an isomorphism $V_\rho \rightarrow V_\mu$ commuting with the action of $\mathcal{B}_{\zeta}$. 

\noindent The algebra $\mathcal{B}_{\zeta}$ is a free module of rank $N$ over its center $\mathcal{Z}$, which is generated by $E^{\pm N}$ and $D^N$. The elements of $\mathcal{Z}$ act as scalar operators on any $\mathcal{B}_{\zeta}$-module $V_{\rho}$, so they define homomorphisms $\chi_\rho: \mathcal{Z} \rightarrow \mc$ called the {\it central characters}. Put $e_{\rho}=\chi_{\rho}(E^N)$ and $d_\rho=\chi_{\rho}(D^N)$. The following lemma is an easy exercise:

\begin{lem}\label{parametrization} Two irreducible cyclic representations $\rho$ and $\mu$ of $\mathcal{B}_{\zeta}$ are equivalent iff $(e_{\rho},d_{\rho})=(e_{\mu},d_{\mu}) \in \mc^* \times \mc^*$.
\end{lem}  

\noindent  We find a nice parametrization of these equivalence classes $[\rho]$ by rewriting $e_{\rho}$ and $d_{\rho}$ as follows. Given non zero complex numbers $t_\rho$ and $x_\rho$ we define a {\it standard} (cyclic) representation $\rho$ of $\mathcal{B}_{\zeta}$ by
\begin{equation}\label{stmat}
\rho(E) = t_\rho^{2}Z\ ,\quad \rho(D) = t_\rho x_\rho X
\end{equation}
\noindent where $X$ and $Z$ are the $N \times N$ matrices with components $X_{ij} = \delta_{i,j+1}$ and $Z_{ij} = \zeta^i \delta_{i,j}$ in the standard basis of $\mc^N$, and $\delta_{i,j}$ is the Kronecker symbol. By Lemma \ref{parametrization} any cyclic irreducible representation of $\mathcal{B}_{\zeta}$ is equivalent to a standard one, and two standard representations $\rho$ and $\mu$ are equivalent iff $t_\rho^{2N}=t_\mu^{2N}$ and $t_\rho^N x_\rho^N= t_\mu^N x_\mu^N$. 

\medskip

\noindent For a regular pair $(\rho,\mu)$, the space $V_\rho \otimes V_\mu$ necessarily splits as the direct sum of $N$ cyclic simple $\mathcal{B}_{\zeta}$-modules. Their central characters are given by $e_{\rho \otimes \mu}$ and $d_{\rho \otimes \mu}$. Then, Lemma \ref{parametrization} implies that these submodules are all isomorphic. We call them the {\it product} submodules, and, abusing of notations, we denote them by $V_{\rho\mu}$. A direct sum decomposition of $V_\rho \otimes V_\mu$ into product submodules is obtained by choosing a linear basis of a characteristic subspace 
$$E_i = {\rm Ker}\left((\rho \otimes \mu)(E) - \zeta^i e_{\rho \otimes \mu}'{\rm id}_{V_\rho \otimes V_\mu}\right)$$
where $e_{\rho \otimes \mu}'$ is some $N$-th root of $e_{\rho \otimes \mu}$. The $\mathcal{B}_{\zeta}$-orbit of any element of that basis is a product submodule. If $\rho$ and $\mu$ are standard we can do these choices in a natural way, by using the standard tensor product basis of $V_\rho$ and $V_\mu$. Now, (\ref{tensrep}) gives $e_{\rho \otimes \mu} = e_\rho e_\mu$ and $d_{\rho \otimes \mu} = e_\rho d_\mu + d_\rho$. For the standard product submodules this reads
\begin{eqnarray} \label{char}
t_{\rho \mu}^{2N} = t_\rho^{2N} t_\mu^{2N} \nonumber \hspace{0.8cm}\\
x_{\rho \mu}^N = t_\rho^N x_\mu^N + x_\rho^N/t_\mu^N \nonumber \quad .
\end{eqnarray}
\noindent So, we conclude that the matrices 
\begin{equation}\label{mat01}
\Psi([\rho])=\left( \begin{array}{ll} 
                t_\rho^N & x_\rho^N \\
                0     & t_\rho^{-N}
                \end{array} \right)
\end{equation}
define a one-to-one correspondence $\Psi$ between the equivalence classes of irreducible cyclic representations of $\mathcal{B}_{\zeta}$, and the set of non diagonal upper triangular matrices of $PSL(2,\mc)= SL(2,\mc)/\{\pm I\}$. (The sign ambiguity is due to the choice of square root of $t_{\rho}^{2N}$). Note that this set is open and dense in the quotient matrix topology of the upper Borel subgoup $B$ of $PSL(2,\mc)$. Moreover, a remarkable feature of the parametrization $\Psi$ is that for any regular pair $(\rho,\mu)$ we have $\Psi([\rho])\Psi([\mu]) = \Psi([\rho\mu])$.

\medskip

\noindent {\bf $6j$-symbols.} We are mainly concerned with the monoidal structure of the spaces of embeddings of cyclic simple $\mathcal{B}_{\zeta}$-modules. We define the \emph{multiplicity module} of two irreducible cyclic representations $\rho$ and $\mu$ as the complex vector space of equivariant maps from $V_{\rho}$ to $V_{\mu}$:
$$M_{\rho,\mu} = {\rm End}_{\mathcal{B}_{\zeta}} \left( V_{\rho},V_\mu \right) = \{ U : V_\rho \rightarrow  V_\mu \vert \ U\rho(a) = \mu(a)U,\ \forall \  a \in \mathcal{B}_{\zeta} \}\quad . $$
\noindent We have seen above that for any regular pair $(\rho,\mu)$, we have ${\rm dim}_\mc(M_{\nu,\rho \otimes \mu})=N$ if $[\nu]=[\rho \mu]$, and zero otherwise. Given a regular triple $(\rho,\mu,\nu)$, consider product representations $\rho\mu$, $\mu\nu$ and $\rho\mu\nu$. Set 
$$\begin{array}{l}
M_{\rho,(\mu,\nu )} = {\rm End}_{\mathcal{B}_{\zeta}} \left(V_{\rho \mu \nu}, V_\rho \otimes \left( V_\mu \otimes V_\nu \right)\right) \\ \\ 
M_{(\rho,\mu),\nu} = {\rm End}_{\mathcal{B}_{\zeta}} \left( V_{\rho \mu \nu},\left( V_{\rho} \otimes V_\mu \right) \otimes V_\nu \right)\quad .
\end{array}$$
\noindent We have vector space isomorphisms
$$\begin{array}{l}
M_{\rho , (\mu ,\nu)} \cong  M_{\rho \mu \nu , \rho \otimes \mu \nu}\otimes  M_{\mu \nu , \mu \otimes \nu}\\ \\
M_{(\rho ,\mu), \nu} \cong   M_{\rho \mu  , \rho \otimes \mu} \otimes M_{\rho \mu \nu , \rho \mu \otimes \nu}\quad .
\end{array}$$
\noindent Moreover, the isomorphism of $\mathcal{B}_{\zeta}$-modules 
$$\alpha_{\rho,\mu,\nu}:V_\rho \otimes \left( V_\mu \otimes V_\nu \right) \longrightarrow \left( V_\rho \otimes V_\mu \right) \otimes V_\nu $$
induces a vector space isomorphism between $M_{\rho , (\mu , \nu)}$ and $M_{(\rho ,\mu ), \nu}$. So we eventually get a linear isomorphism 
$$R(\rho,\mu,\nu) :  M_{\rho \mu \nu , \rho \otimes \mu \nu}\otimes  M_{\mu \nu , \mu \otimes \nu} \longrightarrow M_{\rho \mu  , \rho \otimes \mu} \otimes M_{\rho \mu \nu , \rho \mu \otimes \nu} \quad .$$
\noindent The coherence of the isomorphisms $\alpha_{.,.,.}$ for the tensor product of {\it four} cyclic representations making a regular sequence $(\rho,\mu,\nu,\upsilon)$ implies that 
\begin{equation}\label{basicpentagon}
R_{12}(\rho,\mu,\nu)\ R_{13}(\rho,\mu\nu,\upsilon)\ R_{23}(\mu,\nu,\upsilon) = R_{23}(\rho\mu,\nu,\upsilon)\ R_{12}(\rho,\mu,\nu\upsilon)
\end{equation}
where $R_{12} = R \otimes id$ etc. This $3$-cocycloid relation is called the {\it basic pentagon identity}. We can define $R(\rho,\mu,\nu)$ in another equivalent way. Let $\{ K_\alpha(\rho,\mu)\}_{\alpha=1,\ldots,N}$ denote a linear basis of $M_{\rho\mu, \rho \otimes \mu}$, and similarly for the other multiplicity modules. The families of maps $\{(id \otimes K_\delta(\mu,\nu)) \circ K_\gamma(\rho,\mu \nu)\}_{\delta,\gamma}$ and $\{(K_\alpha(\rho,\mu) \otimes id) \circ K_\beta(\rho\mu,\nu)\}_{\alpha,\beta}$ form two distinct linear basis of the space of embeddings of $V_{\rho\mu\nu}$ into $V_\rho \otimes V_\mu \otimes V_\nu$. Then, the isomorphism $R(\rho,\mu,\nu)$ may be realized as the corresponding change-of-basis matrix: 
\begin{eqnarray}\label{6jdefop}
K_\alpha(\rho,\mu) \ K_\beta(\rho\mu,\nu) = 
\sum_{\delta,\gamma = 0}^{N-1} R(\rho,\mu,\nu)_{\alpha,\beta}
^{\gamma,\delta}\ K_\delta(\mu,\nu) \ K_\gamma(\rho,\mu\nu)\quad .
\end{eqnarray}
\noindent The matrix entries $R(\rho,\mu,\nu)_{\alpha,\beta}^{\gamma,\delta}$ are called $6j$-\emph{symbols}, and the basis vectors $K_\alpha(\rho,\mu)$ are \emph{Clebsch-Gordan operators}. The relation (\ref{6jdefop}) translates the coherence of the isomorphisms $\alpha_{.,.,.}$ cited above. In particular, one may prove (\ref{basicpentagon}) by applying both sides to a suitable composition of Clebsch-Gordan operators, and then using (\ref{6jdefop}) several times.

\medskip

\noindent Let us give a standardized form of the Clebsch-Gordan operators for all multiplicity modules. For that, we restrict to standard representations. By definition, each $K_\alpha(\rho,\mu)$ satisfies $(\rho \otimes \mu)(a)K_\alpha(\rho,\mu) = K_\alpha(\rho,\mu)\rho\mu(a)$, for any $a \in \mathcal{B}_{\zeta}$. These equations are polynomials in the parameters of $\rho$, $\mu$ and $\rho\mu$. So, using the parametrization $\Psi$ defined in (\ref{mat01}), we see that $K_\alpha(\rho,\mu)$ is a matrix valued rational function on a branch of an $N$-fold ramified covering of $B \times B \times B$. Here $B$ is the upper Borel subgroup of $PSL(2,\mc)$. More precisely, a direct computation gives the following result. Recall from (\ref{omeg0}) the definition of the function $\omega(x,y,z\vert n)$.
\begin{lem} \label{CG} Let $(\rho,\mu)$ be a regular pair of standard representations of $\mathcal{B}_{\zeta}$. The set of matrices $\{K_\alpha(\rho,\mu)\}_{\alpha = 0,\ldots,N-1}$ with components
\begin{eqnarray} K_\alpha(\rho,\mu)_{i,j}^k = \zeta^{\alpha j+\frac{\alpha^2}{2}}\ \omega(t_\rho x_\mu,x_\rho/t_\mu,x_{\rho \mu}\vert i-\alpha) \ \delta(i+j-k) \nonumber \end{eqnarray}
\noindent form a linear basis of $M_{\rho \mu, \rho \otimes \mu}$.
\end{lem}
\noindent Put $[x] = N^{-1}(1-x^N)/(1-x)$. Recall from Section \ref{QDILOGSS} the definition of the complex valued functions $g$ and $h$. We have:

\begin{prop} \label{6j} In the normalized basis of Clebsch-Gordan operators formed by the matrices $h(x_{\rho\mu}/t_\rho x_\mu) K_\alpha(\rho,\mu)$, the $6j$-symbols read 
$$R(\rho,\mu,\nu)_{\alpha,\beta}^{\gamma,\delta} = h_{\rho,\mu,\nu} \ 
\zeta^{\alpha\delta+\frac{\alpha^2}{2}}\ \omega(x_{\rho\mu\nu}x_\mu,x_\rho x_\nu,x_{\rho\mu}x_{\mu\nu}
\vert \gamma-\alpha) \ \delta(\gamma + \delta - \beta)$$
\noindent where $h_{\rho,\mu,\nu} = h(x_{\rho\mu}x_{\mu\nu}/x_{\rho\mu\nu}x_\mu)$. The matrix entries of the inverse of $R(\rho,\mu,\nu)$ are given by
$$\bar{R}(\rho,\mu,\nu)_{\gamma,\delta}^{\alpha,\beta} =  
\frac{[\frac{x_{\rho\mu\nu}x_{\mu}}{x_{\rho\mu}x_{\mu\nu}}]}{h_{\rho,\mu,\nu}}\ \zeta^{-\alpha\delta-\frac{\alpha^2}{2}}\ \frac{\delta(\gamma + 
\delta - \beta)}{\omega(\frac{x_{\rho\mu\nu}x_\mu}{\zeta},x_\rho x_\nu,
x_{\rho\mu}x_{\mu\nu}\vert \gamma-\alpha)}\quad .$$
\end{prop}
\noindent Note that the matrices of $6j$-symbols and the normalized Clebsch-Gordan operators have the same form, so that we can write $K_\alpha(\rho,\mu)_{i,j}^k = R(\rho,\mu)_{\alpha,k}^{i,j}$. This explains our choice of the normalization factor $h_{\rho,\mu}$. In fact, one can prove that both are representations of the canonical element of the Heisenberg double of $\mathcal{B}_{\zeta}$, acting on $M_{\rho\mu\nu, \rho\otimes \mu\nu} \otimes M_{\mu\nu,\mu\otimes \nu}$ \cite[\S 3.2-3.3]{B}. This canonical element is called a {\it twisted quantum dilogarithm}.

\medskip

\noindent {\bf Basic pentagon identity and $\Ii$-transits.} We observe that $R(\rho,\mu,\nu)$ is a matrix valued function of $x_\rho x_\nu/x_{\rho\mu}x_{\mu\nu}$ and $x_{\rho\mu\nu}x_\mu/ x_{\rho\mu}x_{\mu\nu}$. Then, let us require that the standard representations $\rho$ used for computing the Clebsch-Gordan operators are defined by taking a {\it same} determination of the $N$-th roots of $t_\rho^{2N}$ and $x_\rho^N$ simultaneously for all $\rho$. The corresponding $6j$-symbols do not depend on the choice of such a determination, because they are homogeneous in the $x$-parameters. Hence, with this convention, we see that $R(\rho,\mu,\nu)$ is a function of, say, $(x_{\rho\mu\nu}x_\mu/x_{\rho\mu}x_{\mu\nu})^N$. Sufficient conditions for the basic pentagon identity (\ref{basicpentagon}) to be true are thus given by the relations between these ratios. We claim that they are just instances of relations between the moduli for the $\Ii$-transit shown in Fig. \ref{idealt}.

\noindent Indeed, associate to the edges $(01)$, $(12)$, $(23)$ and $(34)$ of this figure (for the ordering of the vertices induced, as usual, by the branching) the matrices in (\ref{mat01}) for the representations $\rho$, $\mu$, $\nu$ and $\upsilon$ respectively. Since the sequence $(\rho,\mu,\nu,\upsilon)$ is regular, we can complete this procedure in a unique way on the other edges so that it defines an idealizable Borel valued $1$-cocycle. Now, as explained in Remark \ref{idealrem} (3), the ratios of the form $(x_{\rho\mu\nu}x_\mu/x_{\rho\mu}x_{\mu\nu})^N$ are just the moduli indicated in Fig. \ref{idealt}. So our claim is proved. 

\noindent This discussion shows that the basic pentagon identity holds true when we consider the matrices $R(\rho,\mu,\nu)$ more generally as functions of moduli of idealized hyperbolic tetrahedra, by using the above rule to fix the $N$-th roots of unity. To simplify the notations and also to keep close with those used in \cite{K1} and \cite{B} (where the proofs of the results of this section are given), below we still denote by $R(\rho,\mu,\nu)$ the matrices of $6j$-symbols obtained in Prop. \ref{6j}, which, as we just said, essentially correspond to idealizable Borel valued $1$-cocycles. However, we have to keep in mind the above generalization in terms of moduli.

\medskip

\noindent {\bf Symmetries.} Given a representation $\rho$ of $\mathcal{B}_{\zeta}$, the {\it dual} representation $\bar{\rho}$ is defined by
$$\langle \bar{\rho}(a)\xi,v \rangle = \langle \xi, \rho(S(a))v \rangle$$
where $v \in V_{\rho}$, $\xi \in \bar{V}_{\rho}$ (the dual linear space), $a \in \mathcal{B}_{\zeta}$, $S$ is the antipode of $\mathcal{B}_{\zeta}$, and $\langle .,. \rangle$ is the canonical pairing. In the case where $\rho$ is standard, let us define the {\it inverse} standard representation $\bar{\rho}$ by setting $t_{\bar{p}} = 1/t_p$ and $x_{\bar{p}} = -x_p$. Clearly $\bar{\rho}$ is equivalent to the representation dual to $\rho$ (this explains the abuse of notation). 

\noindent We can rewrite (\ref{symmatrix}) as follows. For any $a$, $c \in \mz/N\mz$ put
$$\begin{array}{l}
R(\rho,\mu,\nu\vert a,c)_{\alpha,\beta}^{\gamma,\delta} = \ \zeta^{c(\gamma - \alpha) - ac/2} \ R(\rho,\mu,\nu)_{\alpha,\beta - a}^{\gamma - 
a ,\delta}\\ 
  \\
\bar{R}(\rho,\mu,\nu\vert a,c)_{\gamma,\delta}^{\alpha,\beta} = \zeta^{c(\gamma - \alpha) + ac/2}\ \bar{R}(\rho,\mu,\nu)_{\gamma+a,\delta}^
{\alpha,\beta+a}\quad .
\end{array}$$
\noindent Note that in (\ref{symmatrix}) we have omitted the index-independent factors $\zeta^{- ac/2}$ and $\zeta^{+ ac/2}$ because of the unavoidable ambiguity of the QHI up to $2N$-th roots of unity (see Remark \ref{phase}). It is easy to verify that 
\begin{eqnarray}\label{factoriz}
R(\rho,\mu,\nu \vert a,c) & = & \zeta^{\frac{ac}{2}} \left( Y_1^{- a}Z_1^{-c} R(\rho,\mu,\nu) Z_1^{c}Z_2^{-a} \right) \nonumber \\ \\ 
\bar{R}(\rho,\mu,\nu \vert a,c) & = & \zeta^{-\frac{ac}{2}} \left( Z_1^{c}Z_2^{-a} \bar{R}(\rho,\mu,\nu) Z_1^{-c}Y_1^{-a} \right) \nonumber
\end{eqnarray}
where $Y_1=Y \otimes id$ etc., and $Y=\zeta^{1/2}XZ$ has components $Y_{m,n} =\omega^{1/2 + n}\delta(m-n-1)$ (the matrices $X$ and $Z$ are defined in (\ref{stmat})). Recall from Section \ref{QDILOGSS} the definition of the matrices $S$ and $T$. Write $\{S^{-1}\}_{m,n}=S^{m,n}$ and so on. Normalizing the scalar factor $\nu$ in $T$ by a certain constant $N$-th root of unity we get:
\begin{prop} \label{symmetry}  Put $b = 1/2 -a -c \in \mathbb{Z}/N\mz$. We have the following {\it symmetry} relations:
$$\begin{array}{l}
\bar{R}(\bar{\rho},\rho\mu,\nu\vert a,b)_{\gamma,\beta}^{\alpha,\delta} = \left(\frac{x_{\rho\mu}x_{\mu \nu}}{x_\mu x_{\rho \mu \nu}}\right)^p\zeta^
{-a/4}\ \sum_{\alpha ',\gamma '=0}^{N-1} R(\rho,\mu,\nu\vert a,c)_{\alpha ',\beta}^
{\gamma ',\delta} \ T_{\gamma,\gamma '} \ T^{\alpha,\alpha'}\\
\\
\bar{R}(\rho\mu,\bar{\mu},\mu\nu\vert b,c)_{\beta,\delta}^{\alpha,\gamma} = \left(\frac{x_{\rho\mu}x_{\mu \nu}}{x_\rho x_{\nu}}\right)^p\zeta^{+c/4}
\ \sum_{\alpha ',\delta '=0}^{N-1} R(\rho,\mu,\nu\vert a,c)_{\alpha ',\beta}^
{\gamma,\delta'} \ T_{\delta,\delta '} \ S^{\alpha,\alpha'}\\
\\
 \bar{R}(\rho,\mu\nu,
\bar{\nu}\vert a,b)_{\alpha,\delta}^{\gamma,\beta} = \left(\frac{x_{\rho\mu}x_{\mu \nu}}{x_\mu x_{\rho \mu \nu}}\right)^p\zeta^{-a/4}\ \sum_{\beta ',\delta '=0}^{N-1} R(\rho,\mu,\nu\vert a,c)_{\alpha,\beta '}^{\gamma,
\delta'} \ S_{\delta,\delta '} \ S^{\beta,\beta '}\quad .
\end{array}$$
\end{prop}
\noindent Note that, for instance, the factor $(x_{\rho\mu}x_{\mu \nu}/x_\mu x_{\rho \mu \nu})^p$ in the first identity is written as $(w_0')^{-p}$ with the notations of Lemma \ref{6jsym}. 

\noindent Given a standard representation $\rho$ define the complex conjugate representation $\rho^*$ by $t_{\rho^*}=(t_\rho)^*$ and $x_{\rho^*}=(x_\rho)^*$.
\begin{prop} \label{unitarite} We have the following unitarity property:
$$\bar{R}(\rho^*,\mu^*,\nu^*\vert a,c)_{\gamma,\delta}^{\alpha,\beta} = 
\left( R(\rho,\mu,\nu\vert a,c)_{-\alpha,-\beta}^{-\gamma,-\delta} \right)^*\quad .
\nonumber$$
\end{prop}
\noindent {\bf Partially symmetrized basic pentagon identity.} Let us use the notations of the proof of Lemma \ref{charge-transit}. Consider the following set of independent charges: $i = c_{01}^4$, $j = c_{01}^2$, $k = c_{12}^0$, $l = c_{23}^1$ and $m = c_{12}^3$. They determine completely the charge transit shown in Fig. \ref{w(e)}. We can easily show that $l + m = c_{13}^2$, $l-i =  c_{23}^0$, $j+k= c_{02}^1$, $i+j =  c_{01}^3$ and $m-k =  c_{12}^4$. Note that the branching in Fig. \ref{w(e)} is the same as the one of Fig. \ref{idealt}. Moreover, we have seen above that the $\Ii$-transit of Fig. \ref{idealt} dominates the basic pentagon identity. The following proposition describes a `charged' generalization of this identity:
\begin{prop} \label{EP} We have:
\begin{eqnarray}
R_{12}(\rho,\mu,\nu \vert i,m-k)\ R_{13}(\rho,\mu\nu,\upsilon \vert j,l+m)\ R_{23}(\mu,\nu,\upsilon 
\vert k,l-i) =  \hspace{3cm} \nonumber \\ 
\hspace{3cm} R_{23}(\rho\mu,\nu,\upsilon \vert j+k,l)\ R_{12}(\rho,\mu,\nu\upsilon 
\vert i+j,m).\nonumber \hspace{2mm}
\end{eqnarray}
\end{prop}

\noindent The proof consists in using the formulas (\ref{factoriz}) and the commutation relations between the matrices $Y$, $Z$ and $R(\rho,\mu,\nu)$ to reduce the statement to the basic pentagon identity.

\bigskip

\noindent {\bf Acknowledgements.} The first author is grateful to the Department of Mathematics of the University of Pisa for its kind hospitality in September 2000 and March 2001, while a part of the present work has been done. It was completed in Pisa when he was a fellow of the EDGE network of the European Community, and later a Marie Curie fellow (proposal No MCFI-2001-00605). The second author profited of the kind invitation of the Laboratoire E. Picard, University P. Sabatier of Toulouse, in December 2000 and June 2001.


\begin{thebibliography}{99}

\bibitem{B} S. Baseilhac.  \emph{Dilogarithme quantique et invariants
de $3$-vari\'{e}t\'{e}s}, Th\`{e}se, Universit\'{e} Paul Sabatier,
Toulouse (France), october 2001;

\bibitem{BB1} S. Baseilhac, R. Benedetti. \emph{QHI, 3-manifold
scissors congruence classes and the volume conjecture}, Geometry and
Topology Monographs, Volume 4: Invariants of Knots and 3-manifolds
(Kyoto 2001), 13-28;

\bibitem{BB2} S. Baseilhac, R. Benedetti.  \emph{Quantum and classical
dilogarithmic invariants of flat $PSL(2,\mc)$-bundles over
$3$-manifolds}, preprint available on arXiv:math.GT/ (2003);

\bibitem{BP1} R. Benedetti and C. Petronio. \emph{Lectures on
Hyperbolic Geometry}, Springer (1992);

\bibitem{BP5} R. Benedetti and C. Petronio.  \emph{A finite graphic
calculus for 3-manifolds}, Manuscripta Math. 88 (1995), 291-310;

\bibitem{BP2} R. Benedetti, C. Petronio. \emph{Branched Standard
Spines of $3$-Manifolds}, Lect. Notes in Math. No 1653, Springer
(1997);

\bibitem{BP3} R. Benedetti and C. Petronio. 
\emph{Combed $3$-manifolds with concave boundary, 
framed links, and pseudo-legendrian links}, J. Knot Th. \&
Ram. Vol. 10, No 1 (2001), 1-35;

\bibitem{BP4} R. Benedetti and C. Petronio. \emph{ Reidemeister
torsion of 3-dimensional Euler structures with simple boundary
tangency and pseudo-legendrian knots,} Manuscripta Math. 106 (2001), 13-61;

\bibitem{BR} V. V. Bazhanov, N. Yu. Reshetikhin. \emph{Remarks on the quantum dilogarithm}, 
J. Phys. A: Math. Gen. No 28 (1995), 2217-2226;

\bibitem{Cas} B. G. Casler. \emph{\it An embedding theorem for connected $3$-manifolds with boundary}, 
Proc. Amer. Math. Soc. No 16 (1965), 559-566;

\bibitem{CP} V. Chari, A. Pressley. {\it A Guide To Quantum Groups}, Cambridge University Press (1994);

\bibitem{D1} J.L. Dupont. \emph{The dilogarithm as a characteristic
class for flat bundles}, J. Pure App. Alg. 44 (1987), 137-164;

\bibitem{D2} J.L. Dupont, \emph{Scissors Congruences, Group Homology and Characteristic Classes}, Nankai Tracts in Math. Vol.1, World Scientific (2001).

\bibitem{DCP} C. De Concini, C. Procesi. \emph{Quantum Groups}, in ``$\Dd$-modules, 
Representation Theory and Quantum Groups'', LNM 1565, Springer (1993);

\bibitem{DS} J. L. Dupont, C-H. Sah. \emph{Scissors congruences II}, J. Pure App. Alg. 44 
(1987), 137-164;

\bibitem{EP} D. B. A. Epstein, R. Penner. \emph{Euclidian
decompositions of non-compact hyperbolic manifolds}, J. Diff. Geom. 27
(1988), 67-80;

\bibitem{FK} L. D. Faddeev, R.M. Kashaev. \emph{Quantum dilogarithm}, Mod. Phys. 
Lett. A Vol. 9, No 5 (1994), 427-434;

\bibitem{K1} R. M. Kashaev. \emph{Quantum dilogarithm as a
6j-symbol}, Mod. Phys. Lett. A Vol. 9, No 40 (1994), 3757-3768;

\bibitem{K2} R. M. Kashaev. \emph{A link invariant from quantum
dilogarithm}, Mod. Phys. Lett. A Vol. 10 (1995), 1409-1418;

\bibitem{K3} R. M. Kashaev. \emph{The algebraic nature of quantum dilogarithm}, Geometry of 
Integrable models (Dubna 1994), World Scientific, River Edge, NJ (1996), 32-51;

\bibitem{K4} R. M. Kashaev. \emph{The hyperbolic volume of knots
from the quantum dilogarithm}, Lett. Math. Phys. 39 (1997), 269-275;

\bibitem{K5} R. M. Kashaev. \emph{Quantum hyperbolic invariants of
knots}, in \emph{Discrete Integrable Geometry and Physics}, A. I. Bobenko \&
R. Seiler ed., Oxford Science Publications (1998);

\bibitem{Mak} A. Yu. Makovetskii. \emph{On transformations of
special spines and special polyhedra}, Math. Notes, Vol. 65, No 3,
1999;

\bibitem{Mat} S. V. Matveev. \emph{Transformations of special spines and the Zeeman 
conjecture}, Math. USSR Izvestia No 31 (1988), 423-434;

\bibitem{MM} H. Murakami, J. Murakami. \emph{The colored Jones
polynomials and the simplicial volume of a knot}, Acta Math. 186 (2001), No. 1, 85-104;

\bibitem{N1} W. D. Neumann. \emph{Combinatorics of triangulations
and the Chern-Simons invariant for hyperbolic $3$-manifolds}, in
\emph{Topology'90} (Columbus, OH, 1990), De Gruyter, Berlin (1992);

\bibitem{N2} W. D. Neumann. {\it Hilbert's 3rd problem and
invariants of 3-manifolds}, Geometry and Topology Monographs 1 (1998), 
The Epstein Birthday Schrift, paper No 19, 383-411;

\bibitem{N3} W. D. Neumann. {\it Extended Bloch group and the
Chern-Simons class (Incomplete working version)}, arXiv:math.GT/0212147;

\bibitem{NY} W.D. Neumann, J. Yang. {\it Bloch invariants of
hyperbolic $3$-manifolds}, Duke Math. Journal Vol. 96, No 1 (1999);

\bibitem{Pi} R. Piergallini. \emph{Standard moves for standard polyhedra and spines}, 
Rend. Circ. Mat. Palermo No 37, suppl. 18 (1988), 391-414;

\bibitem{T} V. Turaev. \emph{Euler structures, nonsingular vector fields, and torsions of Reidemeister type}, Izvestia Ac. Sci. USSR 53 (1989); translated in Math. USSR Izvestia 34 (1990), 627-662;

\bibitem{TV} V. G. Turaev, O. Viro. \emph{State sum invariants of
3-manifolds and 6j-symbols}, Topology 31 (1992), 865-904;

\bibitem{Y} Y. Yokota. {\it On the volume conjecture for hyperbolic knots}, preprint available at http://www.comp.metro-u.ac.jp/$\sim$jojo/volume-conjecture.pdf;  

\bibitem{W} E. Witten. \emph{2+1 Dimensional gravity as an exact soluble system}, 
Nuclear Physics B 323 (1989), 113 - 140;

 
\end{thebibliography}
\end{document}